\documentclass[a4paper,11pt,reqno]{amsart}

\usepackage{graphicx}
\usepackage{amsmath}
\usepackage{amssymb, amsthm}
\usepackage{amsfonts}
\usepackage{comment}
\usepackage{mathrsfs}
\usepackage{graphics}
\usepackage{float}
\usepackage{mathtools}
\usepackage{enumitem}

\usepackage[utf8]{inputenc}
\usepackage[T1]{fontenc}
\usepackage[lowtilde]{url}

\DeclareFixedFont{\ttb}{T1}{txtt}{bx}{n}{8} %%{12} % for bold
\DeclareFixedFont{\ttm}{T1}{txtt}{m}{n}{8}  %% {12}  % for normal

\newcommand\pprof[1]{{\it Proof of {#1}.}}

\graphicspath{{eps_L/}}

 \DeclarePairedDelimiter\norm{\lVert}{\rVert}%

\makeatletter

\newtheorem{theorem}{Theorem}[section]

\@addtoreset{equation}{section}

\newtheorem{remark}[theorem]{Remark}
\newtheorem{proposition}[theorem]{Proposition}
\newtheorem{lemma}[theorem]{Lemma}
\newtheorem{example}[theorem]{Example}

\newtheorem{definition}[theorem]{Definition}

\def\PerfProof{{\it Proof.\ }}

\begin{document}

\title[Decomposition of the longest elements of the Weyl group]
{Decomposition of the longest element of the Weyl group using factors corresponding to the highest roots}
         \author{Rafael Stekolshchik}
\date{}

\begin{abstract}
 Let $\varPhi$ be a root system of a finite Weyl group $W$ with simple roots $\Delta$
and corresponding simple reflections $S$. For $J \subseteq S$,
denote by $W_J$  the standard parabolic subgroup of
$W$ generated by $J$, and by $\Delta_J \subseteq \Delta$ the subset corresponding to $J$.
We show that the longest element of $W$ is decomposed into a product of several ($\le |\Delta|$)
reflections corresponding to mutually orthogonal roots, each of which
is either the highest root of some subset $\Delta_J \subseteq \Delta$ or
is a simple root. For each type of the root system, the factors of the specified decomposition
are listed. The relationship between the longest elements of different types
is found out. The uniqueness of the considered decomposition is shown.
It turns out that subsets of highest roots, which give the decomposition
of longest elements in the Weyl group, coincide with the cascade of orthogonal roots
constructed by B.Kostant and A.Joseph for calculations in the universal enveloping algebra.
\end{abstract}

\maketitle

%% \tableofcontents

\section{\bf Introduction}

\subsection{The longest element}

Let $\varPhi$ be a root system of a finite Weyl group $W$
with simple roots $\Delta = \{\alpha_1, \dots, \alpha_n \}$  and the corresponding simple
reflections $S = \{s_{\alpha_1}, \dots, s_{\alpha_n} \}$, $\varPhi^+$
be the subset of positive roots in $\varPhi$.

For $w \in W$, define $l(w)$ (resp., $l_a(w)$) to be the minimal
number of factors occurring amongst all expressions of $w$ as a product of simple
reflections $S$ (resp., reflections). The function $l$ (resp., $l_a$) is called the standard
(resp., absolute) length function of $(W, S)$. An expression $w = s_1\dots s_n$ with $s_i \in S$
and $n = l(w)$ is called a {\it reduced expression} for $w$.

There exists an element $w_0 \in W$ sending the subset of positive
roots $\varPhi^+$ to the subset of negative roots $\varPhi^-$.
Such an element $w_0$ is unique in $W$.
The length $l(w_0)$ coincides with the number of roots in $\varPhi^+$.
No other element of $W$ has such a large length as $w_0$.  So, the element $w_0$ is
said to be the {\it longest element}. The element $w_0$ transforms the fundamental chamber $C$ to
the chamber $-C$,  $w_0$ is involution and
\begin{equation*}
    l(w_0 w) = l(w_0) - l(w) \text{ for all } w \in W.
\end{equation*}
The element $w_0$ is the unique element $w \in W$ satisfying the condition
\begin{equation*}
    l(w s_{\alpha}) < l(w) \text{ for all } \alpha \in \Delta.
\end{equation*}
The longest element $w$ is $-1$ except for the following $3$ cases\footnotemark[1]:
\footnotetext[1]{In all these cases $\varepsilon$ is the involutive automorphism of the corresponding Dynkin diagram.
For the numbering of vertices, see Fig. \ref{fig_ABCD}.}
\begin{equation*}
\footnotesize
\begin{array}{lll}
  & (1) \; A_n, \; n \geq 2, &  w_0 = -\varepsilon, \text{ where } \;
            \varepsilon(\alpha_i) = \alpha_{n - i + 1},    \\
  & (2) \; D_n, \; n \text{ is odd},  &  w_0 = -\varepsilon, \text{ where }
          \varepsilon(\alpha_{n-1} ) = \alpha_{n}, \; \varepsilon(\alpha_{n-1}) = \alpha_n, \; \\
  & & \varepsilon(\alpha_{i}) = \alpha_i \text{ for other } \alpha_i. \\
  & (3) \; E_6,  & w_0 = -\varepsilon, \text{ where }
     \varepsilon(\alpha_1) = \alpha_6, \; \varepsilon(\alpha_6) = \alpha_1, \\
  &  & \varepsilon(\alpha_3) = \alpha_5, \; \varepsilon(\alpha_5) = \alpha_3, \;
      \varepsilon(\alpha_2) = \alpha_2, \; \varepsilon(\alpha_4) = \alpha_4, \\
\end{array}
\end{equation*}
see \cite[Ch.VI,$\S$1,$n^{\rm o}$6,Cor.6]{Bo02}, \;
   \cite[Plates I-X]{Bo02}, \; \cite[$\S1.8$]{Hu90}, \cite[$\S$2.3]{BB05}.
~\\

\subsection{The main results}
\subsubsection{Decomposition of the longest element}

For any $J \subseteq S$, denote by $W_J$ the subgroup of $W$ generated by $J$,
and by $\Delta_J \subseteq \Delta$ the subset corresponding to $J$.
The subgroups of the form $W_J$ are referred as {\it standard parabolic subgroups}.

\begin{theorem}
  \label{th_factoriz}
   The longest element $w_0 \in W$ is decomposed into a product of several ($\leq n$)
   reflections corresponding to mutually orthogonal roots, each of which is either the highest root of
   some subset $\Delta_J \subseteq S$ or is a simple root,
   see Tables \ref{tab_factoriz_all} and \ref{tab_max_orhog_set}.
\end{theorem}

  \PerfProof The theorem is proved by a case-by-case analysis: see
  Propositions \ref{prop_factoriz_An_2} ($A_n$), \ref{prop_factoriz_Bn} ($B_n$),
  \ref{prop_factoriz_Cn} ($C_n$),
  \ref{prop_1_Dn} ($D_n$), \ref{prop_factoriz_E6} ($E_6$), \ref{prop_factoriz_E7} ($E_7$),
  \ref{prop_factoriz_E8} ($E_8$), \ref{prop_factoriz_F4} ($F_4$), \ref{prop_factoriz_G2} ($G_2$).
   \qed

\subsubsection{Relationship between the longest elements}

Let $W$ be a finite Weyl group and $\varPhi$ be the corresponding root system,
in this case we will write $W = W(\varPhi)$. Denote by $w_0(W)$ the longest element in $W$.
We also use the notation $w_0(\varPhi)$ instead of $w_0(W)$.

\begin{theorem}
 \label{th_connection}
For any root system $\varPhi$ with the Weyl group $W = W(\varPhi)$, there exists
a root subsystem $\varPhi' \subset \varPhi$ with the standard parabolic subgroup
$W' = W(\varPhi')$ in $W$, such that the longest elements $w_0(\varPhi)$ and $w_0(\varPhi')$ are related
as follows:
\begin{equation}
 \label{eq_reccur_rel}
  \begin{array}{lll}
     & w_0(A_n) = w_0(A_{n-2})s_{\alpha_{max}} & \text{ for } n \ge 3, \\
     & w_0(B_n) = w_0(B_{n-2})s_{\alpha_{max}}s_{\alpha_1} & \text{ for } n \ge 4, \\
     & w_0(C_n) = w_0(C_{n-1})s_{\alpha_{max}} & \text{ for } n \ge 3, \\
     & w_0(D_n) = w_0(D_{n-2})s_{\alpha_{max}}s_{\alpha_1} & \text{ for } n \ge 6, \\
     & w_0(E_6) = w_0(A_5)s_{\alpha_{max}}, \\
     & w_0(E_7) = w_0(D_6)s_{\alpha_{max}}, \\
     & w_0(E_8) = w_0(E_7)s_{\alpha_{max}}, \\
     & w_0(F_4) = w_0(C_3)s_{\alpha_{max}}. \\
  \end{array}
\end{equation}
In eq. \eqref{eq_reccur_rel}, the reflection
$s_{\alpha_{max}}$ corresponds to the highest root $\alpha_{max}$ in the root system $\varPhi$.
\end{theorem}

\PerfProof Let us consider, for example, the case $W = W(F_4)$, $W' = W(C_3)$.
\begin{figure}[h]
\centering
   \includegraphics[scale=0.4]{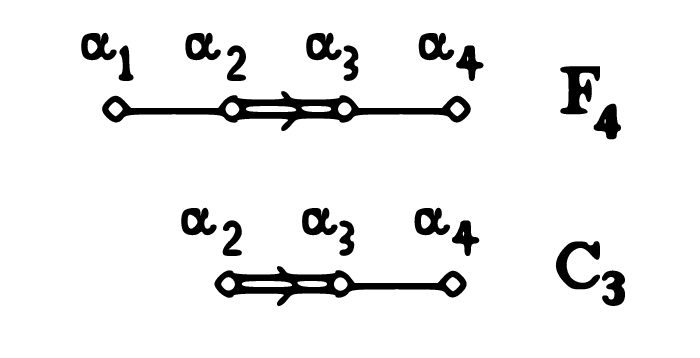}
\caption{Numbering of $C_3$ vertices which is different from Bourbaki's numbering}
%%%%%% The label must come after caption
\label{fig_F4C3}
\end{figure}
By Proposition \ref{prop_factoriz_F4}
\begin{equation*}
  w_0(F_4) = s_{\alpha_2} s_{\alpha_2 + 2\alpha_3} s_{\alpha_2 + 2\alpha_3 + 2\alpha_4}s_{\alpha_{max}}.
\end{equation*}
The roots $\alpha_2 + 2\alpha_3$ and $\alpha_2 + 2\alpha_3 + 2\alpha_4$ are
the highest roots for the root systems $C_2$ and $C_3$, see Fig. \ref{fig_F4C3}.
Denote them as follows:
\begin{equation*}
     \alpha^{c2}_{max} := \alpha_2 + 2\alpha_3, \quad \alpha^{c3}_{max} := \alpha_2 + 2\alpha_3 + 2\alpha_4,
\end{equation*}
although they differ form definitions of $\alpha^{c2}_{max}$ and $\alpha^{c3}_{max}$ for case $C_3$
in Table \ref{tab_factoriz_all}. So,
\begin{equation*}
  w_0(F_4) = s_{\alpha_2} s_{\alpha^{c2}_{max}} s_{\alpha^{c3}_{max}}s_{\alpha_{max}} = w_0(C_3)s_{\alpha_{max}}.
\end{equation*}
Such a difference in definitions of the highest roots also takes place in other cases, see Remark \ref{rem_differ_indices}.
The definition of the highest roots of root subsystems depends on the root system in which this root is considered.

Other cases in eq. \eqref{eq_reccur_rel} are treated in a similar way. \qed

\subsubsection{Uniqueness of decomposition}
  \label{sec_uniqueness}

\begin{definition}{\rm
Let $T = \{\tau_1, \tau_2, \dots, \tau_m \}$ be the subset of distinct roots in $\varPhi$.
This set is said to be {\it max-orthogonal} if
\begin{enumerate}[label=(\alph*)]
  \item the roots of $T$ are {\it mutually orthogonal},
  \item all non-simple roots in $T$ form a {\it linearly ordered subset} in $\varPhi$.
  \item each non-simple root $\tau_i$ is the {\it highest root} in
  some root subsystem $J \subset \varPhi$ corresponding to a standard parabolic
  subgroup  $W_J \subset W$. The root $\tau_1$ is the highest root in $\varPhi$, i.e.,
  $\tau_1 = s_{\alpha_{max}}$.
\end{enumerate}
}
\end{definition}

Note that each root subset used in decompositions
of Theorem \ref{th_factoriz} is max-orthogonal.

\begin{theorem}
 \label{th_uniqueness}
  For any longest element $w_0$, there exists a unique max-orthogonal
  subset $\{\tau_1, \tau_2, \dots, \tau_m \}$, where $m \leq n$, such that
\begin{equation}
  \label{eq_unique}
    w_0 = \prod\limits^m_{i=1}s_{\tau_i}.
\end{equation}
\end{theorem}

\begin{definition}
The decomposition \eqref{eq_unique} corresponding to some max-orthogonal subset $\{\tau_1, \tau_2, \dots, \tau_m \}$
is said to be the {\it max-orthogonal decomposition}.
\end{definition}
We will construct the max-orthogonal decomposition for
each type of root systems, see Tables \ref{tab_factoriz_all} and Table \ref{tab_max_orhog_set}.
~\\

\pprof{Theorem \ref{th_uniqueness}}
In  eq. \eqref{eq_reccur_rel},  for each relation,  there exists a factor
that connects two different longest elements. We refer to this factor as the {\it linking factor}.
Denote it by $\mathcal{L}$. There are two cases for the linking factor $\mathcal{L}$:
\begin{enumerate}[label=(\alph*)]
 \item  $\mathcal{L} = s_{\alpha_{max}}$, this holds for $A_n$, $C_n$, $E_n$ and $F_4$,
 \item $\mathcal{L} = s_{\alpha_{max}}s_{\alpha_1}$, this holds for $B_n$ and $D_n$.
\end{enumerate}
 The uniqueness of the max-orthogonal decomposition \eqref{eq_unique} will be proved
 by induction on the length of the longest element.

 (a) Let $w_0(\varPhi) = w_0(\varPhi')s_{\alpha_{max}}$. By definition of the max-orthogonal decomposition,
 the decomposition \eqref{eq_unique} of $w_0(\varPhi)$  contains the factor $s_{\alpha_{max}}$.
 By induction hypothesis $w_0(\varPhi')$ has an unique max-orthogonal decomposition.
 Then, the max-orthogonal decomposition of $w_0(\varPhi)$ is unique.

\begin{figure}[h]
\centering
   \includegraphics[scale=0.35]{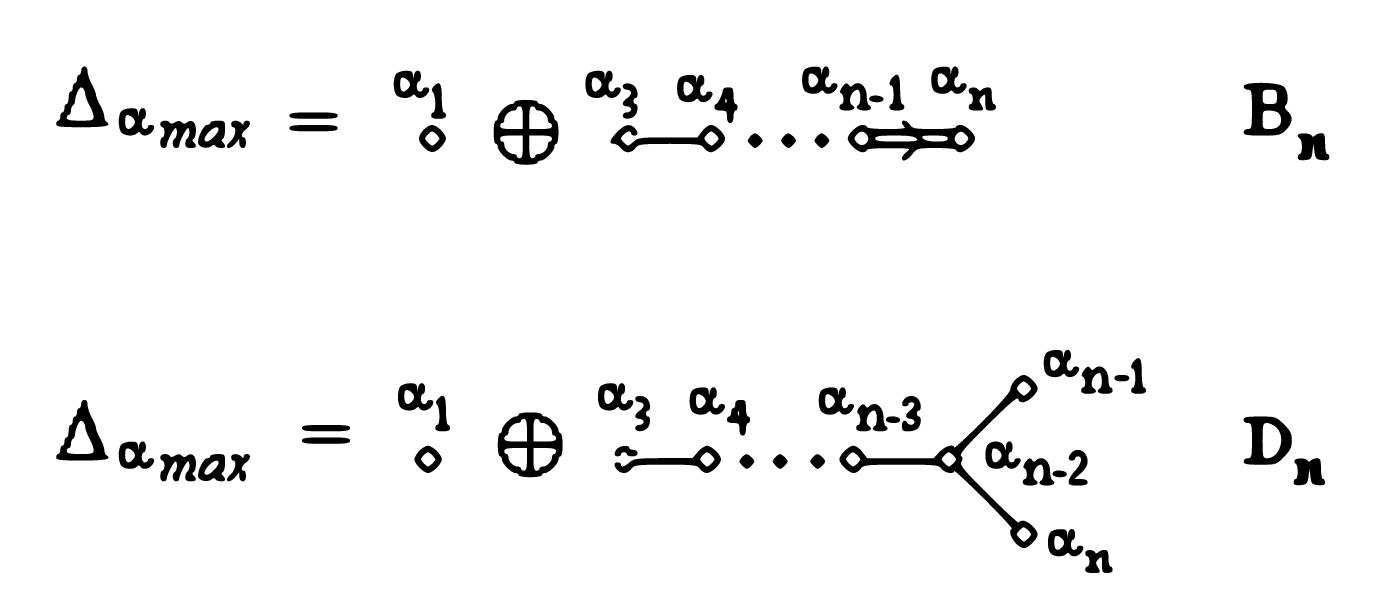}
\caption{Subsets $\Delta_{\alpha_{max}}$ for $B_n$ and $D_n$}
%%%%%% The label must come after caption
\label{fig_Bn_Dn}
\end{figure}

 (b) Let $w_0(\varPhi) = w_0(\varPhi')s_{\alpha_{max}}s_{\alpha_1}$. In the root systems $B_n$ and $D_n$,
 the only simple root non-orthogonal to the highest root $\alpha_{max}$ is $\alpha_2$.
 Consider the root subset $\Delta_{\alpha_{max}}$ consisting
 of roots orthogonal to $\alpha_{max}$. The subset $\Delta_{\alpha_{max}}$ consists of two mutually orthogonal
 subsets: $\{\alpha_1\}$ and $V(\alpha_3, \dots, \alpha_n)$, the subset of roots spanned by
 the roots $\{\alpha_3, \dots, \alpha_n \}$:
 \begin{equation}
    \Delta_{\alpha_{max}} = \{\alpha_i \in \varPhi \mid (\alpha_i, \alpha_{max}) = 0 \} =
       \{\alpha_1\} \oplus V(\alpha_3, \dots, \alpha_n),
 \end{equation}
 see Fig. \ref{fig_Bn_Dn}.
 Since $s_{\alpha_1}$ is the factor in the max-orthogonal decomposition of $w_0(\varPhi)$
 given by \eqref{eq_reccur_rel}
 and $\alpha_1$ is disconnected from $V(\alpha_3, \dots, \alpha_n)$,
 then $s_{\alpha_1}$ is the factor in any other max-orthogonal decomposition of $w_0(\varPhi)$.
 Thus, any max-orthogonal decomposition of $w_0(\varPhi)$ contains $s_{\alpha_{max}}s_{\alpha_1}$.
 Further, we apply induction as in case (a). \qed
~\\

The proof of uniqueness in Theorem \ref{th_uniqueness}  echoes
with the cascade of orthogonal roots constructed  by B.Kostant and A.Joseph.
It turns out that subsets of highest roots in the max-orthogonal decomposition
coincide with the cascade of orthogonal roots
constructed by B.Kostant and A.Joseph for calculations in the universal enveloping algebra,
see $\S$\ref{sec_cascade}.
~\\

\begin{figure}[h]
\centering
   \includegraphics[scale=0.30]{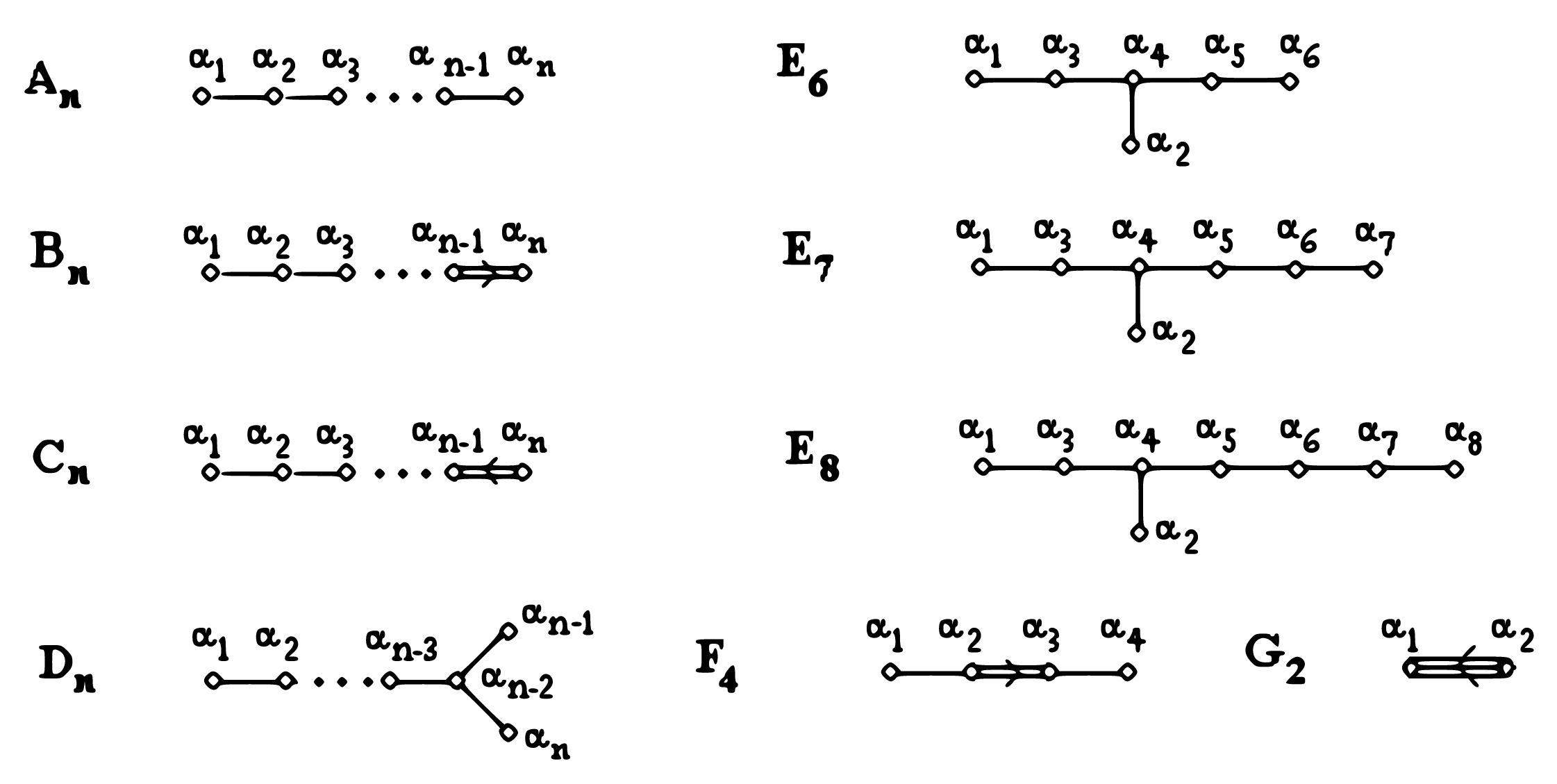}
\caption{Bourbaki's numbering of simple roots, \cite{Bo02}}
%%%%%% The label must come after caption
\label{fig_ABCD}
\end{figure}

\begin{remark}[notations]{\rm
  \label{rem_notations}

(i) We follow Bourbaki's numbering for simple roots in Dynkin diagrams.
see Fig. \ref{fig_ABCD}.

(ii) Table \ref{tab_factoriz_all} contains definitions of highest roots $\alpha^{xi}_{max}$ (= $\alpha^{x,i}_{max}$),
where $x$ is one of the indices $a,b,c,d,e,f,g$, and $i$ is the number of vertices in
the corresponding Dynkin diagram.

(iii) In Tables \ref{tab_highest_roots_ABCD} and \ref{tab_highest_roots_EFG}, for each
root system $\varPhi$, we list the highest roots of the root subsystems used
in Table \ref{tab_factoriz_all}.

(iv) The reflection corresponding to $\alpha_{max}$ is denoted by $s_{\alpha_{max}}$.

(v) Sometimes we prefer to use the notation $s_i$ instead of $s_{\alpha_i}$, which is the same.
}
\end{remark}

\begin{remark}[summands in $\alpha^{xi}_{max}$]{\rm
 \label{rem_differ_indices}
  The indices of simple roots that appears as summands in
  $\alpha^{a5}_{max}$ for the case $E_6$ are different from the indices
  contained in $\alpha^{a5}_{max}$ for the case $A_n$, see Fig. \ref{fig_A7_n_E6}.

\begin{figure}[h]
\centering
   \includegraphics[scale=0.25]{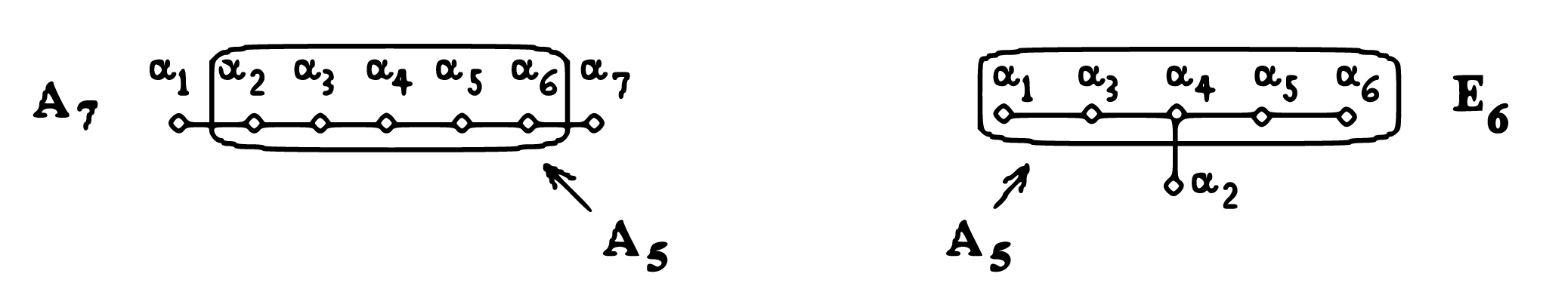}
\caption{\footnotesize Root systems $A_5 \subset A_7$ and $A_5 \subset E_6$. }
%%%%%% The label must come after caption
\label{fig_A7_n_E6}
\end{figure}

    The indices of simple roots that appears as summands in
  $\alpha^{d4}_{max}$ (resp. $\alpha^{d6}_{max}$) for the case $E_7$ are different from the indices
  contained in $\alpha^{d4}_{max}$ (resp. $\alpha^{d6}_{max}$) for
  the case $D_n$, see Fig. \ref{fig_A7_n_E6}.

\begin{figure}[h]
\centering
   \includegraphics[scale=0.25]{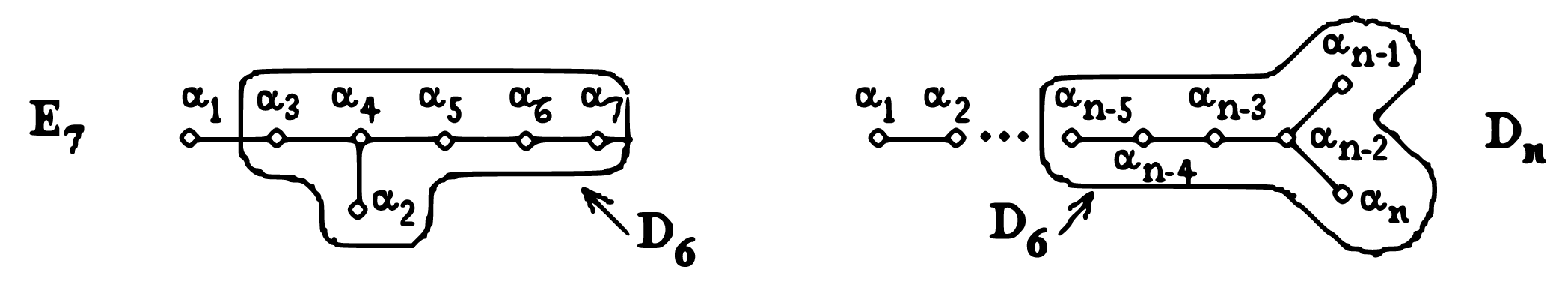}
\caption{\footnotesize Root systems $D_6 \subset E_7$ and $D_6 \subset D_n$. }
%%%%%% The label must come after caption
\label{fig_E7_n_Dn}
\end{figure}
}
\end{remark}

\subsection{Acknowledgments}
I thank Valdemar Tsanov for pointing out on papers \cite{DT16}, \cite{Kos12}, as well as
the connection between {\it max-orthogonal subsets} from Theorem \ref{th_uniqueness} and
{\it Kostant's cascade} construction, see $\S$\ref{sec_cascade}.

\begin{table}[H]
\centering
\renewcommand{\arraystretch}{1.3}
\begin{tabular}{|c|c|c|}
  \hline
       \footnotesize {Weyl}  & \footnotesize{Longest element $w_0$ } & $l(w_0)$  \\
       \footnotesize {group}  & &   \\
  \hline
     \footnotesize $W(A_n)$
     & \footnotesize \qquad $\begin{array}{l}
        \prod\limits_{i=1}^k s_{\alpha^{a,2i}_{max}}, \;  \text{ where } \;
        \alpha^{a,2i}_{max} := \sum\limits_{j = k-i+1}^{k+i} \alpha_j, \;\;
                n = 2k. \\
        \\
        s_{\alpha_{k+1}}\prod\limits_{i=1}^k s_{\alpha^{a,2i+1}_{max}}, \; \text{ where } \;
      \alpha^{a,2i+1}_{max} := \sum\limits_{j = k-i+1}^{k+i+1} \alpha_j, \;\;
                 n = 2k+1, \\
     \end{array}$ & $\frac{n(n+1)}{2}$ \\
  \hline
    \footnotesize $W(B_n)$
      &  \footnotesize \quad $\begin{array}{ll}
  & \qquad \prod\limits_{\stackrel{i=1}{i-odd}}^{n-1}(s_{\alpha_i}s_{\alpha^{b,n-i+1}_{max}})
    \quad \text{ for } n = 2k, \\
    \\
     & s_{\alpha_n}
    \prod\limits_{\stackrel{i=1}{i-odd}}^{n-2}(s_{\alpha_i}s_{\alpha^{b,n-i+1}_{max}})
    \quad \text{ for } n = 2k+1, \\
    & \text{ where } \;
    \alpha^{bi}_{max} := \alpha_{n-i+1} + 2\sum\limits_{j = n-i+2}^n \alpha_j.
     \end{array}$ & \footnotesize $n^2$ \\
  \hline
    \footnotesize $W(C_n)$ &   \footnotesize
      $s_{\alpha_n} \prod\limits_{i=1}^{n-1}s_{\alpha^{c,n-i+1}_{max}}$, \;
  where \; $\alpha^{c,n-i+1}_{max} := \alpha_n + 2\sum\limits_{j = i}^{n-1} \alpha_j$.
      & \footnotesize $n^2$  \\
  \hline
    \footnotesize $W(D_n)$ &   \footnotesize $\begin{array}{cc}
      & s_{\alpha_n}s_{\alpha_{n-1}}
        \prod\limits_{\stackrel{i=1}{i-odd}}^{n-3}(s_{\alpha_i}s_{\alpha^{d,n-i+1}_{max}})
        \quad \text{ for }n = 2k,\\
        \\
     & \qquad \qquad \prod\limits_{\stackrel{i=1}{i-odd}}^{n-2}(s_{\alpha_i}s_{\alpha^{d,n-i+1}_{max}})
        \quad \text{ for }n = 2k+1, \\
     & \text{ where } \;
       \alpha^{di}_{max} := \alpha_{n-i+1} + 2\sum\limits_{j = n-i+2}^{n-2}\alpha_j + \alpha_{n-1} + \alpha_n \text{ for } i \ge 4, \\
     & \alpha^{d3}_{max} := \alpha^{a3}_{max} := \alpha_{n-2} + \alpha_{n-1} + \alpha_n.

           \end{array}$ & \footnotesize $n(n-1)$ \\
  \hline
     \footnotesize  $W(E_6)$ & \footnotesize $\begin{array}{cc}
        & s_{\alpha_4} s_{\alpha^{a3}_{max}}
                       s_{\alpha^{a5}_{max}}s_{\alpha_{max}}, \text{ where} \\ %% \footnotemark[1]
        & \alpha^{a3}_{max} := \alpha_3 + \alpha_4 + \alpha_5, \;\;
        \alpha^{a5}_{max} := \alpha_1 + \alpha_3 + \alpha_4 + \alpha_5 + \alpha_6,     \\
    &  \alpha_{max} = \alpha_1 + 2\alpha_2 + 2\alpha_3 + 3\alpha_4 + 2\alpha_5 + \alpha_6.
       \end{array}$ & \footnotesize $36$  \\
  \hline
     \footnotesize $W(E_7)$ & \footnotesize $\begin{array}{cc}
     & s_{\alpha_2}s_{\alpha_3}s_{\alpha_5} s_{\alpha_7}
        s_{\alpha^{d4}_{max}} s_{\alpha^{d6}_{max}}
        s_{\alpha_{max}},\\
       & \text{ where }
       \alpha^{d4}_{max} := \alpha_2 + \alpha_3 + 2\alpha_4 + \alpha_5,  \\
     & \alpha^{d6}_{max} :=
        \alpha_2 + \alpha_3 + 2\alpha_4 + 2\alpha_5 + 2\alpha_6 + \alpha_7,  \\
     & \alpha_{max} = 2\alpha_1 + 2\alpha_2 + 3\alpha_3 + 4\alpha_4 +
       3\alpha_5 + 2\alpha_6 + \alpha_7.
       \end{array}$ & \footnotesize $63$ \\
  \hline
     \footnotesize $W(E_8)$ & \footnotesize $\begin{array}{cc}
     & s_{\alpha_2}s_{\alpha_3}s_{\alpha_5} s_{\alpha_7}
        s_{\alpha^{d4}_{max}} s_{\alpha^{d6}_{max}} s_{\alpha^{e7}_{max}}
        s_{\alpha_{max}},\\
       & \text{ where }
       \alpha^{d4}_{max} := \alpha_2 + \alpha_3 + 2\alpha_4 + \alpha_5,  \\
     & \alpha^{d6}_{max} :=
        \alpha_2 + \alpha_3 + 2\alpha_4 + 2\alpha_5 + 2\alpha_6 + \alpha_7,  \\
     & \alpha^{e7}_{max} := 2\alpha_1 + 2\alpha_2 + 3\alpha_3 + 4\alpha_4 +
       3\alpha_5 + 2\alpha_6 + \alpha_7,  \\
     & \alpha_{max} = 2\alpha_1 + 3\alpha_2 + 4\alpha_3 + 6\alpha_4 +
       5\alpha_5 + 4\alpha_6 + 3\alpha_7 + 2\alpha_8.
       \end{array}$ & \footnotesize $120$ \\
  \hline
    \footnotesize $W(F_4)$ & \footnotesize $\begin{array}{cc}
        & s_{\alpha_2}s_{\alpha_2 + 2\alpha_3}
             s_{\alpha_2 + 2\alpha_3 + 2\alpha_4} s_{\alpha_{max}},  \\
      & \text{ where } \alpha_{max} = 2\alpha_1 + 3\alpha_2 + 4\alpha_3 + 2\alpha_4.
        \end{array}$&  \footnotesize $24$  \\
  \hline
    \footnotesize $W(G_2)$ &  \footnotesize
           $s_{\alpha_1 + \alpha_2} s_{\alpha_2 + 3\alpha_1}$ & \footnotesize $6$ \\
  \hline
\end{tabular}
\vspace {1mm}
\caption{\footnotesize{Longest elements in Weyl groups,
see Tables \ref{tab_highest_roots_ABCD} and \ref{tab_highest_roots_EFG}. 
See Remark \ref{rem_differ_indices} on definitions of  
$\alpha^{a5}_{max}$}, $\alpha^{d4}_{max}$ and $\alpha^{d6}_{max}$.}
\label{tab_factoriz_all}
\end{table}
  %% \footnotetext[1]{See Remark \ref{rem_differ_indices}.}

\begin{table}[H]
\centering
\renewcommand{\arraystretch}{1.1}
\begin{tabular}{|l|l|c|c|}
  \hline
       \footnotesize {Root}          & \qquad \qquad \qquad \footnotesize{ Max-orthogonal set}
          & \footnotesize$n$ & \footnotesize{$l_a(w_0)$} \\
       \footnotesize {system}        &  &  & \\
  \hline
       $A_n$  &  \footnotesize$\begin{array}{ll}
   &  \alpha^{a2}_{max} < \alpha^{a4}_{max}  < \dots < \alpha^{a,2k}_{max} \\
   & \\
   \quad \alpha_{k+1} < & \alpha^{a3}_{max} < \alpha^{a5}_{max}  < \dots < \alpha^{a,2k+1}_{max} \\
      \end{array}$ &
      \footnotesize$\begin{array}{l}
        n = 2k \\
        \\
        n = 2k+1
      \end{array}$ &
   \footnotesize$\begin{array}{l}
        \; \frac{n}{2} \\
        \\
        \frac{n+1}{2}
      \end{array}$  \\
  \hline
       $B_n$  & \footnotesize$\begin{array}{ll}
   &  \alpha^{b2}_{max} < \alpha^{b4}_{max}  < \dots  < \alpha^{bn}_{max},    \\
   &  \quad \text{ simple roots: } \alpha_1, \alpha_3, \dots, \alpha_{n-1}, \\
   & \\
   &  \alpha^{b3}_{max} < \alpha^{b5}_{max}  < \dots < \alpha^{bn}_{max},  \\
   &  \quad \text{ simple roots: } \alpha_1, \alpha_3, \dots, \alpha_{n-2}, \alpha_n. \\
     \end{array}$ &
      \footnotesize$\begin{array}{l}
        n = 2k \\
        \\
        n = 2k+1 \\
        \\
      \end{array}$ &
         \footnotesize$n$  \\
  \hline
       $C_n$  &  \footnotesize$\begin{array}{lll}
    & \alpha_n < \alpha^{c2}_{max} < \alpha^{c3}_{max}  < \dots  < \alpha^{cn}_{max}\;. &  \\
              \end{array}$  &
      \footnotesize$\begin{array}{l}
        \text{any} \\
      \end{array}$ &
      \footnotesize$n$  \\
  \hline
       $D_n$  &  \footnotesize$\begin{array}{ll}
 &  \alpha^{d4}_{max} < \alpha^{d6}_{max}  < \dots  < \alpha^{d,n-2}_{max} < \alpha^{dn}_{max} \\
 &  \quad \text{ simple roots: } \alpha_1, \alpha_3, \dots, \alpha_{n-3} \text{ and } \alpha_{n-1}, \alpha_n,      \\
 & \\
 &  \alpha^{d3}_{max} < \alpha^{d5}_{max}  < \dots < \alpha^{d,n-2}_{max} < \alpha^{dn}_{max}    \\
 &  \quad \text{ simple roots: }  \alpha_1, \alpha_3, \dots, \alpha_{n-4}, \alpha_{n-2}
      \end{array}$  &
      \footnotesize$\begin{array}{l}
        n = 2k \\
        \\
        n = 2k+1 \\
        \\
      \end{array}$ &
     \footnotesize$\begin{array}{c}
       n  \\
       \\
       n-1 \\
       \\
      \end{array}$
      \\
  \hline
       $E_6$  & \footnotesize$\begin{array}{ll}
                  \quad \alpha_4 < \alpha^{a3}_{max} < \alpha^{a5}_{max} < \alpha^{e6}_{max}.
                    &  \\
                \end{array}$  & \footnotesize$6$ & \footnotesize$4$  \\
  \hline
       $E_7$  & \footnotesize$\begin{array}{ll}
                & \alpha^{d4}_{max} < \alpha^{d6}_{max} < \alpha^{e7}_{max},  \\
                & \quad \text{ simple roots: } \alpha_2, \alpha_3, \alpha_5, \alpha_7
                \end{array}$  & \footnotesize$7$ & \footnotesize$7$  \\
  \hline
       $E_8$ & \footnotesize$\begin{array}{ll}
            & \alpha^{d4}_{max} < \alpha^{d6}_{max} < \alpha^{e7}_{max} < \alpha^{e8}_{max}.  \\
            & \quad \text{ simple roots: } \alpha_2, \alpha_3, \alpha_5,  \alpha_7,
          \end{array}$  & \footnotesize$8$ & \footnotesize$8$ \\
  \hline
      $F_4$ &  \footnotesize$\begin{array}{ll}
            & \alpha_2 <  \alpha^{c2}_{max} < \alpha^{c3}_{max} < \alpha_{max}.  \\
              \end{array}$  & \footnotesize$4$ & \footnotesize$4$ \\
  \hline
      $G_2$ &  \footnotesize$\begin{array}{ll}
            & \alpha_1 <  \alpha_{max}   \\
            \end{array}$  & \footnotesize$2$ & \footnotesize$2$ \\
   \hline
\end{tabular}
\vspace {1mm}
\caption{\footnotesize{Max-orthogonal set and absolute length function $l_a(w_0)$.}}
\label{tab_max_orhog_set}
\end{table}

\subsection{B.Kostant and A.Joseph: The cascade of orthogonal roots}
  \label{sec_cascade}

\subsubsection{The cascade construction}
A subset $\{\beta_1, \beta_2, ..., \beta_r \in \varPhi \}$ is said to be a {\it strongly orthogonal}
set of roots if $\beta_i \pm \beta_j$ is not a root  for all pairs $\{i,j\}$.
Denote by $\Delta_{\lambda}$ the subset of roots orthogonal to $\lambda$.
The sequence of roots obtained by taking the highest root $\alpha_{max}$ of $\varPhi$,
the corresponding highest roots of components of $\Delta_{\alpha_{max}}$ and so on
(see $\S$\ref{sec_cascade_example}), is a strongly orthogonal set.
In \cite[p.5]{Jo76},  Joseph  refers to a private communication with Kostant,
where the latter notes that: {\it ``\dots any orthogonal
set of roots determines a strongly orthogonal set (by taking sums and differences)
and any maximal strongly orthogonal set is unique up to W''}.
In his following paper \cite{Jo77}, Joseph finds maximal such sets for each
root system and uses them for calculations in the universal
enveloping algebra $U(\mathfrak{g})$, where $\mathfrak{g}$ is a semisimple Lie algebra.
Apparently, these maximal subsets were first presented in \cite[Table III]{Jo77}.
Today they are known as {\it Kostant's cascade}. 

\begin{figure}[h]
\centering
   \includegraphics[scale=0.35]{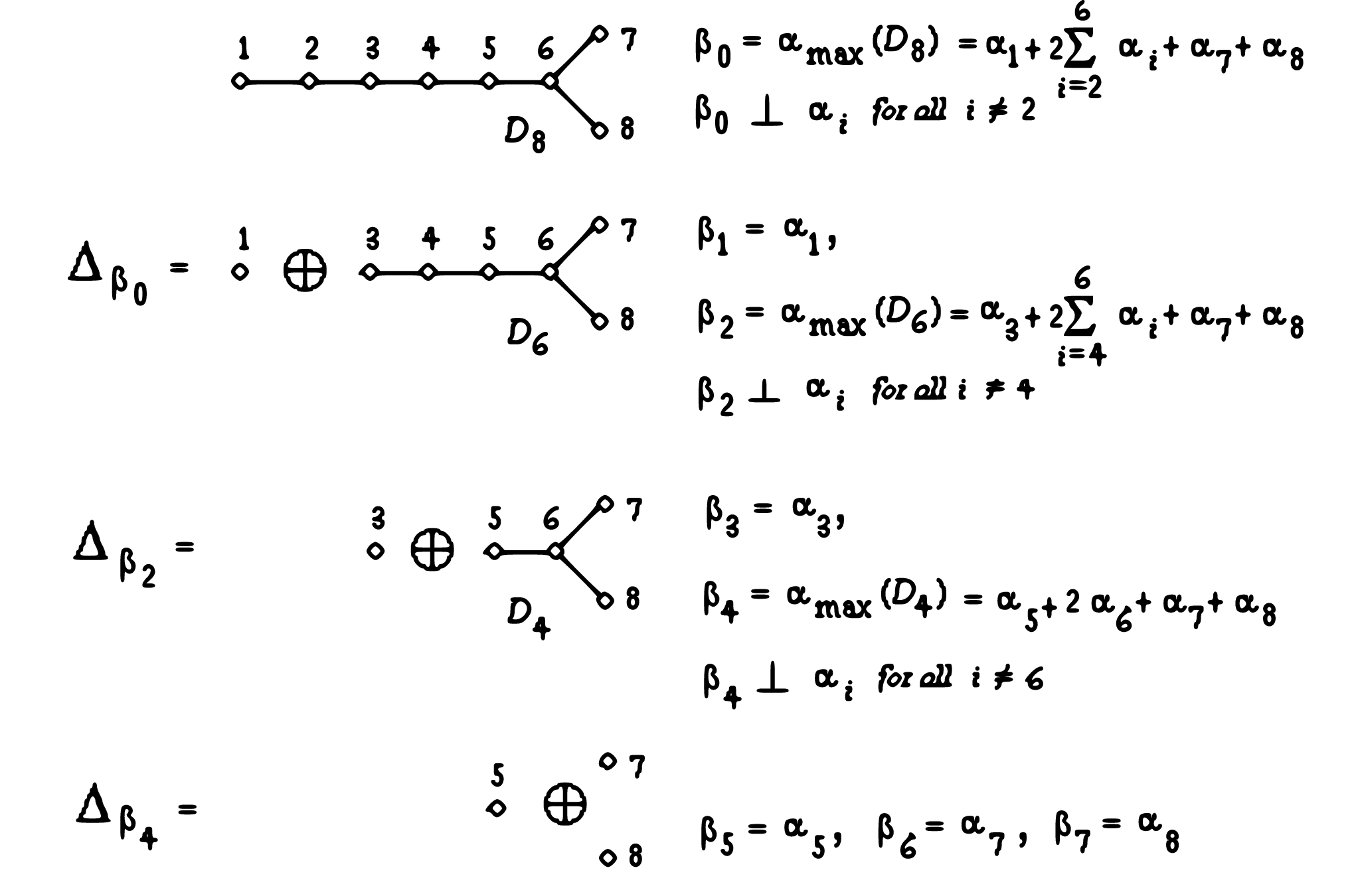}
\caption{The $D_8$-cascade: $\{\beta_0, \beta_1, \beta_2, \beta_3, \beta_4, \beta_5, \beta_6, \beta_7\}$ }
%%%%%% The label must come after caption
\label{fig_cascade_D8}
\end{figure}

     Kostant and Joseph  using the cascade construction, see \cite[$\S$1.1]{Kos13},
independently of each other, using very different methods, obtained a number of structure theorems
for the center $U(\mathfrak{n})$, where $\mathfrak{n}$  is the span of
the positive root spaces in $\mathfrak{g}$. In \cite{DT16},  Dimitrov and Tsanov
obtained a complete list of homogeneous hypercomplex structures on the compact Lie groups
using root subsets called {\it stem}, which were later recognized as the cascade constructed
by Kostant and Joseph. The cascade is also used in \cite{LW82} by Lipsman and Wolf for constructing
certain elements in symmetric algebra
$S(\mathfrak{g})$ and by Panyushev in \cite{Pa22} for classification of a certain class
parabolic subalgebras of $\mathfrak{g}$.

\subsubsection{Example: the $D_8$-cascade}
  \label{sec_cascade_example}
The $D_8$-cascade of orthogonal roots is constructed as follows:

(1) The element $\beta_0$ of cascade is the highest root $\beta_0 = \alpha_{max}(\varPhi)$.

(2) Consider the root subset $\Delta_{\beta_0}$ consisting of roots orthogonal to $\beta_0$:
\begin{equation*}
     \Delta_{\beta_0} = \{\alpha_i \in \varPhi \mid (\alpha_i, \beta_0) = 0 \}.
\end{equation*}
 This subset splits into $2$ connected components:
\begin{equation*}
     \Delta_{\beta_0} = \bigcup_{i_1} \Delta_{0i_1}, \text{ where } \Delta_{00} = \{\alpha_1\} \text{ and }
        \Delta_{01} = \varPhi(D_6).
\end{equation*}
see Fig. \ref{fig_cascade_D8}.

(3) For each  $\Delta_{0i_1}$, we recursively repeat steps (1) and (2). In other words,
we find the highest root $\beta_{0i_1}$ for each component $\Delta_{0i_1}$
which splits into components $\bigcup_{i_2} \Delta_{0i_1i_2}$ and so on.
For details, see \cite[$\S$2.2]{Jo77} or \cite[$\S$3]{LW82}.

(4) For construction of $D_8$-cascade, see Fig. \ref{fig_cascade_D8}. The result $D_8$-cascade is as follows:
\begin{equation*}
  \begin{split}
   & \beta_0 = \alpha_{max}(D_8), \; \beta_1 = \alpha_1, \;  \beta_2 = \alpha_{max}(D_6), \; \beta_3 = \alpha_3, \\
   & \beta_4 = \alpha_{max}(D_4), \; \beta_5 = \alpha_5, \; \beta_6 = \alpha_6, \; \beta_7 = \alpha_7.
  \end{split}
\end{equation*}

\subsection{Bibliographic references}

\subsubsection{Applications of the longest element}
 Let us list a few applications of the longest element.
 \begin{itemize}
 \item In the theory of algebraic groups, the longest element $w_0$  connects
 the Borel subgroup $B$ with the opposite Borel subgroup $B^{-}$ as follows:
 $B^{-} = w_0 B w_0$, \cite{Hu90}.
 \item In the representation theory of Lie algebras,
 $w_0$ maps the highest weight $\lambda$ (corresponding to simple representation $L_{\lambda}$)
 to the highest weight $-w_0\lambda$ (corresponding to the dual representation $L_{\lambda}^{*}$)
 \footnotemark[1].
 \footnotetext[1]{See  J.~Humphreys' article ``Longest element of a finite Coxeter group'' (2015)
  available at http://people.math.umass.edu/$\sim$yeh/pub/longest.pdf}
 \item The longest element of the finite Weyl group plays a special role in the decomposition formula
 of elements $C_w$ of canonical basis (or Kazhdan-Lusztig basis) in the Hecke algebra
 $\mathscr{H}(W)$, see \cite{Bl09}, \cite{Xie17}.
 \item In the context of  the q-analogue of the universal enveloping algebra $\mathscr{U}$,
 it was defined the canonical basis $\mathscr{B}$ of $\mathscr{U}^{+}$
 (the positive part of $\mathscr{U}$), which carries a natural coloured graph structure $\mathscr{G}$.
 Lusztig gave a description of $\mathscr{G}$ using
 the Weyl group longest element, see  \cite{Gr02}, \cite{Ka90}, \cite{Lu90a}, \cite{Lu90b}, \cite{Re95}.
 \item  In the quiver varieties area, Nakajama introduced  the reflection functor for the longest element
 of the Weyl group and identified it with the Lusztig's {\it new symmetry}
 for cases $ADE$, see \cite{Ki16}, \cite{Na03}, \cite{Lu2000}. %%
 \item For more references on applications of the longest elements, see
 \cite{BB05}, \cite{YLi19}, \cite{Lu15}, \cite{Stem97}.
\end{itemize}

\subsubsection{Decompositions}
\begin{itemize}
\item Information about the reduced expressions of the longest elements
is obtained from my calculations, see \cite{St22}, and partly from \cite{YY21}, \cite{Li98} and
\cite[Table 1]{BKOP14}.
\item Perhaps Carter was the first who noticed that {\it each involution} in the Weyl
group can be decomposed  into a product of reflections corresponding to mutually orthogonal
roots, and the number of these reflections can be very small.
However, Carter's decomposition is not unique.
\cite[Lemma 5]{Ca72}.
\item In \cite{Stan84}, Stanley proposed the simple formula for the number of reduced decompositions
of the longest  element in $W(A_n)$. The number of these decompositions reaches
quite large numbers. For example, for $W(A_5)$ this number is $292864$.
\item  The so-called {\it totally orthogonal set} of roots was introduced by Deodhar.
It was proved that for any  involution $w \in W$, there exists a unique totally orthogonal set
$\{\varphi_1, \dots, \varphi_k \}$, such that $w_0 = s_{\varphi_1} \dots s_{\varphi_k}$,
\cite[Theorem 5.4]{Deo82}.  Deodhar's decomposition assumes the existence
of some lexicographic order in the root system.
\item After all the calculations for this paper were completed, I realized that the idea
of the decomposition implemented here is close to that considered by Springer in \cite{Sp82}:
  {\it Let $w \in W$ be an involution. There is a unique special parabolic subsystem $J$
  such that $w$ is the longest element of $W_J$}, see the proof in \cite[Prop.3]{Sp82}.
\end{itemize}
~\\

\section{\bf Relations in the Weyl groups}

\begin{lemma}[$\Lambda\rm{V}$-relation]
  \label{lem_lambda_V}
  Let $\{\alpha_k, \alpha_{k+1}, ..., \alpha_{n-1}, \alpha_n\}$ is the sequence of simple roots
  corresponding to the connected subdiagram $A_{n-k+1}$ of the Dynkin diagram $A_l$, such that
  only roots with adjacent indices are non-orthogonal. Let $s_i$ be reflections
  in the Weyl group $W(A_{n-k+1})$ associated with roots $\alpha_i$, where $k \leq i \leq n$. Then,
\begin{equation}
 \label{lambda_V_rel}
    s_k s_{k+1} \dots s_{n-1} s_n s_{n-1}\dots s_{k+1}s_k  =
         s_n s_{n-1} \dots s_{k+1}s_k s_{k+1} \dots s_{n-1} s_n.
\end{equation}
\end{lemma}

\begin{figure}[h]
\centering
   \includegraphics[scale=0.18]{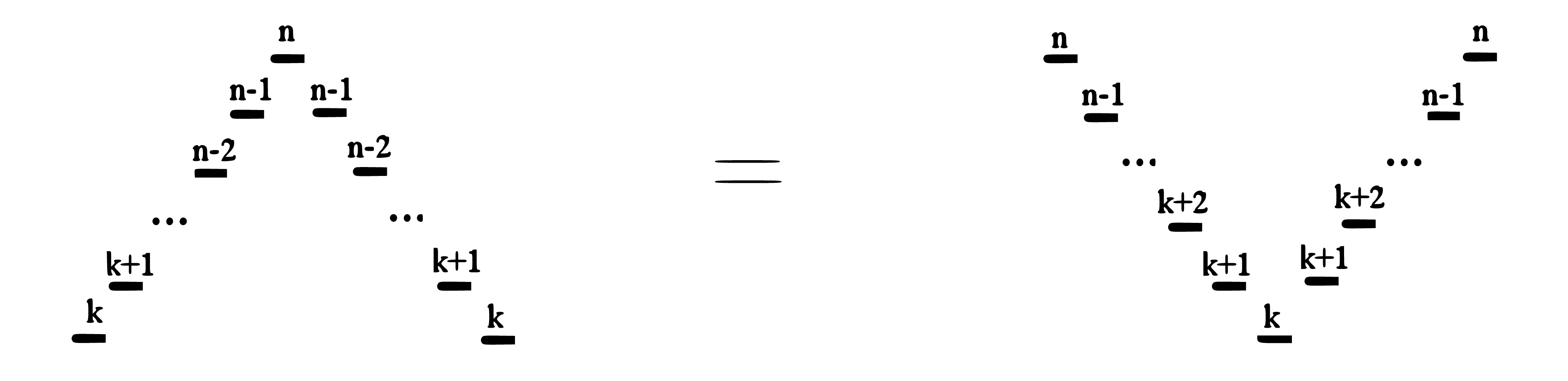}
\caption{The stepwise shape of indices in $\Lambda\rm{V}$-relations.}
%%%%%% The label must come after caption
\label{fig_lambda_V}
\end{figure}

 Relation \eqref{lambda_V_rel} is said to be the $\Lambda\rm{V}${\it-relation}
in accordance with the stepwise form by which the indices are formed, see Fig. \ref{fig_lambda_V}.

\begin{example}
 \label{exam_braid_1}
\begin{equation*}
  \footnotesize
\begin{split}
n = 1,\quad & s_{n-1} s_n s_{n-1} = s_{n} s_{n-1} s_{n}, \;{\rm braid \; relation}, \\
n = 2,\quad & s_{n-2} s_{n-1} s_n s_{n-1} s_{n-2} = s_{n} s_{n-1} s_{n-2} s_{n-1}s_{n}, \; k = n-2, \\
n = 3,\quad & s_1 s_2 s_3 s_4 s_3 s_2 s_1 = s_4 s_3 s_2 s_1 s_2 s_3 s_4, \; k = 1, n = 4, \; {\rm since}\\
        s_1 s_2 (s_3 & s_4 s_3) s_2 s_1 = s_1 s_2 s_{\alpha_3 + \alpha_4} s_2 s_1 =
     s_1 s_{\alpha_2 + \alpha_3 + \alpha_4} s_1 = s_{\alpha_1 + \alpha_2 + \alpha_3 + \alpha_4},
             \;{\rm and } \\
     s_4 s_3 (s_2 & s_1 s_2) s_3 s_4 = s_4 s_3 s_{\alpha_1 + \alpha_2} s_3 s_4 =
          s_4 s_{\alpha_1 + \alpha_2 + \alpha_3} s_4 = s_{\alpha_1 + \alpha_2 + \alpha_3 + \alpha_4}. \qed \\
\end{split}
\end{equation*}
\end{example}

\pprof{Lemma \ref{lem_lambda_V}}
 Note that expressions standing in both sides of
the $\Lambda\rm{V}$-relation coincide with the reflection associated with the root
\begin{equation*}
  \alpha_k + \alpha_{k+1} + \dots + \alpha_{n-1} + \alpha_n.
\end{equation*}
Eq. \eqref{lambda_V_rel} is equivalent to $2$ following equalities:
\begin{equation}
  \label{lambda_V_rel_root}
  \begin{split}
    &  s_k s_{k+1} \dots s_{n-1} s_n s_{n-1}\dots s_{k+1}s_k  =
               s_{\alpha_k + \alpha_{k+1} + \dots + \alpha_{n-1} + \alpha_n} \\
    &  s_n s_{n-1} \dots s_{k+1}s_k s_{k+1} \dots s_{n-1} s_n  =
               s_{\alpha_k + \alpha_{k+1} + \dots + \alpha_{n-1} + \alpha_n} \\
   \end{split}
\end{equation}
The induction step for the first (resp. second) equality in \eqref{lambda_V_rel_root} is
conjugation with $s_k$ (resp. $s_n$) and using Proposition \ref{prop_fact_1}:
\begin{equation*}
  \begin{split}
   & s_k s_{\alpha_{k+1} + \dots + \alpha_{n-1} + \alpha_n}s_k  =
               s_{\alpha_k + \alpha_{k+1} + \dots + \alpha_{n-1} + \alpha_n}, \\
   & s_n s_{\alpha_{k} + \alpha_{k+1} + \dots + \alpha_{n-1}} s_n  =
               s_{\alpha_k + \alpha_{k+1} + \dots + \alpha_{n-1} + \alpha_n}. \qed \\
   \end{split}
\end{equation*}

\begin{lemma}[on permutation]
  \label{lem_permute_An_1}
  Let $\{ \alpha_k,\alpha_{k+1}, \dots, \alpha_n \}$ be the sequence of simple roots in
  the root system $A_l$ corresponding to the set of consecutive vertices $\{k, k+1, \dots, n \}$ in
  the Dynkin diagram $A_l$, $k \leq n \leq l$. Then,
\begin{equation}
  \label{eq_prop_An_1}
   s_{\alpha_k+\dots \alpha_{n-2}+\alpha_{n-1}} \prod_{n \geq j \geq k}s_j =
      \big( \prod_{n-1 \geq j \geq {k+1}} s_j \big )s_{\alpha_k+\alpha_{k+1}+\dots{\alpha_{n-1}+\alpha_n}}
\end{equation}
\end{lemma}

\begin{example}
 \label{exam_An_1}
   {\rm Let us check eq.} \eqref{eq_prop_An_1} {\rm for} $n = 5$, $k = 1$.
\begin{equation*}
\begin{split}
  {\rm Since }\;
  & s_{\alpha_1+\alpha_2+\alpha_3+\alpha_4} = s_4s_3s_2s_1s_2s_3s_4, \; {\rm then } \\
  & s_{\alpha_1+\alpha_2+\alpha_3+\alpha_4} s_5s_4s_3s_2s_1 = 
      s_4s_3s_2(s_1s_2s_3s_4s_5s_4s_3s_2s_1) = \qquad \qquad \qquad \\
        & s_4s_3s_2s_{\alpha_1+\alpha_2+\alpha_3+\alpha_4+\alpha_5}. \qed
\end{split}
\end{equation*}
\end{example}

\pprof{Lemma \ref{lem_permute_An_1}}  As in Example \ref{exam_An_1}, using Lemma \ref{lem_lambda_V},
we have

\begin{equation*}
\begin{split}
   &  s_{\alpha_k+\dots \alpha_{n-2}+\alpha_{n-1}} \prod_{n \geq j \geq k}s_j \; = \;
      s_{\alpha_k + \dots + \alpha_{n-2}+\alpha_{n-1}}  s_ns_{n-1} \dots s_k = \\
   & s_{\alpha_{n-1}}s_{\alpha_{n-2}}\dots s_{\alpha_{k+1}}s_{\alpha_k}
        s_{\alpha_{k+1}} \dots s_{\alpha_{n-1}}
        (s_{\alpha_n}s_{\alpha_{n-1}}\dots s_{k+1}s_k) = \\
   &  (s_{\alpha_{n-1}}s_{\alpha_{n-2}}\dots s_{\alpha_{k+1}})
         (s_{\alpha_k}s_{\alpha_{k+1}} \dots s_{\alpha_{n-1}}
         s_{\alpha_n}s_{\alpha_{n-1}}\dots s_{k+1}s_k) = \\
   & \big( \prod_{n-1 \geq j \geq {k+1}} s_j \big )s_{\alpha_k+\alpha_{k+1}+\dots{\alpha_{n-1}+\alpha_n}}. \qed
\end{split}
\end{equation*}

\section{\bf Decomposition in $W(A_n)$}
\begin{figure}[h]
\centering
   \includegraphics[scale=0.4]{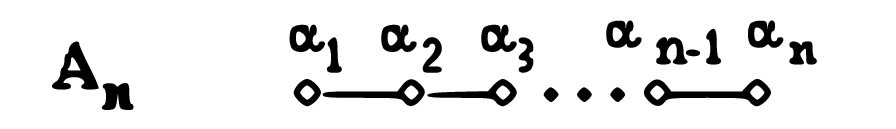}
\caption{\footnotesize{Numbering of simple roots in $A_n$}}
%%%%%% The label must come after caption
\label{fig_An}
\end{figure}

Let us introduce the notation for the highest roots in the root subsystems $A_i \subset A_n$,
see Table \ref{tab_factoriz_all}. The numbering of the vertices of the Dynkin diagram $A_n$ is shown
in Fig. \ref{fig_An}.
\begin{equation}
  \label{eq_def_An}
      \alpha^{a,2i}_{max} := \sum\limits_{j = k-i+1}^{k+i} \alpha_j, \;\;
      \alpha^{a,2i+1}_{max} := \sum\limits_{j = k-i+1}^{k+i+1} \alpha_j, \;\;
      \text{ where } i \ge 1. \\
\end{equation}
see Remark \ref{rem_notations}. For $i = 0$, we put $\alpha^{a1}_{max} := \alpha_{k+1}$.
By \eqref{eq_def_An}, the root $\alpha^{ap}_{max}$ is the \underline{highest root} in the root system $A_p$
for any $p \ge 0$.

\begin{proposition}
  \label{prop_factoriz_An}
  {\rm{(i)}}
The element
\begin{equation}
  \label{prop_lst_in_An}
     w_0 = s_1(s_2s_1)(s_3s_2s_1)\dots (s_ns_{n-1}{\dots}s_2s_1)
\end{equation}
is the longest element in $W(A_n)$.

  {\rm{(ii)}}  The following equality holds:
\begin{equation}
 \begin{split}
  \label{prop_lst_in_An_2}
     s_1(s_2s_1)(s_3s_2s_1)\dots & (s_ns_{n-1}{\dots}s_2s_1) = \\
     & s_n(s_{n-1}s_n)(s_{n-2}s_{n-1}s_n)\dots (s_1s_2{\dots}s_{n-1}s_n).
 \end{split}
\end{equation}
\end{proposition}

\PerfProof
(i) Let us denote the factors entering in \eqref{prop_lst_in_An} as follows:
\begin{equation}
  \label{multipl_lst_in_An}
  w_1 = s_1, \; w_2 = s_2s_1, \; \dots, \;  w_{n-1} = s_{n-1}{\dots}s_2s_1,  \;
  w_{n} = s_{n}{\dots}s_2s_1.
\end{equation}

Then, $w_0 = w_1w_2\dots w_n$. For any integer $k$ such that $1 \leq k\leq n$, we have
\begin{equation}
  \label{eq_lst_in_An_2}
  \begin{split}
    & w_n: \alpha_k \longrightarrow \alpha_{k-1}, \\
    & w_{n-1}: \alpha_{k-1} \longrightarrow \alpha_{k-2}, \\
    & \dots \\
    & w_{n-k+2}:  \alpha_2 \longrightarrow \alpha_1, \\
    & w_{n-k+1}:  \alpha_1 \longrightarrow -(\alpha_1+\alpha_2+ \dots + \alpha_{n-k+1}), \\
    & w_{n-k}:   -(\alpha_1+\alpha_2+ \dots + \alpha_{n-k+1}) \longrightarrow -\alpha_{n-k+1}. \\
    & \text{Elements } w_{n-k-1}, w_{n-k-2}, \dots, w_1 \text{ preserve } \alpha_{n-k+1}.
  \end{split}
\end{equation}
Let us explain \eqref{eq_lst_in_An_2}. We have
\begin{equation*}
  \footnotesize
  \label{eq_lst_in_An_expl2}
  \begin{array}{ll}
  & w_n(\alpha_k) = s_{\alpha_n}s_{\alpha_{n-1}} \dots s_{\alpha_{k}}s_{\alpha_{k-1}} \dots s_{\alpha_2}s_{\alpha_1} (\alpha_k) = \\
  & \qquad  s_{\alpha_n}s_{\alpha_{n-1}} \dots s_{\alpha_{k}}s_{\alpha_{k-1}}(\alpha_k) =
    s_{\alpha_n}s_{\alpha_{n-1}}\dots s_{\alpha_{k}}(\alpha_k+\alpha_{k-1}) = \\
  & \qquad s_{\alpha_n}s_{\alpha_{n-1}}\dots s_{\alpha_{k+1}}(\alpha_{k-1}) = \alpha_{k-1}. \\
  &  \text{Similarly, } \\
  & w_{n-1}(\alpha_{k-1}) = \alpha_{k-2}, \dots, w_{n-k+2}(\alpha_2) = \alpha_1,  \\
  & w_{n-k+1}(\alpha_{1})  =  s_{\alpha_{n-k+1}}\dots s_{\alpha_1} (\alpha_1) =
   s_{\alpha_{n-k+1}} \dots s_{\alpha_2}(-\alpha_1) = \\
  & \qquad  -(\alpha_1 +\dots + \alpha_{n-k+1}).\\
  & w_{n-k}(-(\alpha_1 +\dots + \alpha_{n-k+1})) =
    s_{\alpha_{n-k}}\dots  s_{\alpha_1} (-(\alpha_1 + \alpha_2 + \dots + \alpha_{n-k+1}) = \\
  & \qquad s_{\alpha_{n-k}}\dots  s_{\alpha_2} (-(\alpha_2 +\dots + \alpha_{n-k+1})) = \dots = \\
  & \qquad  s_{\alpha_{n-k}}s_{\alpha_{n-k-1}}(-(\alpha_{n-k-1} + \alpha_{n-k} + \alpha_{n-k+1})) = \\
  & \qquad  s_{\alpha_{n-k}}(-(\alpha_{n-k} + \alpha_{n-k+1})) = -\alpha_{n-k+1}.
  \end{array}
\end{equation*}

By \eqref{eq_lst_in_An_2} $w_0$ transforms $\alpha_i$ to $-\alpha_{n-i+1}$ for $i = 1,2\dots,n $  
and $n \geq 2$, see \cite[Plate 1]{Bo02}. Because the longest element is unique in $W(A_4)$, $w_0$  is the longest.
~\\

(ii) Here, we just rename $s_i$ to $s_{n-i+1}$.
 \qed

\begin{example} \footnotesize
\begin{equation*}
  \begin{split}
    & n = 2, \quad w_0 = s_1(s_2s_1) = \underline{s_{\alpha_1 + \alpha_2}}, \\
    & n = 3, \quad w_0 = s_1(s_2s_1)(s_3s_2s_1) = (s_2s_1s_2)(s_3s_2s_1) =
        \underline{s_{\alpha_2} s_{\alpha_1 + \alpha_2 + \alpha_3}}, \\
    & n = 4, \quad \rm{by\; Lemma }\; \ref{lem_permute_An_1}, \rm{\;using\; the\; relation\; for\;} 
                n=3: \\
    & \qquad w_0 = s_1(s_2s_1)(s_3s_2s_1)(s_4s_3s_2s_1) =
               s_{\alpha_2} s_{\alpha_1 + \alpha_2 + \alpha_3}(s_4s_3s_2s_1) = \\
    & \qquad \quad  s_2(s_3s_2) s_{\alpha_1 + \alpha_2 + \alpha_3 + \alpha_4} =
          \underline{s_{\alpha_2 + \alpha_3} s_{\alpha_1 + \alpha_2 + \alpha_3 + \alpha_4}}, \\
    &  n = 5, \quad \rm{by\; Lemma }\;\ref{lem_permute_An_1}, \rm{\;using\; the\; relation\; for\;}
                n=4: \\
    &  \quad w_0 = s_1(s_2s_1)(s_3s_2s_1)(s_4s_3s_2s_1)(s_5s_4s_3s_2s_1) = \\
    &  \qquad s_{\alpha_2 + \alpha_3}s_{\alpha_1 + \alpha_2 + \alpha_3 + \alpha_4}(s_5s_4s_3s_2s_1) = \\
    &  \qquad s_{\alpha_2 + \alpha_3}(s_4s_3s_2)
                       s_{\alpha_1 + \alpha_2 + \alpha_3 + \alpha_4 + \alpha_5} =
  \underline{s_{\alpha_3}s_{\alpha_2 + \alpha_3 + \alpha_4}
       s_{\alpha_1 + \alpha_2 + \alpha_3 + \alpha_4 + \alpha_5}}. \\
  \end{split}
\end{equation*}
\end{example}
~\\
\begin{proposition}[decomposition in $W(A_n)$]
  \label{prop_factoriz_An_2}
{\rm{(i)}}
The longest element in $W(A_n)$ is decomposed into $\displaystyle \big [ \frac{n+1}{2} \big ]$
factors as follows:  
\begin{equation}
  \label{eq_factors_An_1}
  \begin{split}
     w_0 = s_1 & (s_2s_1) (s_3s_2s_1)\dots (s_ns_{n-1}{\dots}s_2s_1) = \\
      \Biggl \{ & \begin{array}{lll}
           & \prod\limits_{i=1}^k s_{\alpha^{a,2i}_{max}} & \text{ for } n = 2k, \\
           s_{\alpha_{k+1}} & \prod\limits_{i=1}^k s_{\alpha^{a,2i+1}_{max}} & \text{ for } n = 2k+1.            
       \end{array}    
   \end{split}
\end{equation}

{\rm{(ii)}} The roots appearing in the decomposition
\eqref{eq_factors_An_1} are  \underline{mutually orthogonal}. Thus,
the factors in \eqref{eq_factors_An_1} commute with each other.
The roots $\alpha^{ap}_{max}$ are ordered as follows:
\begin{equation*}
 \label{eq_highest_roots_An}
 \begin{array}{lll}
   &  \alpha^{a2}_{max} < \alpha^{a4}_{max}  < \dots < \alpha^{a,2k}_{max}& \text{for } \; n = 2k,  \\
   \alpha_{k+1} < & \alpha^{a3}_{max} < \alpha^{a5}_{max}  < \dots < \alpha^{a,2k+1}_{max}& \text{for } \; n = 2k+1.  \\
 \end{array}
\end{equation*}
\end{proposition}
~\\

\PerfProof (i)
The product of eq. \eqref{eq_factors_An_1} is $w_0 = w_1...w_n$, see \eqref{prop_lst_in_An}, \eqref{multipl_lst_in_An}.
Let us discard the last factor $w_n$: $w_0' = w_1 \dots w_{n-1}$.
Using the induction step, we replace $w_0'$ with the second line of \eqref{eq_factors_An_1},
while $n$ is replaced by $n-1$ and $k+1$ by $k$.
Further, we use relation of Lemma \ref{lem_permute_An_1} and eq. \eqref{eq_def_An}:
\begin{equation*}
  \begin{split}
     & s_1(s_2s_1)(s_3s_2s_1)\dots (s_ns_{n-1}{\dots}s_2s_1) = \\
     & s_{\alpha_k} s_{\alpha_{k-1} + \alpha_k + \alpha_{k+1}} \dots
           s_{\alpha_2 + \dots +\alpha_{n-2}} s_{\alpha_1 + \dots +\alpha_{n-1}}(s_ns_{n-1}{\dots}s_2s_1) =  \\
     & s_{\alpha_k}s_{\alpha_{k-1} + \alpha_k + \alpha_{k+1}} \dots
       s_{\alpha_3 + \dots +\alpha_{n-3}} s_{\alpha_2 + \dots +\alpha_{n-2}}(s_{n-1}{\dots}s_2) 
       s_{\alpha_{max}} = \\
     & s_{\alpha_k}s_{\alpha_{k-1} + \alpha_k + \alpha_{k+1}} \dots  
       s_{\alpha_3 + \dots +\alpha_{n-3}}(s_{n-2}{\dots}s_3)
                \prod\limits_{i={k-1}}^k s_{\alpha^{a,2i}_{max}} = \\
     & s_{\alpha_k}s_{\alpha_{k-1} + \alpha_k + \alpha_{k+1}} \dots
       s_{\alpha_4 + \dots +\alpha_{n-4}}(s_{n-3}{\dots}s_4)
                 \prod\limits_{i={k-2}}^k s_{\alpha^{a,2i}_{max}} = \\
     &  \dots \\
     & s_{\alpha_k}s_{\alpha_{k-1} + \alpha_k + \alpha_{k+1}}(s_{k+2}s_{k+1}s_k s_{k-1}) 
       \prod\limits_{i=3}^k s_{\alpha^{a,2i}_{max}} = \\
     & s_{\alpha_k}(s_{k+1}s_k )s_{\alpha_{k-1} + \alpha_k + \alpha_{k+1} + \alpha_{k+2}}
       \prod\limits_{i=3}^k s_{\alpha^{a,2i}_{max}} = 
       s_{\alpha_k}(s_{k+1}s_k )s_{\alpha^{a4}_{max}}
       \prod\limits_{i=3}^k s_{\alpha^{a,2i}_{max}} = \\
     & (s_k s_{k+1}s_k)\prod\limits_{i=2}^k s_{\alpha^{a,2i}_{max}} = 
       s_{\alpha_k + \alpha_{k+1}}\prod\limits_{i=2}^k s_{\alpha^{a,2i}_{max}} =
       \prod\limits_{i=1}^k s_{\alpha^{a,2i}_{max}}.
   \end{split}
\end{equation*}
 Similarly, the second equality of \eqref{eq_factors_An_1} is proved.
~\\

(ii) The orthogonality relations are verified directly.
\qed
~\\

\section{\bf Decomposition in $W(B_n)$}

\begin{figure}[h]
\centering
   \includegraphics[scale=0.4]{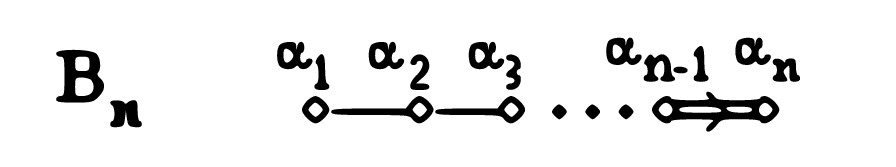}
\caption{Numbering of simple roots in $B_n$}
%%%%%% The label must come after caption
\label{fig_Bn}
\end{figure}

Let us introduce the notation for the highest roots in the root subsystem $B_i \subset B_n$,
see Table \ref{tab_factoriz_all}. The numbering of vertices of the Dynkin diagram $B_n$ is shown
in Fig. \ref{fig_Bn}.
\begin{equation*}
  \label{eq_not_Bn_Cn_1}
     \alpha^{bi}_{max} := \alpha_{n-i+1} + 2\sum\limits_{j = n-i+2}^n \alpha_j,
      \text{ where } i \ge 2.
\end{equation*}
see Remark \ref{rem_notations}. The root $s_{\alpha^{bp}_{max}}$ is the \underline{highest root}
in the root system $B_p$ for any $p \geq 2$. Let
\begin{equation*}
  \{ \varepsilon_i \mid i = 1,2, \dots n \}
\end{equation*}
be the orthogonal basis of unit vectors in the linear space $\mathbb{E}_n$.
The simple roots $\{ \alpha_i \mid i = 1,\dots, n \}$ in $B_n$ are as follows:
\begin{equation}
 \label{eq_compare_with_Bn}
 \begin{split}
    &  \alpha_i = \varepsilon_i - \varepsilon_{i+1}
       \text{ for } i = 1, \dots,n-1, \text{ and } \alpha_n = \varepsilon_n, \\
    & \norm{\alpha_i} = \sqrt{2}  \; \text { for }  \; i = 1,\dots,n-1, \text{ and }
       \norm{\alpha_n} = 1,
 \end{split}
\end{equation}
see \cite[Plate II]{Bo02}.
Then,
\begin{equation}
  \label{eq_def_ei_in_Bn_1}
    \varepsilon_i = \alpha_n + \alpha_{n-1} + \dots + \alpha_i, \;\; \text{ for } i < n, \; \text{ and } \;
    \varepsilon_n = \alpha_n.
\end{equation}

\begin{proposition}
   The element
\begin{equation}
 \label{eq_lnst_Bn_1}
   w_0 = s_n(s_{n-1}s_ns_{n-1})(s_{n-2}s_{n-1}s_ns_{n-1}s_{n-2}) \dots
     (s_1s_2\dots s_n \dots s_2s_1)
\end{equation}
is the longest element in $W(B_n)$.
\end{proposition}

\PerfProof Denote the factors of \eqref{eq_lnst_Bn_1} as follows:
\begin{equation*}
  w_n = s_n, \; w_{n-1} = s_{n-1}s_ns_{n-1}, \; \dots, \;
      w_1 = s_1s_2 \dots s_n \dots s_2s_1.
\end{equation*}
The following relations hold:
\begin{equation}
  \label{eq_wi_in_Bn}
  \begin{array}{lll}
    & w_i = s_i w_{i+1} s_i &\text{ for } \quad i=1,\dots,n-1, \\
    & w_i = s_{\alpha_n + \alpha_{n-1} + \dots + \alpha_i}
              & \text{ for } \quad i=1,\dots,n.
  \end{array}
\end{equation}
By \eqref{eq_def_ei_in_Bn_1} $w_i = s_{\varepsilon_i}$ for $i=1,\dots,n$.
The roots $\varepsilon_i$ are short, and we use
Proposition \ref{prop_fact_1} for conjugation in \eqref{eq_wi_in_Bn}.  Then,
\begin{equation}
  \label{eq_not_fin_factoriz_Bn}
   w_0 = \prod\limits_{i=1}^n s_{\varepsilon_i}.
\end{equation}
where all $s_{\varepsilon_i}$  commute with each other.
The element $w_0$ acts on $\varepsilon_i$ as $s_{\varepsilon_i}$:
\begin{equation*}
   w_0(\varepsilon_i) = s_{\varepsilon_i}(\varepsilon_i) = -\varepsilon_i.
\end{equation*}
Thus,  $w_0 = -1$ on $\mathbb{E}_n$, see \cite[Plate II]{Bo02}.
Since the longest element in the Weyl group is unique, $w_0$ is the longest. \qed
~\\

The roots $\varepsilon_i = \alpha_n + \dots + \alpha_i$ are not highest roots, so
the decomposition \eqref{eq_not_fin_factoriz_Bn} is not max-orthogonal, see $\S$\ref{sec_uniqueness}.
We consider another decomposition for the longest element in $W(B_n)$ 
that satisfies the max-orthogonality conditions (Proposition \ref{prop_factoriz_Bn}).
Let us start with the following example using Propositions \ref{prop_fact_1} and \ref{prop_fact_2}.

\begin{example}
 \label{examp_Bn_1}
 \footnotesize
\begin{equation*}
 \begin{split}
    n = 2, & \quad w_0 = (s_2s_1s_2) s_1 = \underline{s_{\alpha_1 + 2\alpha_2} s_{\alpha_1}},\\
    n = 3, & \quad w_0 = s_3s_2(s_3s_2s_1s_2s_3)s_2s_1 =
              s_3(s_2s_{\alpha_1 + \alpha_2 + 2\alpha_3}s_2)s_1 =
              \underline{s_3(s_{\alpha_1 + 2\alpha_2 + 2\alpha_3})s_1}, \\
    n = 4, & \quad w_0 = s_4(s_3s_4s_3)s_2s_3(s_4s_3s_2s_1s_2s_3s_4)s_3s_2s_1) =  \\
           & (s_4s_3s_4)s_3(s_2s_3
               (s_{\alpha_1 + \alpha_2 + \alpha_3 +2\alpha_4})s_3s_2)s_1 =
    \underline{s_{\alpha_3 + 2\alpha_4}s_{\alpha_3}
    s_{\alpha_1 + 2\alpha_2 + 2\alpha_3 +2\alpha_4}s_{\alpha_1}}, \\
    n = 5, & \quad w_0 = s_5(s_4s_5s_4)(s_3s_4s_5s_4s_3)(s_2s_3s_4s_5s_4s_3s_2)
          (s_1s_2s_3s_4s_5s_4s_3s_2s_1) = \\
         & s_5(s_4(s_5s_4s_3s_4s_5)s_4)(s_3)(s_2s_3s_4(s_5s_4s_3s_2s_1s_2s_3s_4s_5)
          s_4s_3s_2)s_1  = \\
         &  s_5(s_4s_{\alpha_3 + \alpha_4 + 2\alpha_5}s_4)(s_3)(s_2s_3s_4
             s_{\alpha_1 + \alpha_2 + \alpha_3 + \alpha_4 + 2\alpha_5}s_4s_3s_2)s_1 = \\
 &  \underline{s_5(s_{\alpha_3 + 2\alpha_4 + 2\alpha_5})
     s_3(s_{\alpha_1 + 2\alpha_2 + 2\alpha_3 + 2\alpha_4 + 2\alpha_5})s_1}. \\
  \end{split}
\end{equation*}
\end{example}
~\\

\begin{proposition}[decomposition in $W(B_n)$]
  \label{prop_factoriz_Bn}
{\rm{(i)}} The longest element in $W(B_n)$ is decomposed into $n$ factors as follows:
\begin{equation}
  \label{eq_factors_Bn_1B}
  \footnotesize
  \begin{split}
    w_0 = s_n & (s_{n-1}s_ns_{n-1})(s_{n-2}s_{n-1}s_ns_{n-1}s_{n-2}) \dots 
       (s_1s_2\dots s_n \dots s_2s_1) = \\       
    & \Biggl \{ 
       \begin{array}{lll}
        & \prod\limits_{\stackrel{i=1}{i - odd}}^{n-1} 
          (s_{\alpha_i}s_{\alpha^{b,n-i+1}_{max}}) 
        &      \text{ for } n = 2k, \\
        \\
        s_{\alpha_n} & \prod\limits_{\stackrel{i=1}{i - odd}}^{n-2}
           (s_{\alpha_i}s_{\alpha^{b,n-i+1}_{max}})
        &      \text{ for } n = 2k+1.
       \end{array}    
  \end{split}
\end{equation}

{\rm{(ii)}} The roots appearing in the decomposition
\eqref{eq_factors_Bn_1B} are \underline{mutually orthogonal}.
Thus, the factors in \eqref{eq_factors_Bn_1B} commute with each other.
 The roots $s_{\alpha^{bp}_{max}}$ are ordered as follows:
\begin{equation*}
 \label{eq_highest_roots_Bn}
 \begin{array}{lll}
   &  \alpha^{b2}_{max} < \alpha^{b4}_{max}  < \dots  < \alpha^{bn}_{max} & \text{ for } n = 2k,  \\
   &  \alpha^{b3}_{max} < \alpha^{b5}_{max}  < \dots < \alpha^{bn}_{max}  & \text{ for } n = 2k+1.  \\
 \end{array}
\end{equation*}
\end{proposition}

\PerfProof  (i) Let us divide all factors in \eqref{eq_factors_Bn_1B}
into groups consisting of pairs of neighboring factors. Consider one such pair:
\begin{equation*}
 \footnotesize
  \label{eq_divide_to_pairs}
   p_i :=
   (s_{i+1}s_{i+2}\dots s_{n-1}s_ns_{n-1} \dots s_{i+1})(s_i s_{i+1}\dots s_{n-1}s_ns_{n-1} 
      \dots s_{i+2}s_{i+1} s_i).
\end{equation*}
where $i = 1=1,3,\dots,n-1$ for $n = 2k$ \ and
$i = 1,3,\dots,n-2$ for $n=2k+1$, see eq. \eqref{eq_factors_Bn_1B}. The element $p_i$
is transformed as follows:
\begin{equation}
 \label{eq_Bn_pi}
 \footnotesize
  \begin{split}
   p_i = & s_{i+1}s_{i+2}\dots s_{n-1}\big ( s_ns_{n-1} \dots s_{i+1}s_i s_{i+1}\dots s_{n-1}s_n \big ) s_{n-1} \dots s_{i+2}s_{i+1} s_i = \\
   & s_{i+1}s_{i+2}\dots s_{n-1}
     \big ( s_n (s_{\alpha_i + \alpha_{i+1} + \dots + \alpha_{n-1}}) s_n \big )
        s_{n-1}\dots s_{i+2}s_{i+1} s_i = \\
    & s_{i+1}s_{i+2}\dots s_{n-1}
      \big ( s_{\alpha_i + \alpha_{i+1} + \alpha_{i+2} + \dots + \alpha_{n-1} +
         2\alpha_n} \big )  s_{n-1}\dots s_{i+2}s_{i+1} s_i = \\
    & s_{i+1}s_{i+2}\dots
      \big ( s_{\alpha_i + \alpha_{i+1} + \alpha_{i+2} + \dots + 2\alpha_{n-1} +
          2\alpha_n} \big )  \dots s_{i+2} s_{i+1} s_i = \\
    & s_{i+1}
       \big ( s_{\alpha_i + \alpha_{i+1} + 2\alpha_{i+2} + \dots + 2\alpha_{n-1} +
          2\alpha_n} \big ) s_{i+1} s_i = \\
    &  s_{\alpha_i + 2\alpha_{i+1} + 2\alpha_{i+2} + \dots + 2\alpha_{n-1} +
          2\alpha_n}  s_i.
  \end{split}
\end{equation}
~\\
In eq. \eqref{eq_Bn_pi}, we use the following two relations:
\begin{equation*}
 \footnotesize
  \begin{split}
  (1) \qquad & (\alpha_i + \alpha_{i+1} + 2\alpha_{i+2} + \dots + 2\alpha_n, \alpha_{i+1}) = -1 + 2 -2 = -1, \text{ then }  \\
     & s_{i+1}( s_{\alpha_i + \alpha_{i+1} + 2\alpha_{i+2} + \dots + 2\alpha_n)} s_{i+1} = 
        s_{(\alpha_i + \alpha_{i+1} + 2\alpha_{i+2} + \dots +  2\alpha_n) + \alpha_{i+1}} = \\
     &   s_{\alpha_i + 2\alpha_{i+1} + 2\alpha_{i+2} + \dots + 2\alpha_n},   \\
  (2) \qquad & (\alpha_i + 2\alpha_{i+1} + 2\alpha_{i+2} + \dots +  2\alpha_n, \alpha_i) = 
         2 - 2 = 0, \text{ then }   \\
     &  s_i \text{ and }s_{\alpha_i + 2\alpha_{i+1} + 2\alpha_{i+2} + \dots + 2\alpha_n} \text{ commute. }     
  \end{split}
\end{equation*}

(ii) By \eqref{eq_def_ei_in_Bn_1} we get
\begin{equation*}
  \label{eq_def_ei_ei1}
  \begin{split}
   & \varepsilon_i + \varepsilon_{i+1} =
       \alpha_i + 2\alpha_{i+1}+ \dots + 2\alpha_{n-1} + 2\alpha_n, \\
   & p_i = s_{\varepsilon_i + \varepsilon_{i+1}}s_{\alpha_i}.
  \end{split}
\end{equation*}
The roots $\varepsilon_i + \varepsilon_{i+1}$ are highest in the root subsystem $B_{n-i+1}$.
Then, for  $n = 2k$, relation \eqref{eq_factors_Bn_1B} has $n$ factors and looks as follows:
\begin{equation}
  \label{eq_product_Bn_even}
  w_0 = \prod\limits_{i = 1,3,...,n-1}s_{\varepsilon_i + \varepsilon_{i+1}}s_{\alpha_i},
\end{equation}
For $n = 2k+1$, relation \eqref{eq_factors_Bn_1B} also has $n$ factors
and is equivalent to the following:
\begin{equation}
  \label{eq_product_Bn_odd}
  w_0 =  s_{\alpha_n}\prod\limits_{i = 1,3,...,n-2}
         s_{\varepsilon_i + \varepsilon_{i+1}}s_{\alpha_i},
\end{equation}
~\\
The orthogonality of roots corresponding to reflections in \eqref{eq_product_Bn_even}
(resp. \eqref{eq_product_Bn_odd}) is verified directly.
For example, since $(\varepsilon_i, \varepsilon_j) = 0$ for $i \neq j$, we have
\begin{equation*}
   (\varepsilon_i + \varepsilon_{i+1}, \varepsilon_{i+2} + \varepsilon_{i+3}) = 0, \; i + 3 \leq n.
\end{equation*}
In addition, for $i + 1 < n$, the following relation holds:
\begin{equation*}
   (\varepsilon_i + \varepsilon_{i+1}, \alpha_i) =
          (\alpha_i,  \alpha_i) + 2(\alpha_i,  \alpha_{i+1}) = 2 - 2 = 0, \\
\end{equation*}
and for $i + 3 \leq n$,
\begin{equation*}
   (\varepsilon_i + \varepsilon_{i+1}, \alpha_{i+2}) =
    2(\alpha_{i+1} + \alpha_{i+2}  + \alpha_{i+3}, \; \alpha_{i+2}) = 0. \qed
\end{equation*}

\section{\bf Decomposition in $W(C_n)$}

\begin{figure}[h]
\centering
   \includegraphics[scale=0.4]{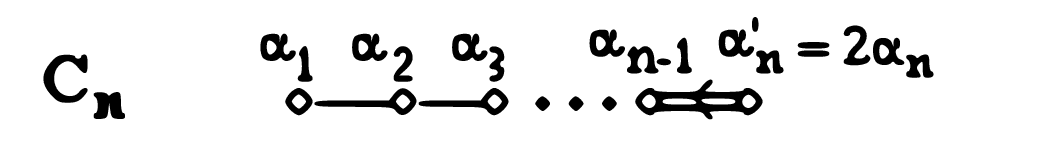}
\caption{Numbering of simple roots in $C_n$. Here, $\alpha_n$ (resp. $\alpha'_n$) 
is the simple root in $B_n$ (resp. $C_n$)}
%%%%%% The label must come after caption
\label{fig_Cn}
\end{figure}

We introduce the notation for the highest roots in the root subsystems $C_i \subset C_n$
see Table \ref{tab_factoriz_all}. The numbering of vertices of $C_n$ is shown
in Fig. \ref{fig_Cn}.
\begin{equation}
  \label{eq_not_Bn_Cn_2}
      \alpha^{ci}_{max} := 2\sum\limits_{j = n-i+1}^{n-1} \alpha_j + \alpha'_n,
      \text{ where } i \ge 2, \text{ and } \alpha^{c1}_{max} := \alpha'_n. \\
\end{equation}
see Remark \ref{rem_notations}. The root $\alpha^{ci}_{max}$ is the \underline{highest root}
in the root system $C_i$ for any $i \geq 2$.
~\\

Recall that the Weyl groups $W(B_n)$ and $W(C_n)$ coincide,
however the decompositions of the longest elements are specialized
for the cases $B_n$ and $C_n$. Consider the root systems
$B_n$ and $C_n$  in the same linear space $\mathbb{E}_n$ with
an orthogonal basis of unit vectors
\begin{equation*}
  \{ \varepsilon_i \mid i = 1,2, \dots n \}.
\end{equation*}
The simple roots in $B_n$ are denoted by $\alpha_i$, see \eqref{eq_compare_with_Bn}.
Denote simple roots in $C_n$ by $\alpha'_i$. The simple roots $\{ \alpha'_i \mid \ i =1,\dots, n  \}$ in $C_n$
are as follows:
\begin{equation*}
 \begin{split}
    &  \alpha'_i = \varepsilon_i - \varepsilon_{i+1}
       \text{ for } i = 1, \dots,n-1, \text{ and } \alpha'_n = 2\varepsilon_n, \\
    & \norm{\alpha_i} = \sqrt{2}  \text { for } i = 1,\dots,n-1, \text{ and }
       \norm{\alpha_n} = 2.
 \end{split}
\end{equation*}
We have
\begin{equation*}
  \alpha'_i = \alpha_i \text{ for } i < n, \text{ and } \alpha'_n = 2\alpha_n
\end{equation*}
Put $\varepsilon'_i = 2\varepsilon_i$ for $i = 1,\dots, n$. Then,
\begin{equation}
  \label{eq_def_ei_in_Bn_2}
  \begin{split}
   & \varepsilon_i = \alpha_n + \alpha_{n-1} + \dots + \alpha_i, \;\; \text{ for } i < n, \\
   &  \varepsilon'_i = \alpha'_n + 2\alpha_{n-1} + \dots + 2\alpha_i \;\; \text{ for } i < n, \\
   &  \varepsilon_n = \alpha_n, \;\; \varepsilon'_n = 2\varepsilon_n = 2\alpha_n = \alpha'_n.
  \end{split}
\end{equation}
By \eqref{eq_not_Bn_Cn_2} and \eqref{eq_def_ei_in_Bn_2}, $\varepsilon'_i$ is the highest root in $C_i$, i.e.,
\begin{equation*}
   \varepsilon'_i = \alpha^{c,n-i+1}_{max} \; \text{ for } \; i \geq 2, \text{ and }
   \varepsilon'_n = \alpha^{c1}_{max} = \alpha'_n.
\end{equation*}
~\\
\begin{proposition}[decomposition in $W(C_n)$]
  \label{prop_factoriz_Cn}
  The element $w_0$ is decomposed into $n$ factors as follows:
\begin{equation}
  \label{eq_factoriz_Cn_1}
   w_0 = \prod\limits_{i=1}^n  s_{\varepsilon'_i} = s_{\alpha'_n}\prod\limits_{i=1}^{n-1} s_{\alpha^{c,n-i+1}_{max}} \;\;.
\end{equation}
The roots corresponding to factors \eqref{eq_factoriz_Cn_1} are \underline{mutually orthogonal}.
Thus, the factors in \eqref{eq_factoriz_Cn_1} commute with each other.
The roots $\varepsilon'_i$  are ordered as follows:
\begin{equation}
 \label{eq_highest_roots_Cn}
 \begin{array}{lll}
    \alpha'_n < \alpha^{c2}_{max} < \alpha^{c3}_{max}  < \dots  < \alpha^{cn}_{max}\;.  \\
 \end{array}
\end{equation}
\end{proposition}

\PerfProof  Since the Weyl groups $W(B_n)=W(C_n)$ coincide, the longest element $w_0$ is the same for
$B_n$ and $C_n$. Recall that $\varepsilon'_i = 2\varepsilon_n$, 
i.e., reflections $s_{\varepsilon_i}$ and $s_{\varepsilon'_i}$ coincide for $i=1,\dots,n$.
Then, taking into account \eqref{eq_not_fin_factoriz_Bn}, we obtain
\begin{equation*}
 %% \label{eq_not_fin_factoriz_Cn}
   w_0 = \prod\limits_{i=1}^n s_{\varepsilon_i} =
      \prod\limits_{i=1}^n s_{\varepsilon'_i},
\end{equation*}
where all $s_{\varepsilon_i}$ (resp. $s_{\varepsilon'_i}$) commute with each other.
The element $w_0$ acts on $\varepsilon'_i$ as $s_{\varepsilon'_i}$:
\begin{equation*}
   w_0(\varepsilon_i) = s_{\varepsilon'_i}(\varepsilon'_i) = -\varepsilon'_i.
\end{equation*}
Thus,  $w_0 = -1$ on $\mathbb{E}_n$, see \cite[Plate III]{Bo02}.
Since the longest element is unique, then $w_0$ is the longest.   \qed

\section{\bf Decomposition in $W(D_n)$, \; $n \geq 4$}

\begin{figure}[h]
\centering
   \includegraphics[scale=0.4]{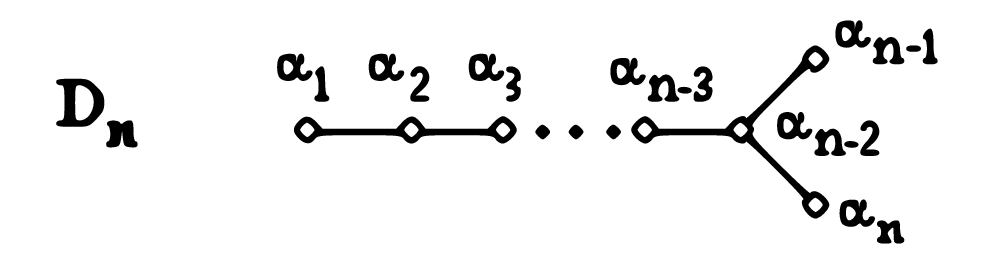}
\caption{\footnotesize{Numbering of simple roots in $D_n$}}
%%%%%% The label must come after caption
\label{fig_Dn}
\end{figure}

Let us introduce the notation for the highest roots in the root subsystem $D_i \subset D_n$,
see Table \ref{tab_factoriz_all}. The numbering of vertices in the Dynkin diagram $D_n$ is shown
in Fig. \ref{fig_Dn}.
\begin{equation}
  \label{eq_not_Dn_1}
 \begin{split}
     & \alpha^{di}_{max} :=
        \alpha_{n-i+1} + 2\sum\limits_{j = n-i+2}^{n-2}\alpha_j + \alpha_{n-1} + \alpha_n
           \; \text{ for } \; i \ge 4, \\
     & \alpha^{d3}_{max} := \alpha^{a3}_{max} := \alpha_{n-2} + \alpha_{n-1} + \alpha_n.
  \end{split}
\end{equation}
see Remark \ref{rem_notations}. The root $\alpha^{dp}_{max}$ is the \underline{highest root}
in the root system $D_p$ for any $p \geq 3$.
We shall regard the following element $w_0$ in $W(D_n)$:
\begin{equation}
  \label{eq_how_looks_w0_Dn}
   w_0 = s_n s_{n-1}(s_{n-2}s_n s_{n-1}s_{n-2}) \dots
       (s_1s_2 \dots s_{n-2}s_n s_{n-1}s_{n-2} \dots s_2s_1).
\end{equation}
In Proposition \ref{prop_1_Dn}, we find the max-orthogonal decomposition of $w_0$ and
prove that $w_0$ is the longest element in $W(D_n)$.
First of all, let us start with examples of decomposition of $w_0$ for small values of $n$.

\begin{example}
 \label{exmp_D4_1}
  \footnotesize
  \begin{equation*}
    \begin{split}
       n & =  4, \quad w_0 =  s_4s_3(s_2s_4s_3s_2)(s_1s_2s_4s_3s_2s_1) =
           s_4s_3s_2(s_4s_3(s_2s_1s_2)s_4s_3)s_2s_1 = \\
         &  s_4s_3s_2s_{\alpha_1 + \alpha_2 + \alpha_3 + \alpha_4}s_2s_1 =
           s_4s_3s_{\alpha_1 + 2\alpha_2 + \alpha_3 + \alpha_4}s_1 =
          \underline{ s_3s_4s_{\alpha_{max}}s_1}, \\
       n & = 5,  \quad w_0 =
           s_5s_4(s_3s_5s_4s_3)(s_2s_3s_5s_4s_3s_2)(s_1s_2s_3s_5s_4s_3s_2s_1) = \\
         & (s_5s_4s_3s_5s_4)s_3(s_2s_3s_5s_4(s_3s_2s_1s_2s_3)s_5s_4s_3s_2)s_1 = \\
         & s_{\alpha_3 + \alpha_4 + \alpha_5} s_3
            s_{\alpha_1 + 2\alpha_2 + 2\alpha_3 + \alpha_4 + \alpha_5}s_1 =
      \underline{(s_{\alpha_3 + \alpha_4 +\alpha_5})s_3(s_{\alpha_{max}})s_1}, \\
        n & = 6,  \quad w_0 =
           s_6s_5(s_4s_6s_5(s_4s_3s_4)s_6s_5s_4)s_3 \; \times \\
   & \quad (s_2s_3s_4(s_6s_5s_4(s_3s_2s_1s_2s_3)s_4s_6s_5)s_4s_3s_2)s_1 = \\
     & \quad  s_5s_6(s_{\alpha_3 + 2\alpha_4 + \alpha_5 + \alpha_6})s_3
  (s_{\alpha_1 + 2\alpha_2 +2\alpha_3 + 2\alpha_4 + \alpha_5 + \alpha_6})s_1 =  \\
     & \quad  \underline{s_5s_6(s_{\alpha_3 + 2\alpha_4 + \alpha_5 + \alpha_6})s_3
  (s_{\alpha_{max}})s_1},  \\
        n & = 7,  \quad w_0 =
       s_7s_6(s_5s_7s_6s_5)(s_4s_5s_7s_6s_5s_4)(s_3s_4s_5s_7s_6s_5s_4s_3) \times \\
     & (s_2s_3s_4s_5s_7s_6s_5s_4s_3s_2)(s_1s_2s_3s_4s_5s_7s_6s_5s_4s_3s_2s_1) = \\
    & \dots (s_2s_3s_4s_5
 (s_{\alpha_1 + \alpha_2 + \alpha_3 + \alpha_4 + \alpha_5 + \alpha_6 + \alpha_7})
        s_5s_4s_3s_2)s_1 = \\
    & s_7s_6(s_5s_7s_6s_5)(s_4s_5s_7s_6(s_5s_4s_3s_4s_5)s_7s_6s_5s_4)s_3
 (s_{\alpha_1 + 2\alpha_2 + 2\alpha_3 + 2\alpha_4 + 2\alpha_5 + \alpha_6 + \alpha_7})
     s_1 = \\
 & (s_7s_6s_5s_7s_6)s_5(s_{\alpha_3 + 2\alpha_4 + 2\alpha_5 + \alpha_6 +
  \alpha_7})s_3
  (s_{\alpha_1 + 2\alpha_2 + 2\alpha_3 + 2\alpha_4 + 2\alpha_5 + \alpha_6
 + \alpha_7})s_1 = \\
 & {(s_{\alpha_5 + \alpha_6+ \alpha_7})s_5(s_{\alpha_3 + 2\alpha_4 + 2\alpha_5 + \alpha_6 +
  \alpha_7})s_3
  (s_{\alpha_1 + 2\alpha_2 + 2\alpha_3 + 2\alpha_4 + 2\alpha_5 + \alpha_6
 + \alpha_7})s_1} = \\
 & \underline{{(s_{\alpha_5 + \alpha_6+ \alpha_7})s_5(s_{\alpha_3 + 2\alpha_4 + 2\alpha_5 + \alpha_6 +
  \alpha_7})s_3
  (s_{\alpha_{max}})s_1}}. \\
   \end{split}
  \end{equation*}
\end{example}
Example \ref{exmp_D4_1} prompts us the general form of the element $w_0$.
\qed
~\\

\begin{proposition}[decomposition in $W(D_n)$]
 \label{prop_1_Dn}
{\rm{(i)}}   The element $w_0$ in $W(D_n)$ can be decomposed into $n$ (resp. $n-1$) factors
for $n=2k$ (resp. $n= 2k+1$): 
\begin{equation}
  \label{eq_D4_factor_even}
  w_0 =  
   \Biggl \{ 
     \begin{array}{lll}
       s_{\alpha_n}s_{\alpha_{n-1}}
       & \prod\limits_{\stackrel{i=1}{i-odd}}^{n-3}(s_{\omega_i}s_{\alpha_i}), \;
          & \text{ where }\; \omega_i = \alpha^{d,n-i+1}_{max}.  \\
          \\
       & \prod\limits_{\stackrel{i=1}{i-odd}}^{n-2}(s_{\omega_i}s_{\alpha_i}), \;
          & \text{ where }\; \omega_i = \alpha^{d,n-i+1}_{max}. \\
     \end{array} 
\end{equation}
~\\
The root corresponding to factors in \eqref{eq_D4_factor_even} are as follows:
\begin{equation}
  \label{eq_D4_max_roots_1}
 \footnotesize
\begin{split}
 \text{ For } n = 2k: \quad & 
  \begin{cases}
     & \alpha^{d4}_{max}, \; \alpha^{d6}_{max}, \; \dots, \;
       \alpha^{d,n-2}_{max}, \;\alpha_{max}, \\
     & \alpha_1, \alpha_3, \dots, \alpha_{n-5}, \alpha_{n-3}, \alpha_{n-1}, \alpha_{n}. \\
   \end{cases} \\
 & \\
  \text{ For } n = 2k+1: \quad  & 
  \begin{cases}
     & \alpha^{d3}_{max}, \; \alpha^{d5}_{max}, \; \alpha^{d7}_{max}, \; \dots, \; 
       \alpha^{d,n-2}_{max}, \;  \alpha_{max},\\
     & \alpha_1, \alpha_3, \dots, \alpha_{n-4}, \alpha_{n-2}.
   \end{cases} 
\end{split}
\end{equation}
~\\

{\rm{(ii)}}  The roots in the decomposition \eqref{eq_D4_factor_even}
are \underline{mutually orthogonal}. Thus, factors in \eqref{eq_D4_factor_even} commute with each other. 
The roots $\alpha^{dp}_{max}$  are ordered as follows:

\begin{equation}
 \label{eq_highest_roots_Dn}
 \begin{array}{lll}
   &  \alpha^{d4}_{max} < \alpha^{d6}_{max}  < \dots  < \alpha^{d,n-2}_{max} < \alpha^{dn}_{max} = \alpha_{max} &  \text{ for } n = 2k,  \\
   &  \alpha^{d3}_{max} < \alpha^{d5}_{max}  < \dots < \alpha^{d,n-2}_{max}  < \alpha^{dn}_{max} = \alpha_{max} & \text{ for } n = 2k+1.  \\
 \end{array}
\end{equation}
Here, $\alpha^{d3}_{max} = \alpha^{a3}_{max}$, see \eqref{eq_not_Dn_1}.
~\\

{\rm{(iii)}}
 The element $w_0$ acts as $-1$ if $n=2k$, and as $-\varepsilon$ if $n=2k+1$,
where $\varepsilon$ is the automorphism that interchanges roots
$\alpha_n$ and $\alpha_{n-1}$ and fixes all other simple roots.
In other words, is the longest element in $W(D_n)$,  see  \cite[Plate IV]{Bo02}.
\end{proposition}

\begin{figure}[h]
 \centering
   \includegraphics[scale=0.32]{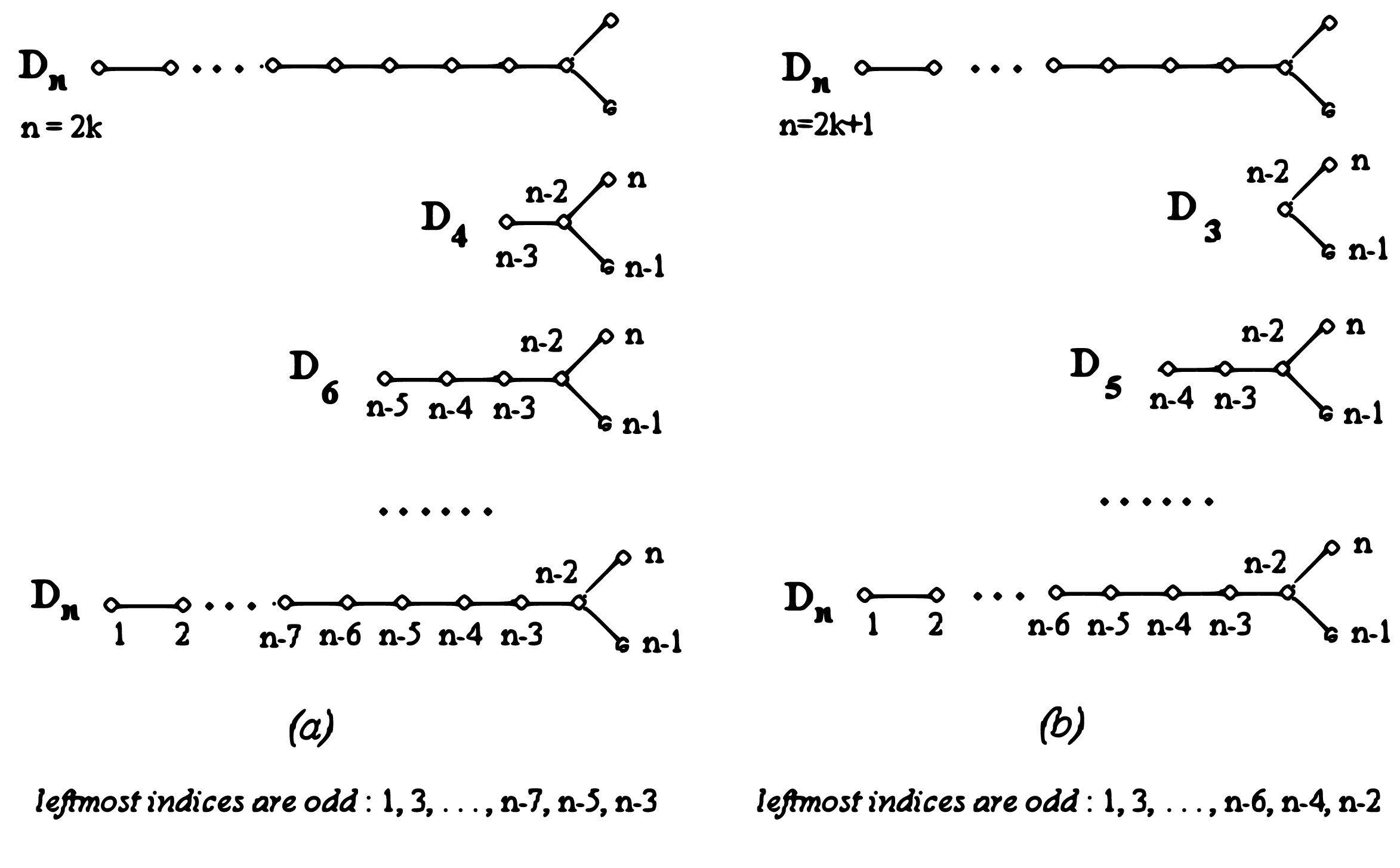}
\caption{\footnotesize The leftmost indices for
diagrams $D_{2p}$ (side ({\it a})) and $D_{2p+1}$ (side ({\it b})) are odd}
%%%%%% The label must come after caption
\label{fig_Dn_max_alpha}
\end{figure}

Before proving Proposition \ref{prop_1_Dn}, we need the following lemma:

\begin{lemma}
 \label{lem_only_2_roots}
{\rm{(i)}}  Let $D_k \subset D_n$ be the root system corresponding to the highest root
$\alpha^{dk}_{max}$ in \eqref{eq_highest_roots_Dn}.
The index of the leftmost root, see Fig. \ref{fig_Dn_max_alpha}, for the root system $D_k$ is odd and it
is equal to $n-k+1$.

{\rm{(ii)}} Let $i$ be even.
The roots $\alpha^{d,n-i}_{max}$ and $\alpha^{d,n-i+2}_{max}$
are the only $2$ roots among the highest roots $\alpha^{dp}_{max}$ in \eqref{eq_highest_roots_Dn}
that are not orthogonal to $\alpha_i$.

{\rm{(iii)}} Let $i$ be odd. All roots $\alpha^{dp}_{max}$ in \eqref{eq_highest_roots_Dn} 
are orthogonal to $\alpha_i$.

\end{lemma}

\PerfProof (i) In the list of \eqref{eq_D4_max_roots_1} with $n=2k$,
the leftmost root of $D_{2p}$ is $n-2p+1$, i.e., all leftmost roots have odd indices,
see Fig. \ref{fig_Dn_max_alpha}(a).
In the list of \eqref{eq_D4_max_roots_1} with $n=2k+1$,
the leftmost of $D_{2p+1}$ is $n-2p$, i.e., again, all leftmost roots have odd indices,
see Fig. \ref{fig_Dn_max_alpha}(b).

(ii) Adjacent roots for $\alpha_i$, where $i$ is even, are
$\alpha_{i-1}$ and $\alpha_{i+1}$. By (i) they correspond to
the leftmost roots of $D_{n-i+2}$ and $D_{n-i}$.
Among all the highest roots of the list \eqref{eq_highest_roots_Dn}
only the highest roots corresponding to $D_{n-i+2}$
and $D_{n-i}$ are not orthogonal to $\alpha_i$. 
These roots are $\alpha^{d,n-i}_{max}$ and $\alpha^{d,n-i+2}_{max}$.

(iii) It holds because the only simple root not orthogonal 
to $\alpha^{dp}_{max}$ is the root adjacent to the leftmost root with index $n-p+1$, 
but such a root has even index $n-p+2$. %% since n-p is always even
\qed
~\\

\pprof{Proposition \ref{prop_1_Dn}}

(i) Let us divide all the factors in brackets in \eqref{eq_how_looks_w0_Dn} 
into pairs of neighboring factors. Each such pair is converted into 
two neighboring factors in \eqref{eq_D4_factor_even}.
Taking into account that $s_n$ and $s_{n-1}$ commute and repeatedly applying 
Proposition \ref{prop_fact_1}, we get:
\begin{equation*}
 \footnotesize
  \begin{split}
  & (s_{i+1} \dots s_{n-2}s_n s_{n-1}s_{n-2} \dots s_{i+1})
           (s_is_{i+1} \dots s_{n-2}s_n s_{n-1}s_{n-2} \dots s_{i+1}s_i) =  \\
  & s_{i+1} \dots s_{n-2}(s_n s_{n-1}s_{n-2} \dots s_{i+1}
           s_i s_{i+1} \dots s_{n-2} s_{n-1}s_n)s_{n-2} \dots s_{i+1}s_i =  \\
  & s_{i+1} \dots s_{n-3}s_{n-2}(s_{\alpha_n +\alpha_{n-1} + \alpha_{n-2} +
      \alpha_{n-3} + \dots \alpha_{i+1} +
    \alpha_i} ) s_{n-2}s_{n-3} \dots s_{i+1}s_i =  \\
  & s_{i+1} \dots s_{n-3}(s_{\alpha_n +\alpha_{n-1} + 2\alpha_{n-2} +
      \alpha_{n-3} + \dots \alpha_{i+1} +
    \alpha_i} ) s_{n-3} \dots s_{i+1}s_i =  \\
  & s_{i+1} \dots (s_{\alpha_n +\alpha_{n-1} + 2\alpha_{n-2} +
      2\alpha_{n-3} + \dots \alpha_{i+1} +
    \alpha_i} )  \dots s_{i+1}s_i = \dots \\
  & (s_{\alpha_n +\alpha_{n-1} + 2\alpha_{n-2} +
      2\alpha_{n-3} + \dots + 2\alpha_{i+1} + \alpha_i} )s_i. \\
  \end{split}
\end{equation*}

(ii)  The only simple root not orthogonal to $\alpha_{max}$ is
$\alpha_2$, which is not a component in any root of list \eqref{eq_D4_max_roots_1}. 
Then, the highest root $\alpha_{max} = \alpha^{dn}_{max}$ is orthogonal to all other roots 
of list \eqref{eq_D4_max_roots_1}. Further, the only simple root not orthogonal to 
$\alpha^{d,n-2}_{max}$ is $\alpha_4$, which is not a component in any roots of list
\eqref{eq_D4_max_roots_1}, except of $\alpha_{max}$. Then, the root $\alpha^{d,n-2}_{max}$
is orthogonal to all other roots of list \eqref{eq_D4_max_roots_1}, except of $\alpha_{max}$.
But orthogonality $\alpha_{max}$ and $\alpha^{d,n-2}_{max}$ has already been proved.  
Continuing in this way, we obtain orthogonality for all remaining roots in \eqref{eq_D4_max_roots_1}.
~\\

(iii)
We consider the action of $w_0$ on the following $4$ types of simple roots $\alpha_i$ in $D_n$.

  (a) \underline{Actions on $\alpha_i$ with odd $i < n-2$}.
   By Lemma \ref{lem_only_2_roots}(iii) among the factors associated with
   \eqref{eq_D4_max_roots_1}  
   only $s_{\alpha_i}$ acts not trivially on $\alpha_i$, i.e.,
 \begin{equation*}
   \begin{split}
      & w_0(\alpha_i) = s_{\alpha_i}(\alpha_i) = -\alpha_i.
    \end{split}
 \end{equation*}

 (b) \underline{Actions on $\alpha_i$ with even $i < n-2$}.
  In this section, we use a shorter notation:  
\begin{equation*}
   s_{n-i,  max} \text{ instead of } s_{\alpha^{d,n-i}_{max}} \text{ and } 
   s_{n-i+2,  max} \text{ instead of } s_{\alpha^{d,n-i+2}_{max}}.
\end{equation*}  
 By Lemma \ref{lem_only_2_roots}, among roots \eqref{eq_D4_max_roots_1}
 only the following reflections non-trivially act on $\alpha_i$:
\begin{equation*}
   s_{\alpha_{i-1}}, s_{\alpha_{i+1}}, s_{n-i+2, max}, s_{n-i,  max}
\end{equation*}
Using eq. \eqref{eq_not_Dn_1} we get the following:
\begin{equation*}
 \begin{split}
   & \alpha^{d,n-i} = \alpha_{i+1} + 2(\alpha_{i+2} + \dots \alpha_{n-2}) + \alpha_{n-1} + \alpha_n, \\
   & \alpha^{d,n-i+2} = \alpha_{i-1} + 2(\alpha_i + \dots \alpha_{n-2}) + \alpha_{n-1} + \alpha_n.
 \end{split}
\end{equation*}
Then,
\begin{equation}
  \label{eq_Dn_max_alpha}
 \begin{split}
 &  (\alpha^{d,n-i}_{max}, \alpha_i) = (\alpha_{i+1}, \alpha_i) = -1, \\
 & (\alpha^{d,n-i+2}_{max}, \alpha_i) = (\alpha_{i-1} + 2\alpha_i + 2\alpha_{i+1}, \alpha_i)
    = -1+4-2 = 1, \\
 &  \alpha^{d,n-i+2}_{max} - \alpha^{d,n-i}_{max} =  \alpha_{i-1} + 2\alpha_i + \alpha_{i+1}.
 \end{split}
\end{equation}
see Fig. \ref{fig_Dn_max_alpha}. By \eqref{eq_Dn_max_alpha} we have
\begin{equation*}
  \begin{split}
    w_0(\alpha_i) = & s_{n-i, max} s_{n-i+2,m} s_{i-1} s_{i+1} (\alpha_i) = \\
    &  s_{n-i, max}s_{n-i+2, max}  (\alpha_i + \alpha_{i-1} + \alpha_{i+1}).
  \end{split}
\end{equation*}
Since $i+1$ and $i-1$ are odd, by Lemma \ref{lem_only_2_roots}(iii)
reflections $s_{n-i, max}$ and $s_{n-i+2, max}$ preserve $\alpha_{i-1} + \alpha_{i+1}$.
Then, by \eqref{eq_Dn_max_alpha}
\begin{equation*}
  \begin{split}
   w_0(\alpha_i) =
    &  \alpha_{i-1} + \alpha_{i+1} +  s_{n-i, max}s_{n-i+2, max}(\alpha_i) = \\
    & \alpha_{i-1} + \alpha_{i+1} + s_{n-i, max}(\alpha_i - \alpha^{d,n-i+2}_{max}) = \\
    & \alpha_{i-1} + \alpha_{i+1} + \alpha_i + (\alpha^{d,n-i}_{max} - \alpha_{max}^{n-i+2}) = \\
   & \alpha_{i-1} + \alpha_i + \alpha_{i+1}  -
                    (\alpha_{i-1} + 2\alpha_i + \alpha_{i+1}) = -\alpha_i. \\
  \end{split}
\end{equation*}

(c) \underline{Actions on $\alpha_{n-1}$ and $\alpha_{n}$}.
Let $n$ be even.  Among the factors associated with \eqref{eq_D4_max_roots_1} 
only $s_{\alpha_{n-1}}$ acts on $\alpha_{n-1}$, and only $s_{\alpha_{n}}$ acts on
$\alpha_{n}$. Thus,
 \begin{equation*}
   \begin{split}
      & w_0(s_{\alpha_{n-1}}) = -s_{\alpha_{n-1}}, \quad
      w_0(s_{\alpha_{n}}) = -s_{\alpha_{n}}.
    \end{split}
 \end{equation*}
 Let $n$ be odd.  Among the factors associated with \eqref{eq_D4_max_roots_1} only
 $s_{\alpha_{n-2}}$ and $\alpha^{d3}_{max}$  acts on $\alpha_{n-1}$, $\alpha_{n}$.
 Then,
 \begin{equation*}
   \begin{split}
      w_0(\alpha_{n-1}) = & s_{\alpha^{d3}_{max}}s_{n-2}(\alpha_{n-1}) =
         s_{\alpha^{d3}_{max}}(\alpha_{n-1} + \alpha_{n-2}) = \\
      & s_{\alpha_{n} + \alpha_{n-1} + \alpha_{n-2}}(\alpha_{n-1} + \alpha_{n-2}) =
        -\alpha_{n}, \text{ and }\\
       w_0(\alpha_{n}) = & s_{\alpha^{d3}_{max}}s_{n-2}(\alpha_{n}) =
         s_{\alpha^{d3}_{max}}(\alpha_{n} + \alpha_{n-2}) = \\
      & s_{\alpha_{n} + \alpha_{n-1} + \alpha_{n-2}}(\alpha_{n} + \alpha_{n-2}) =
        -\alpha_{n-1}.
    \end{split}
 \end{equation*}
   Thus, if $n$ is odd, $w_0= -\varepsilon$, where $\varepsilon$ acts on
   as the automorphism that interchanges $\alpha_{n}$ and $\alpha_{n-1}$, see \cite[Plate IV]{Bo02}.
~\\

   (d) \underline{Actions on the branch point $\alpha_{n-2}$}.
   Let $n$ is even.
   Among the factors associated with \eqref{eq_D4_max_roots_1} only
   $s_{\alpha_{n}}$,  $s_{\alpha_{n-1}}$, $s_{\alpha_{n-3}}$ and
   $s_{\alpha^{d4}_{max}}$ act on  $\alpha_{n-2}$.  Note that
   $\alpha^{d4}_{max} \perp \{\alpha_{n-1}, \alpha_{n}, \alpha_{n-3}\}$, then
 \begin{equation*}
   \begin{split}
     w_0(\alpha_{n-2}) = &
     s_{\alpha^{d4}_{max}}s_{\alpha_{n}}s_{\alpha_{n-1}}s_{\alpha_{n-3}}(\alpha_{n-2}) = \\
    & s_{\alpha^{d4}_{max}}(\alpha_{n-2} + \alpha_{n-3} + \alpha_{n-1} + \alpha_{n}) =  \\
    & \alpha_{n-3} + \alpha_{n-1} + \alpha_{n} + s_{\alpha^{d4}_{max}}(\alpha_{n-2})= \\
    &  \alpha_{n-3} + \alpha_{n-1} + \alpha_{n} + (\alpha_{n-2} - \alpha^{d4}_{max}) =
      -\alpha_{n-2}.
    \end{split}
 \end{equation*}
   Let $n$ is odd.
   Among the factors of \eqref{eq_D4_max_roots_1} only
   $s_{\alpha_{n-2}}$ acts on $\alpha_{n-2}$. Then,
   $w_0(\alpha_{n-2}) = s_{\alpha_{n-2}}(\alpha_{n-2}) = -\alpha_{n-2}$.  \qed
~\\

\section{\bf Decomposition in $W(E_6)$}

Let us denote the highest roots in the root subsystems $A_3 \subset A_5 \subset E_6$ as follows:
\begin{equation}
  \label{eq_not_E6_1}
  \begin{split}
     & \alpha^{a3}_{max} := \alpha_3 + \alpha_4 + \alpha_5, \\
     & \alpha^{a5}_{max} := \alpha_1 +  \alpha_3 + \alpha_4 + \alpha_5 + \alpha_6,
   \end{split}
\end{equation}
see Table \ref{tab_factoriz_all} and Remark \ref{rem_notations}.
Note that indices of simple roots for $\alpha^{a5}_{max}$ in \eqref{eq_not_E6_1} differ from the indices
of simple roots that appears as summands in $\alpha^{a5}_{max}$ 
defined in \eqref{eq_def_An}, see Remark \ref{rem_differ_indices}. 

\begin{figure}[h]
\centering
   \includegraphics[scale=0.4]{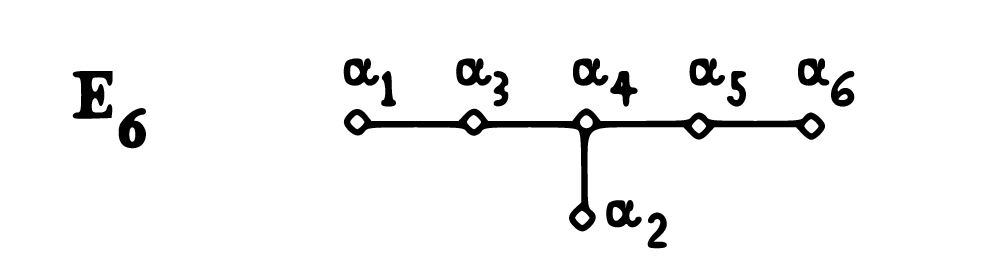}
\caption{Numbering of simple roots in $E_6$}
%%%%%% The label must come after caption
\label{fig_E6}
\end{figure}

The numbering of vertices of the Dynkin diagram $E_6$ is shown
in Fig. \ref{fig_E6}.
Consider the following two elements in the Weyl group $W(E_6)$:

\begin{equation*}
 \label{eq_w1_w2_for_E6}
  \begin{split}
   & w_1 = s_2 s_4 s_5 s_3 s_4 (s_6 s_5 s_2 (s_4 s_3 s_1 s_3 s_4) s_2 s_5 s_6) s_4 s_3 s_5 s_4 s_2, \\
   & w_2 = s_6 (s_5 s_6) (s_4 s_5 s_6) (s_3 s_4 s_5 s_6) (s_1 s_3 s_4 s_5 s_6). \\
  \end{split}
\end{equation*}

\begin{proposition}[decomposition in $W(E_6)$]
 \label{prop_factoriz_E6}
  {\rm{(i)}}
   The element $w_1$ is the reflection associated with the highest root
   $\alpha_{max}$:
\begin{equation}
   \label{eq_how_w1_looks_E6}
   w_1 = s_{\alpha_{max}}.
\end{equation}

  {\rm{(ii)}}
   The element $w_2$ is the product of $3$ reflections associated with positive roots
 $\alpha_4$, $\alpha^{a3}_{max}$, $\alpha^{a5}_{max}$:
\begin{equation}
  \label{eq_for_n_is_5}
  w_2 = s_{\alpha_4} s_{\alpha^{a3}_{max}} s_{\alpha^{a5}_{max}}.
\end{equation}

  {\rm{(iii)}}
   The element $w_0 = w_1w_2$ is the product of $4$ reflections associated with
   \underline{mutually orthogonal} roots
\begin{equation}
 \label{eq_4roots_E6}
   \alpha_4, \quad \alpha^{a3}_{max}, \quad \alpha^{a5}_{max}, \quad \alpha_{max}.
\end{equation}
   The roots \eqref{eq_4roots_E6} are ordered as follows:
\begin{equation*}
    \alpha_4 < \alpha^{a3}_{max} < \alpha^{a5}_{max} < \alpha^{e6}_{max}.  \\
\end{equation*}
{\rm{(iv)}}
   The element $w_0$ is the longest element in $W(E_6)$.
\end{proposition}

\PerfProof (i) We consequentially apply Proposition \ref{prop_fact_1} from $\S$\ref{sec_app_1}:

\begin{equation*}
  \begin{split}
   & s_3 s_1 s_3 = s_{\alpha_1 + \alpha_3}, \\
   & s_4 s_3 s_1 s_3 s_4 = s_4 s_{\alpha_1 + \alpha_3} s_4 =
                    s_{\alpha_1 + \alpha_3 + \alpha_4}, \\
   & s_2 s_4 s_3 s_1 s_3 s_4 s_2 = s_2 s_{\alpha_1 + \alpha_3 + \alpha_4} s_2 =
                    s_{\alpha_1 + \alpha_3 + \alpha_4 + \alpha_2}, \\
   & s_5 s_2 s_4 s_3 s_1 s_3 s_4 s_2 s_5 =
                    s_5  s_{\alpha_1 + \alpha_3 + \alpha_4 + \alpha_2} s_5 =
                    s_{\alpha_1 + \alpha_3 + \alpha_4 + \alpha_2 + \alpha_5}, \\
   & s_6 s_5 s_2 s_4 s_3 s_1 s_3 s_4 s_2 s_5 s_6 =
                    s_{\alpha_1 + \alpha_3 + \alpha_4 + \alpha_2 + \alpha_5 + \alpha_6}, \\
   & s_4s_6 s_5 s_2 s_4 s_3 s_1 s_3 s_4 s_2 s_5 s_6s_4 =
                    s_{\alpha_1 + \alpha_3 + 2\alpha_4 + \alpha_2 + \alpha_5 + \alpha_6}, \\
   & s_3s_4s_6 s_5 s_2 s_4 s_3 s_1 s_3 s_4 s_2 s_5 s_6s_4s_3 =
                    s_{\alpha_1 + 2\alpha_3 + 2\alpha_4 + \alpha_2 + \alpha_5 + \alpha_6}, \\
   & s_5s_3s_4s_6 s_5 s_2 s_4 s_3 s_1 s_3 s_4 s_2 s_5 s_6s_4s_3s_5 =
                    s_{\alpha_1 + 2\alpha_3 + 2\alpha_4 + \alpha_2 + 2\alpha_5 + \alpha_6}, \\
   & s_4s_5s_3s_4s_6 s_5 s_2 s_4 s_3 s_1 s_3 s_4 s_2 s_5 s_6s_4s_3s_5s_4 =
                    s_{\alpha_1 + 2\alpha_3 + 3\alpha_4 + \alpha_2 + 2\alpha_5 + \alpha_6}, \\
   & s_2s_4s_5s_3s_4s_6 s_5 s_2 s_4 s_3 s_1 s_3 s_4 s_2 s_5
     s_6s_4s_3s_5s_4s_2 =
 s_{\alpha_1 + 2\alpha_3 + 3\alpha_4 + 2\alpha_2 + 2\alpha_5 + \alpha_6}. \\
  \end{split}
\end{equation*}
The last reflection corresponds to $\alpha_{max}$.
~\\

(ii) Applying  Proposition \ref{prop_factoriz_An}(ii)
and eq.\eqref{eq_factors_An_1} from Proposition \ref{prop_factoriz_An_2} to the sequence
of reflections $\{s_1, s_3, s_4, s_5, s_6\}$, we get
\begin{equation*}
 %% w_2 = s_{\alpha_4} s_{\scriptscriptstyle {a3}, max} s_{\scriptscriptstyle {a5}, max}.
  w_2 = s_{\alpha_4} s_{\alpha^{a3}_{max}} s_{\alpha^{a5}_{max}}.
\end{equation*}

(iii) Since $\alpha_{max}$ is orthogonal to all simple roots $\alpha_i$ except $\alpha_2$, then $\alpha_{max}$ is orthogonal to $\alpha_4$, $\alpha^{a3}_{max}$, and $\alpha^{a5}_{max}$, see \eqref{eq_not_E6_1}.
The orthogonality relations for roots $\alpha_4$, $\alpha^{a3}_{max}$
and $\alpha^{a5}_{max}$  are verified directly.
~\\

(iv) It suffices to prove that $w_0$ acts on simple roots in the same way as the longest element.
The reflection $s_{\alpha_{max}}$ acts only on the simple root $\alpha_2$,
the reflection $\alpha^{a5}_{max}$ acts only
on simple roots $\alpha_1$, $\alpha_2$, $\alpha_6$.
\begin{equation*}
  \begin{split}
     & s_{\alpha^{a5}_{max}} s_{\alpha^{a3}_{max}}s_{\alpha_4}(\alpha_1) =
     s_{\alpha^{a5}_{max}}(\alpha_1 + \alpha_3 + \alpha_4 + \alpha_5) =  -\alpha_6, \\
     & s_{\alpha^{a5}_{max}}s_{\alpha^{a3}_{max}}s_{\alpha_4}(\alpha_3) =
        s_{\alpha^{a5}_{max}} s_{\alpha^{a3}_{max}}(\alpha_3 + \alpha_4) =
          s_{\alpha^{a5}_{max}}(-\alpha_5) = -\alpha_5, \\
     & s_{\alpha^{a5}_{max}} s_{\alpha_3 + \alpha_4 + \alpha_5}s_{\alpha_4}(\alpha_4) =
          s_{\alpha^{a5}_{max}}(-\alpha_4) = -\alpha_4, \\
  \end{split}
\end{equation*}
By (iii) $w_0$ is an involution, so we do not have to check the action of $w_0$
on $\alpha_5$ and $\alpha_6$. It remains only to check the action of $w_0$ on $\alpha_2$:
\begin{equation*}
      s_{\alpha^{a3}_{max}}s_4(\alpha_2) =
        s_{\alpha^{a3}_{max}}(\alpha_2 + \alpha_4) =
        \alpha_2 + \alpha_3 + 2\alpha_4 + \alpha_5 \\
\end{equation*}
It is easy to check that  $(\alpha^{a5}_{max}\; , \; \alpha_2 + \alpha_3 + 2\alpha_4 + \alpha_5) = -1$, then
\begin{equation*}
      s_{\alpha^{a5}_{max}}(\alpha_2 + \alpha_3 + 2\alpha_4 + \alpha_5) = \alpha^{a5}_{max} + \alpha_2 +
        \alpha_3 + 2\alpha_4 + \alpha_5 = \alpha_{max} - \alpha_2.
\end{equation*}
Further, since $(\alpha_{max}, \alpha_2) = 1$, we get
\begin{equation*}
  \begin{split}
   & s_{\alpha_{max}}(\alpha_{max} - \alpha_2) = -\alpha_{max} - s_{\alpha_{max}}(\alpha_2) = \\
   & -\alpha_{max} - (\alpha_2 -(\alpha_{max}, \alpha_2)\alpha_{max}) = -\alpha_2.
  \end{split}
\end{equation*}
Thus, $w_0$ maps
\begin{equation*}
  \alpha_1,  \alpha_2, \alpha_3, \alpha_4, \alpha_5, \alpha_6 \;\; \text{ to } \;
  -\alpha_6, -\alpha_2, -\alpha_5, -\alpha_4, -\alpha_3, -\alpha_1,
\end{equation*}
see \cite[Plate V]{Bo02}. Since the longest element is unique in $W(E_6)$, %% Bo02
$w_0$ is the longest.
\qed

\section{\bf Decomposition in $W(E_7)$}
  \label{sec_E7}

\begin{figure}[h]
\centering
   \includegraphics[scale=0.4]{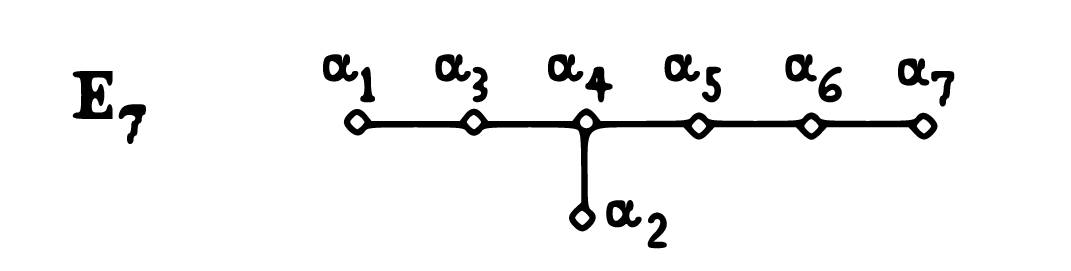}
\caption{Numbering of simple roots in $E_7$}
%%%%%% The label must come after caption
\label{fig_E7}
\end{figure}

Let us denote the \underline{highest roots} in the root subsystems $D_4 \subset D_6 \subset E_7$ as follows:
\begin{equation}
  \label{eq_not_E7_1}
  \begin{split}
     & \alpha^{d4}_{max} := \alpha_2 + \alpha_3 + 2\alpha_4 + \alpha_5, \\
     & \alpha^{d6}_{max} := \alpha_2 + \alpha_3 + 2\alpha_4 + 2\alpha_5 + 2\alpha_6 + \alpha_7,
   \end{split}
\end{equation}
see Table \ref{tab_factoriz_all} and Remark \ref{rem_notations}.
In \eqref{eq_not_E7_1}, the indices of simple roots  differ from the indices
of simple roots that appears as summands in $\alpha^{d4}_{max}$ and $\alpha^{d6}_{max}$
defined in \eqref{eq_not_Dn_1}, see Remark \ref{rem_differ_indices}. 
The numbering of vertices of the Dynkin diagram $E_7$ is shown
in Fig. \ref{fig_E7}.
We construct the longest element $w_0$ as the product of two following
elements in the Weyl group $W(E_7)$:
\begin{equation}
  \begin{split}
    w_1 =
 & s_7 (s_6 s_7) (s_5 s_6 s_7) (s_4 s_5 s_6 s_7) (s_3 s_4 s_5 s_6 s_7) (s_2 s_4 s_5 s_6 s_7)
    (s_3 s_4 s_5 s_6) \times \\
        & (s_2 s_4 s_5 s_3 s_4 s_2), \\
    w_2 =
   & s_1 s_3 s_4 (s_5 s_2 s_4 s_6 (s_5 s_7 s_6 s_3 (s_4 s_5 s_2 (s_4 s_3 s_1 s_3 s_4) s_2 s_5 s_4)
        s_3 s_6 s_7 s_5) \times \\
          & s_6 s_4 s_2 s_5 ) s_4 s_3 s_1. \\
  \end{split}
\end{equation}

\begin{proposition}
 \label{prop_6_orthog}
  {\rm{(i)}} The element $w_1$ is the product of $6$ reflections associated with
  the following mutually orthogonal roots:
\begin{equation}
  \label{eq_6_roots_of_E7}
    \alpha_4, \;\; \alpha_7, \;\; \alpha^{d6}_{max}, \;\;
    \alpha_3 + \alpha_4 + \alpha_5, \;\; \alpha_2 + \alpha_4 + \alpha_5.
         \;\; \alpha_2 + \alpha_4 + \alpha_3,
\end{equation}

  {\rm{(ii)}} The element $w_2$ is the reflection associated with the highest root:
\begin{equation*}
            w_2 = s_{\alpha_{max}}.
\end{equation*}
\end{proposition}

\PerfProof
(i) Let us denote by $w^{(1)}_1$, $w^{(2)}_1$ $w^{(3)}_1$, $w^{(4)}_1$ parts of the element $w_1$:
\begin{equation*}
  \begin{split}
    & w^{(1)}_1  :=  s_7 (s_6 s_7)(s_5 s_6 s_7)(s_4 s_5 s_6 s_7)(s_3 s_4 s_5 s_6 s_7), \\
    & w^{(2)}_1  :=  s_2 s_4 s_5 s_6 s_7, \\
    & w^{(3)}_1  :=  s_3 s_4 s_5 s_6, \\
    & w^{(4)}_1  :=  s_2 s_4 s_5 s_3 s_4 s_2, \\
  \end{split}
\end{equation*}
Then,
\begin{equation*}
   w_1 = w^{(1)}_1 w^{(2)}_1 w^{(3)}_1 w^{(4)}_1.
\end{equation*}
Using Proposition \ref{prop_factoriz_An},  eq. \eqref{prop_lst_in_An_2} and 
Proposition \ref{eq_factors_An_1}, eq. \eqref{eq_factors_An_1}, we get
\begin{equation*}
  \begin{split}
   w^{(1)}_1  =  & s_7 (s_6 s_7)(s_5 s_6 s_7)(s_4 s_5 s_6 s_7)(s_3 s_4 s_5 s_6 s_7) = \\
   & s_{\alpha_5} s_{\alpha_4+\alpha_5+\alpha_6} s_{\alpha_3 + \alpha_4+\alpha_5+\alpha_6 + \alpha_7} \, .
  \end{split}
\end{equation*}
~\\
Then, repeatedly applying Proposition \ref{prop_fact_1}, we obtain
\begin{equation*}
  \begin{split}
   w^{(1)}_1w^{(2)}_1  = & s_{\alpha_5} s_{\alpha_4+\alpha_5+\alpha_6}
       s_{\alpha_3 + \alpha_4+\alpha_5+\alpha_6 + \alpha_7} (s_2 s_4 s_5 s_6 s_7) = \\
     & s_{\alpha_5} s_{\alpha_4+\alpha_5+\alpha_6} (s_2 s_4 s_5 s_6 s_7)
       s_{\alpha_2 + \alpha_3 + 2\alpha_4+2\alpha_5+2\alpha_6 + \alpha_7} = \\
     & s_{\alpha_5}s_{\alpha_4}s_{\alpha_5} (s_6 s_5 s_4 s_2 s_4 s_5 s_6) s_7
            s_{\alpha_2 + \alpha_3 + 2\alpha_4+2\alpha_5+2\alpha_6 + \alpha_7} = \\
     &  s_{\alpha_4+\alpha_5} s_{\alpha_2+\alpha_4+\alpha_5+\alpha_6} s_{\alpha_7}
            s_{\alpha_2 + \alpha_3 + 2\alpha_4+2\alpha_5+2\alpha_6 + \alpha_7} = \\
     &  s_{\alpha_7}s_{\alpha_4+\alpha_5} s_{\alpha_2+\alpha_4+\alpha_5+\alpha_6+\alpha_7}
            s_{\alpha_2 + \alpha_3 + 2\alpha_4+2\alpha_5+2\alpha_6 + \alpha_7} \,.
  \end{split}
\end{equation*}
~\\
Further, by \eqref{eq_not_E7_1}
\begin{equation*}
  \begin{split}
    w^{(1)}_1 & w^{(2)}_1 w^{(3)}_1 = \\
      & s_{\alpha_7}s_{\alpha_4+\alpha_5} s_{\alpha_2+\alpha_4+\alpha_5+\alpha_6+\alpha_7}
       s_{\alpha_2 + \alpha_3 + 2\alpha_4+2\alpha_5+2\alpha_6 + \alpha_7}(s_3 s_4 s_5 s_6) = \\
      & s_{\alpha_7}s_{\alpha_4+\alpha_5} s_{\alpha_2+\alpha_4+\alpha_5+\alpha_6+\alpha_7}
       s_{\alpha^{d6}_{max}}(s_3 s_4 s_5 s_6). \\
  \end{split}
\end{equation*}
All simple roots of $E_7$, except $\alpha_1$ and $\alpha_6$,
are orthogonal to $\alpha^{d6}_{max}$. Since
\begin{equation*}
  \begin{split}
     & s_{\alpha^{d6}_{max}}s_6 = s_6 s_{\alpha^{d6}_{max} - \alpha_6}, \\
     & s_6 (s_{\alpha^{d6}_{max} - \alpha_6})s_6 = s_{\alpha^{d6}_{max}},
  \end{split}   
\end{equation*}
we get,
\begin{equation*}
  \begin{split}
      w^{(1)}_1 & w^{(2)}_1 w^{(3)}_1 = \\
         & s_{\alpha_7}s_{\alpha_4+\alpha_5} s_{\alpha_2+\alpha_4+\alpha_5+\alpha_6+\alpha_7}
          (s_3 s_4 s_5) s_{\alpha^{d6}_{max}} s_6 =  \\
         & s_{\alpha_7}s_{\alpha_4+\alpha_5} s_{\alpha_2+\alpha_4+\alpha_5+\alpha_6+\alpha_7}
          (s_3 s_4 s_5 s_6) s_{\alpha^{d6}_{max}-\alpha_6}  =  \\
         & s_{\alpha_7}s_{\alpha_4+\alpha_5}(s_3 s_4 s_5)
               s_{\alpha_2+\alpha_3+2\alpha_4+2\alpha_5+\alpha_6+\alpha_7}
          s_6 s_{\alpha^{d6}_{max}-\alpha_6}  =  \\
         & s_{\alpha_7}s_{\alpha_4+\alpha_5}(s_3 s_4 s_5)
               (s_{\alpha^{d6}_{max}-\alpha_6}(s_6) s_{\alpha^{d6}_{max}-\alpha_6})  =  \\
         & s_{\alpha_7}s_{\alpha_4+\alpha_5}(s_3 s_4 s_5)s_{\alpha^{d6}_{max}}.
  \end{split}
\end{equation*}
By Lemma \ref{lem_permute_An_1}, we get
\begin{equation}
  \label{eq_roots_w1_1}
  w^{(1)}_1 w^{(2)}_1 w^{(3)}_1 = s_{\alpha_7}s_{\alpha_4} s_{\alpha_3 + \alpha_4+\alpha_5}s_{\alpha^{d6}_{max}}
\end{equation}
It is easy to see the following decomposition:
\begin{equation}
  \label{eq_roots_w1_2}
  \begin{split}
  w^{(4)}_1 = & s_2 s_4 s_5 s_3 s_4 s_2 = (s_2s_4s_5 s_4s_2)(s_2s_4s_3 s_4s_2) = \\
    & s_{\alpha_2 + \alpha_4 + \alpha_5} s_{\alpha_2 + \alpha_4 + \alpha_3}.
   \end{split}
\end{equation}
Since $\alpha^{d6}_{max} \perp \alpha_2, \alpha_3, \alpha_4, \alpha_5$,  the reflection
$s_{\alpha^{d6}_{max}}$ commute with $w^{(4)}_1$.
Putting together \eqref{eq_roots_w1_1} and \eqref{eq_roots_w1_2}, we get
\begin{equation*}
   w_1 =  w^{(1)}_1 w^{(2)}_1 w^{(3)}_1w^{(4)}_1 =
       s_{\alpha_7}s_{\alpha_4} s_{\alpha_3 + \alpha_4+\alpha_5}
       s_{\alpha_2 + \alpha_4 + \alpha_5}
       s_{\alpha_2 + \alpha_4 + \alpha_3}s_{\alpha^{d6}_{max}}.
\end{equation*}
This decomposition corresponds to the roots
\eqref{eq_6_roots_of_E7}. Their mutual orthogonality is easily verified.
~\\

 (ii) We consequentially apply Proposition \ref{prop_fact_1} to the element $w_2$:
\begin{equation}
 \label{eq_alpha_max_E7}
 \footnotesize
 \begin{split}
   & w_2 = \\
   & s_1 s_3 s_4 s_5 s_2 s_4 s_6 (s_5 s_7 s_6 s_3 (s_4 s_5 s_2 (s_4 s_3 s_1 s_3 s_4) s_2 s_5 s_4)
        s_3 s_6 s_7 s_5)s_6 s_4 s_2 s_5 s_4 s_3 s_1 = \\
   & s_1 s_3 s_4 s_5 s_2 s_4 s_6 (s_5 s_7 s_6 s_3 (s_4 s_5 s_2 (s_{\alpha_1 + \alpha_3 + \alpha_4})
        s_2 s_5 s_4) s_3 s_6 s_7 s_5) s_6 s_4 s_2 s_5 s_4 s_3 s_1  = \\
   & s_1 s_3 s_4 s_5 s_2 s_4 s_6 (s_5 s_7 s_6 s_3
          (s_{\alpha_1 + \alpha_2 + \alpha_3 + 2\alpha_4 + \alpha_5}) s_3 s_6 s_7 s_5)
             s_6 s_4 s_2 s_5  s_4 s_3 s_1  = \\
   & s_1 s_3 s_4 s_5 s_2 s_4 s_6 (s_5 s_7
          (s_{\alpha_1 + \alpha_2 + 2\alpha_3 + 2\alpha_4 + \alpha_5 +\alpha_6}) s_7 s_5)
             s_6 s_4 s_2 s_5  s_4 s_3 s_1  = \\
   & s_1 s_3 s_4 (s_5 s_2 s_4 s_6
          (s_{\alpha_1 + \alpha_2 + 2\alpha_3 + 2\alpha_4 + 2\alpha_5 +\alpha_6+\alpha_7})
            s_6 s_4 s_2 s_5 ) s_4 s_3 s_1  = \\
   & s_1 s_3 s_4 (s_5 s_2
          (s_{\alpha_1 + \alpha_2 + 2\alpha_3 + 3\alpha_4 + 2\alpha_5 +2\alpha_6+\alpha_7})
            s_2 s_5 ) s_4 s_3 s_1  = \\
   & s_1 s_3 s_4
          (s_{\alpha_1 + 2\alpha_2 + 2\alpha_3 + 3\alpha_4 + 3\alpha_5 +2\alpha_6+\alpha_7})
            s_4 s_3 s_1  = \\
   & s_1  (s_{\alpha_1 + 2\alpha_2 + 3\alpha_3 + 4\alpha_4 + 3\alpha_5 +2\alpha_6+\alpha_7})
            s_1  = \\
   & s_{2\alpha_1 + 2\alpha_2 + 3\alpha_3 + 4\alpha_4 + 3\alpha_5 +2\alpha_6+\alpha_7}.
    %% s_{\alpha_{max}}.
 \end{split}
\end{equation}
In \eqref{eq_alpha_max_E7}, the last reflection
corresponds to the highest root $\alpha_{max}$ of $E_7$, see Table \ref{tab_highest_roots_EFG}.
\qed
~\\

\begin{proposition}
  \label{prop_factoriz_E7_case_0}
  {\rm{(i)}} The element $w_0 = w_1 w_2$ is the product of 7 reflections associated
with the following mutually orthogonal positive roots:
\begin{equation}
  \label{eq_7_roots_of_E7}
     \alpha_4, \; \alpha_7, \; \alpha^{d6}_{max}, \; \alpha_{max}, \;
      \alpha_3 + \alpha_4 + \alpha_5, \; \alpha_2 + \alpha_4 + \alpha_5,
         \; \alpha_2 + \alpha_4 + \alpha_3.
\end{equation}

{\rm{(ii)}} The element $w_0$ is the longest element in $W(E_7)$.
\end{proposition}

\PerfProof
(i) The highest root $\alpha_{max}$ is orthogonal to all simple roots except $\alpha_1$,
i.e.,  $\alpha_{max}$ is orthogonal to all other roots of \eqref{eq_7_roots_of_E7}.
Taking Proposition \ref{prop_6_orthog}(i) into account, we obtain that all roots in
 \eqref{eq_7_roots_of_E7} are mutually orthogonal.
~\\

(ii) We put
\begin{equation*}
  \alpha_{345} := \alpha_3 + \alpha_4 + \alpha_5,\; \alpha_{245} := \alpha_2 + \alpha_4 + \alpha_5, \;
  \alpha_{234} := \alpha_2 + \alpha_3 + \alpha_4.
\end{equation*}
Consider the determinant constructed by roots of \eqref{eq_7_roots_of_E7},
see (\ref{eq_determ_E7}a). We transform the rows of the determinant  as follows:
\begin{equation*}
 \begin{array}{llllll}
   & \alpha_{345} & \rightarrow & \alpha_{345} - \alpha_{245}+ \alpha_{234} -\alpha_4 & = &(0,0,2,0,0,0,0,0), \\
   & \alpha_{245} & \rightarrow & \alpha_{245} - \alpha_{234}+ \alpha_{345} -\alpha_4 & = & (0,0,0,0,2,0,0,0), \\
   & \alpha_{234} & \rightarrow & \alpha_{245} - \alpha_{345}+ \alpha_{234} -\alpha_4 & = & (0,2,0,0,0,0,0,0), \\
   & \alpha^{d6}_{max} & \rightarrow &
            \alpha^{d6}_{max} - \alpha_{245} - \alpha_{345} - \alpha_7 & = & (0,0,0,0,0,0,2,0), \\
   & \alpha_{max} & \rightarrow & \alpha_{max} -  \alpha^{d6}_{max} - \alpha_{345} - \alpha_{234} & = &
               (2,0,0,0,0,0,0,0).
 \end{array}
\end{equation*}
The obtained determinant (\ref{eq_determ_E7}b) is non-degenerate.
\begin{equation}
  \label{eq_determ_E7}
 \footnotesize
 \begin{split}
  &\left | \left |
  \begin{array}{lllllll}
     0 & 0 & 1 & 1 & 1 & 0 & 0 \\
     0 & 1 & 0 & 1 & 1 & 0 & 0 \\
     0 & 1 & 1 & 1 & 0 & 0 & 0 \\
     0 & 0 & 0 & 1 & 0 & 0 & 0 \\
     0 & 1 & 1 & 2 & 2 & 2 & 1 \\
     0 & 0 & 0 & 0 & 0 & 0 & 1 \\
     2 & 2 & 3 & 4 & 3 & 2 & 1 \\
  \end{array}
  \right | \right |
  \begin{array}{l}
     \alpha_{345} \\
     \alpha_{245} \\
     \alpha_{234} \\
     \alpha_4 \\
     \alpha^{d6}_{max} \\
     \alpha_7 \\
     \alpha_{max} \\
  \end{array}
     \Longrightarrow
  \left | \left |
  \begin{array}{llllllll}
     0 & 0 & 2 & 0 & 0 & 0 & 0 \\
     0 & 0 & 0 & 0 & 2 & 0 & 0 \\
     0 & 2 & 0 & 0 & 0 & 0 & 0 \\
     0 & 0 & 0 & 1 & 0 & 0 & 0 \\
     0 & 0 & 0 & 0 & 0 & 2 & 0 \\
     0 & 0 & 0 & 0 & 0 & 0 & 1 \\
     2 & 0 & 0 & 0 & 0 & 0 & 0 \\
  \end{array}
  \right | \right | \\
  & \qquad \qquad \quad \quad (a)
    \qquad \qquad \qquad \qquad \qquad \qquad \qquad \qquad (b)
  \end{split}
\end{equation}
~\\
Thus, the roots \eqref{eq_7_roots_of_E7} are linearly independent
and constitutes the basis of the whole space.
Since the routs \eqref{eq_7_roots_of_E7} are mutually orthogonal,
$w_0$ acts as \;$-1$ on each of these roots. Hence $w_0$ acts as $-1$
on the whole space, i.e., $w_0$ acts in the same way
as the longest element. Thus $w_0$ is the longest in $W(E_7)$. \qed
~\\

The roots \eqref{eq_7_roots_of_E7} are not satisfy to the max-orthogonality conditions, 
see $\S$\ref{sec_uniqueness}, since
$\{\alpha_3 + \alpha_4 + \alpha_5, \, \alpha_2 + \alpha_4 + \alpha_5, \, \alpha_2 + \alpha_4 + \alpha_3\}$
are not linearly ordered. The following proposition provides another decomposition of $w_0$,
which solves this issue.

\begin{proposition}[decomposition in $W(E_7)$]
  \label{prop_factoriz_E7}
 The element $w_0$ can be decomposed as the product of the following $7$
reflections  associated with the \underline{mutually orthogonal} roots:
\begin{equation}
  \label{eq_7_roots_of_E7_2}
    \alpha_2, \; \alpha_3, \; \alpha_5 \; \alpha_7, \;
              \alpha^{d4}_{max}, \; \alpha^{d6}_{max}, \; \alpha_{max}.
\end{equation}
The highest roots in \eqref{eq_7_roots_of_E7_2} are ordered as follows:
\begin{equation*}
     \alpha^{d4}_{max} < \alpha^{d6}_{max} < \alpha_{max}.  \\
\end{equation*}
\end{proposition}

\PerfProof
  We apply conjugation $s_{\alpha_3 + \alpha_4}$  to
  the reflections associated with roots \eqref{eq_7_roots_of_E7}.
  First, this conjugation does not change $s_7$, $s_{\alpha^{d6}_{max}}$ and
  $s_{\alpha_{max}}$ since
  $\{\alpha_7, \alpha^{d6}_{max}, \alpha_{max}\} \perp \{\alpha_3, \alpha_4\}$.
  For the remaining roots, we have as follows:
\begin{equation*}
  \begin{split}
  & s_{\alpha_3 + \alpha_4}s_{\alpha_4}s_{\alpha_3 + \alpha_4} = s_{\alpha_3}, \\
  & s_{\alpha_3 + \alpha_4}s_{\alpha_3 + \alpha_4 + \alpha_5}s_{\alpha_3 +
     \alpha_4} = s_{\alpha_5}, \\
  & s_{\alpha_3 + \alpha_4}s_{\alpha_2 + \alpha_4 + \alpha_5}s_{\alpha_3 + \alpha_4}
   = s_{\alpha_2 + \alpha_3 +2\alpha_4 + \alpha_5} = s_{\alpha^{d4}_{max}}, \\
  & s_{\alpha_3 + \alpha_4}s_{\alpha_2 + \alpha_4 + \alpha_3}s_{\alpha_3 +
     \alpha_4} = s_{\alpha_2}. \\
  \end{split}
\end{equation*}
\qed
~\\

\section{\bf Decomposition in $W(E_8)$}

\begin{figure}[h]
\centering
   \includegraphics[scale=0.4]{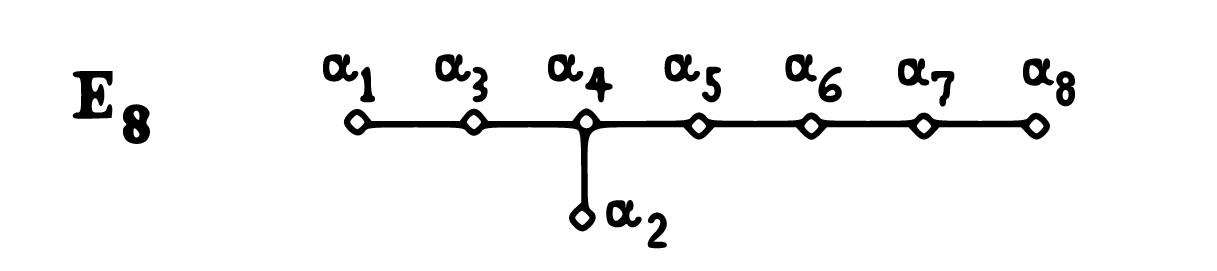}
\caption{Numbering of simple roots in $E_8$}
%%%%%% The label must come after caption
\label{fig_E8}
\end{figure}

We denote the \underline{highest roots}
in the root subsystems $D_4 \subset D_6 \subset E_7 \subset E_8$  as follows:
\begin{equation}
  \label{eq_not_E8_1}
  \begin{split}
     & \alpha^{d4}_{max} := \alpha_2 + \alpha_3 + 2\alpha_4 + \alpha_5, \\
     & \alpha^{d6}_{max} := \alpha_2 + \alpha_3 + 2\alpha_4 + 2\alpha_5 + 2\alpha_6 + \alpha_7,\\
     & \alpha^{e7}_{max} :=  2\alpha_1 + 2\alpha_2 + 3\alpha_3 + 4\alpha_4 +
                 3\alpha_5 + 2\alpha_6 + \alpha_7,
   \end{split}
\end{equation}
see Table \ref{tab_factoriz_all} and Remark \ref{rem_notations}.
In \eqref{eq_not_E8_1}, the indices of simple roots differ from the indices
of simple roots that appears as summands in $\alpha^{d4}_{max}$ and $\alpha^{d6}_{max}$
defined in \eqref{eq_not_Dn_1}, see Remark \ref{rem_differ_indices} and note after \eqref{eq_not_E7_1}. 
The numbering of vertices of the Dynkin diagram $E_8$ is shown
in Fig. \ref{fig_E8}.

\subsection{Reducing $w_0$ from $120$ factors to $14$}

Let us consider the following $4$ elements in the Weyl group $W(E_8)$:
\begin{equation}
 \label{eq_elems_1234}
  \begin{split}
    w_1 =
 & s_8 (s_7 s_8) (s_6 s_7 s_8) (s_5 s_6 s_7 s_8) (s_4 s_5 s_6 s_7 s_8) (s_3 s_4 s_5 s_6 s_7 s_8), \\
    w_2 =
 & s_2 s_4 s_5 s_6 (s_3 s_4 s_5 (s_7 s_6 s_2 s_4 (s_3 s_5 s_4 s_2 (s_8 s_7 s_6 s_5 (s_4 s_3 s_1 s_3 s_4 ) \times \\
 &   s_5 s_6 s_7 s_8)s_2 s_4 s_5 s_3) s_4 s_2 s_6 s_7) s_5 s_4 s_3) s_6 s_5 s_4 s_2, \\
    w_3 =
 & s_1 s_3 (s_4 s_2 s_5 s_4)( s_3 s_6 s_5 s_4 )(s_7 s_6 s_8 s_7)
         (s_5 s_6) (s_2 s_4 s_5) (s_3 s_4) s_2, \\
    w_4 =
 & s_1 (s_3 s_4 s_5 s_6 (s_2 s_4 s_5 s_3 (s_4 s_7 s_6 s_5
 (s_2 s_4 s_3 s_8 s_1 s_3 s_4 s_2 )s_5 s_6 s_7 s_4 ) \times \\
 &   s_3 s_5 s_4 s_2) s_6 s_5 s_4 s_3) s_1.
  \end{split}
\end{equation}
Denote by $w_0$ the product of elements \eqref{eq_elems_1234}:
\begin{equation}
  \label{eq_w0_def}
  w_0 = w_1 w_2 w_3 w_4
\end{equation}
We will prove that $w_0$ is the longest element in $W(E_8)$.
But first we have to transform them so that they have as few factors as possible.
\begin{proposition}
  Elements $w_1$, $w_2$, $w_3$, $w_4$ are transformed as follows:
\begin{equation}
 \label{eq_transf_14_factors}
  \begin{split}
    w_1 = &s_{\alpha_5 + \alpha_6}
      s_{\alpha_4 + \alpha_5 + \alpha_6 + \alpha_7}
      s_{\alpha_3 + \alpha_4 + \alpha_5 + \alpha_6 + \alpha_7 + \alpha_8}, \\
    w_2 = &s_{\alpha_1 + 3\alpha_2 + 3\alpha_3 + 5\alpha_4 +
                   4\alpha_5 +3\alpha_6+ 2\alpha_7 + \alpha_8}, \\
    w_3 = &s_1 s_{\alpha_2 + \alpha_4}
            s_{\alpha_3 + \alpha_4 + \alpha_5} s_{\alpha_2 + \alpha_3 + 2\alpha_4 +
               \alpha_5 + \alpha_6} \\
          & s_{\alpha_2 + \alpha_3 + \alpha_4 + \alpha_5 +\alpha_6 + \alpha_7}
             s_{\alpha_7 + \alpha_8} s_6 s_{\alpha_2 + \alpha_4 + \alpha_5}, \\
                %% & \footnotesize \text{(8 factors)} \\
    w_4= &s_{\alpha_1 + \alpha_3 + \alpha_4 + \alpha_5 + \alpha_6 + \alpha_7+\alpha_8} s_{2\alpha_1 + 2\alpha_2 + 3\alpha_3 + 4\alpha_4 +  3\alpha_5 + 2\alpha_6 + \alpha_7}.
  \end{split}
\end{equation}
 \end{proposition}

\underline{Element $w_1$}.  We use Propositions \ref{prop_factoriz_An} and
\ref{prop_factoriz_An_2} with one small change: we take $w_0$ for $A_6$ and
increment all indices by $2$.  Let us denote this operation by
$(\dots)_{i \rightarrow i+2}$.
\begin{equation*}
\begin{split}
  w_1 = & s_8 (s_7 s_8) (s_6 s_7 s_8) (s_5 s_6 s_7 s_8) (s_4 s_5 s_6 s_7 s_8) (s_3 s_4 s_5 s_6 s_7 s_8) = \\
        & (s_6 (s_5 s_6) (s_4 s_5 s_6) (s_3 s_4 s_5 s_6) (s_2 s_3 s_4 s_5 s_6) (s_1 s_2 s_3 s_4 s_5 s_6))_{i \rightarrow i+2} = \\
    & (s_{\alpha_3 + \alpha_4}
    s_{\alpha_2 + \alpha_3 + \alpha_4 + \alpha_5}
    s_{\alpha_1 + \alpha_2 + \alpha_3 + \alpha_4 + \alpha_5 + \alpha_6})_{i \rightarrow i+2} = \\
   & s_{\alpha_5 + \alpha_6}
    s_{\alpha_4 + \alpha_5 + \alpha_6 + \alpha_7}
    s_{\alpha_3 + \alpha_4 + \alpha_5 + \alpha_6 + \alpha_7 + \alpha_8}.
\end{split}
\end{equation*}
The element $w_1$ is decomposed into $3$ factors.
~\\

\underline{Element $w_2$}.
As in Proposition \ref{prop_factoriz_E6}, we repeatedly apply Proposition \ref{prop_fact_1}.
For the beginning, let us transform $5$ factors in parentheses lying in the center of the element $w_2$:
\begin{equation}
 \label{eq_w2_1_E8}
 \begin{split}
  & s_3 s_5 s_4 s_2 (s_8 s_7 s_6 s_5 (s_4 s_3 s_1 s_3 s_4)
    s_5 s_6 s_7 s_8)s_2 s_4 s_5 s_3 = \\
  & s_3 s_5 s_4 s_2(s_8 s_7 s_6 s_5 (s_{\alpha_1 + \alpha_3 + \alpha_4})s_5 s_6 s_7 s_8)
     s_2 s_4 s_5 s_3 = \\
  & s_3 s_5 s_4 s_2 (s_{\alpha_1 + \alpha_3 + \alpha_4 + \alpha_5 +\alpha_6+ \alpha_7 +
     \alpha_8}) s_2 s_4 s_5 s_3 = \\
  & s_{\alpha_1 + \alpha_2 + 2\alpha_3 + 2\alpha_4 + 2\alpha_5 +\alpha_6+ \alpha_7 +
     \alpha_8}.
  \end{split}
\end{equation}
By \eqref{eq_w2_1_E8} we get
\begin{equation*}
 \begin{split}
   w_2 & =
   s_2 s_4 s_5 s_6 (s_3 s_4 s_5 (s_7 s_6 s_2 s_4 (s_{\alpha_1 + \alpha_2 + 2\alpha_3 + 2\alpha_4 + 2\alpha_5 +\alpha_6+ \alpha_7 +
     \alpha_8}) \times \\
  & s_4 s_2 s_6 s_7)s_5 s_4 s_3) s_6 s_5 s_4 s_2 = \\
  & s_2 s_4 s_5 s_6 ( s_3 s_4 s_5 (s_{\alpha_1 + 2\alpha_2 + 2\alpha_3 + 3\alpha_4 + 2\alpha_5 +2\alpha_6+ 2\alpha_7 + \alpha_8}) s_5 s_4 s_3)
   s_6 s_5 s_4 s_2 = \\
  &  s_2 s_4 s_5 s_6 (s_{\alpha_1 + 2\alpha_2 + 3\alpha_3 + 4\alpha_4 + 3\alpha_5 +2\alpha_6+ 2\alpha_7 + \alpha_8}) s_6 s_5 s_4 s_2 = \\
  &  s_2 (s_{\alpha_1 + 2\alpha_2 + 3\alpha_3 + 5\alpha_4 + 4\alpha_5 +3\alpha_6+ 2\alpha_7 + \alpha_8})s_2 = \\
  & s_{\alpha_1 + 3\alpha_2 + 3\alpha_3 + 5\alpha_4 + 4\alpha_5 +3\alpha_6+ 2\alpha_7 + \alpha_8}. \\
  \end{split}
\end{equation*}
The element $w_2$ has only $1$ factor.
~\\

\underline{Element $w_3$}.
The following identities hold:
\begin{equation}
  \label{eq_ident_1_E8}
 \begin{split}
   & s_4 s_2 s_5 s_4 = (s_4 s_2 s_4)(s_4 s_5 s_4) =
                s_{\alpha_2 + \alpha_4}s_{\alpha_5 + \alpha_4},  \\
   & s_7 s_6 s_8 s_7 = (s_7 s_6 s_7) (s_7 s_8 s_7) =
                s_{\alpha_6 + \alpha_7}s_{\alpha_7 + \alpha_8},  \\
   & s_2 s_4 s_5 s_3 s_4 s_2 = (s_2 s_4 s_5 s_4 s_2) (s_2 s_4 s_3 s_4 s_2) =
       s_{\alpha_2 + \alpha_4 + \alpha_5} s_{\alpha_2 + \alpha_4 + \alpha_3}
 \end{split}
\end{equation}
By \eqref{eq_ident_1_E8} we have
\begin{equation*}
 \begin{split}
   &(s_4 s_2 s_5 s_4)(s_3 s_6 s_5 s_4) =  s_{\alpha_2 + \alpha_4}s_{\alpha_5 + \alpha_4}
            (s_3 s_6 s_5 s_4) = \\
   & s_{\alpha_2 + \alpha_4}s_5s_4 s_5 s_3 s_6 s_5 s_4 =
     s_{\alpha_2 + \alpha_4}s_5s_4 s_3 s_5 s_6 s_5 s_4 = \\
   &  s_{\alpha_2 + \alpha_4}s_5 s_4 s_3 s_4 (s_4 s_5 s_6 s_5 s_4) = \\
   & s_{\alpha_2 + \alpha_4}s_5 s_{\alpha_3 + \alpha_4} (s_{\alpha_4 + \alpha_5 + \alpha_6}).
 \end{split}
\end{equation*}
Then,
\begin{equation*}
 \begin{split}
  w_3 &=  s_1 s_3 (s_4 s_2 s_5 s_4)( s_3 s_6 s_5 s_4)(s_7 s_6 s_8 s_7)
         (s_5 s_6) (s_2 s_4 s_5) (s_3 s_4) s_2 = \\
        & s_1 s_3 s_{\alpha_2 + \alpha_4}s_5 s_{\alpha_3 + \alpha_4}
               s_{\alpha_4 + \alpha_5 + \alpha_6}
               s_{\alpha_6 + \alpha_7}s_{\alpha_7 + \alpha_8}
         s_5 s_6 s_{\alpha_2 + \alpha_4 + \alpha_5} s_{\alpha_2 + \alpha_4 + \alpha_3} = \\
       & s_1 s_3 s_{\alpha_2 + \alpha_4}
       \big (s_5 s_{\alpha_3 + \alpha_4} s_{\alpha_4 + \alpha_5 + \alpha_6}
           s_{\alpha_6 + \alpha_7} s_5 \big )
     s_{\alpha_7 + \alpha_8} s_6 s_{\alpha_2 + \alpha_4 + \alpha_5}
     s_{\alpha_2 + \alpha_4 + \alpha_3}.
 \end{split}
\end{equation*}
Applying conjugation with $s_5$ to the parenthesized expression in \eqref{eq_ident_1_E8}
we reduce the number of factors by $2$:
\begin{equation*}
 \begin{split}
  w_3 & =  s_1 s_3 s_{\alpha_2 + \alpha_4}
           s_{\alpha_3 + \alpha_4 + \alpha_5} s_{\alpha_4 + \alpha_5 + \alpha_6}
           s_{\alpha_5 +\alpha_6 + \alpha_7}
           s_{\alpha_7 + \alpha_8} \times \\           
    &  s_6s_{\alpha_2 + \alpha_4 + \alpha_5} s_{\alpha_2 + \alpha_3 + \alpha_4}. \\
 \end{split}
\end{equation*}
 Since $s_3 s_{\alpha_2 + \alpha_4} = s_{\alpha_2 + \alpha_4} s_{\alpha_2 + \alpha_3 + \alpha_4}$
 and $s_6$ commute with $s_{\alpha_2 + \alpha_3 + \alpha_4}$ we get the following:
\begin{equation*}
 \begin{split}        
   w_3 & =   s_1 s_{\alpha_2 + \alpha_4} s_{\alpha_2 + \alpha_3 + \alpha_4}
           s_{\alpha_3 + \alpha_4 + \alpha_5} s_{\alpha_4 + \alpha_5 + \alpha_6}
           s_{\alpha_5 +\alpha_6 + \alpha_7}
           s_{\alpha_7 + \alpha_8}  \times \\
    &  s_{\alpha_2 + \alpha_3 + \alpha_4}s_6s_{\alpha_2 + \alpha_4 + \alpha_5}.
 \end{split}
\end{equation*}
Further, we can swap $s_{\alpha_2 + \alpha_3 + \alpha_4}$ and $s_{\alpha_3 + \alpha_4 + \alpha_5}$ 
since the roots $\alpha_2 + \alpha_3 + \alpha_4$ and $\alpha_3 + \alpha_4 + \alpha_5$ are orthogonal.
For reflections $s_{\alpha_2 + \alpha_3 + \alpha_4}$  and
 $s_{\alpha_4 + \alpha_5 + \alpha_6}$, we have
\begin{equation*}
    s_{\alpha_2 + \alpha_3 + \alpha_4}s_{\alpha_4 + \alpha_5 + \alpha_6} =
    s_{\alpha_2 + \alpha_3 + 2\alpha_4 +\alpha_5 +\alpha_6}s_{\alpha_2 + \alpha_3 + \alpha_4}.
\end{equation*}
Then,
\begin{equation}
 \label{eq_ident_3_E8}
 %%\footnotesize
 \begin{split}
  w_3 =  & s_1 s_{\alpha_2 + \alpha_4}
         s_{\alpha_3 + \alpha_4 + \alpha_5} s_{\alpha_2 + \alpha_3 + 2\alpha_4 +
           \alpha_5 + \alpha_6} \\
        & \big ( s_{\alpha_2 + \alpha_3 + \alpha_4}s_{\alpha_5 +\alpha_6 + \alpha_7}
           s_{\alpha_2 + \alpha_3 + \alpha_4} \big ) s_{\alpha_7 + \alpha_8} s_6
           s_{\alpha_2 + \alpha_4 + \alpha_5}.
   \end{split}
\end{equation}
  Applying conjugation with $s_{\alpha_2 + \alpha_3 + \alpha_4}$
to the expression in brackets in \eqref{eq_ident_3_E8}
we again reduce the number of factors by $2$:
\begin{equation*}
 \begin{split}
  w_3 =  & s_1 s_{\alpha_2 + \alpha_4}
         s_{\alpha_3 + \alpha_4 + \alpha_5} s_{\alpha_2 + \alpha_3 + 2\alpha_4 +
           \alpha_5 + \alpha_6} \\
        & s_{\alpha_2 + \alpha_3 + \alpha_4 + \alpha_5 +\alpha_6 + \alpha_7}
           s_{\alpha_7 + \alpha_8} s_6 s_{\alpha_2 + \alpha_4 + \alpha_5}.
   \end{split}
\end{equation*}
The element $w_3$ has $8$ factors.
~\\

\underline{Element $w_4$}.
Let us transform $3$ factors in parentheses lying in the middle of the element $w_4$:
\begin{equation*}
 \begin{split}
    & s_4 s_7 s_6 s_5 (s_2 s_4 s_3 s_8 s_1 s_3 s_4 s_2 )s_5 s_6 s_7 s_4 = \\
    & s_4 s_7 s_6 s_5 (s_8s_{\alpha_1 + \alpha_3 + \alpha_4 + \alpha_2})s_5 s_6 s_7 s_4 = \\
    & s_{\alpha_7+\alpha_8} s_{\alpha_1 + \alpha_3 + 2\alpha_4 + \alpha_2 + \alpha_5 + \alpha_6 + \alpha_7}. \\
  \end{split}
\end{equation*}
Let us add one more factor on each side, see $w_4$ in \eqref{eq_elems_1234}:
\begin{equation*}
 \begin{split}
   & s_2 s_4 s_5 s_3 (s_{\alpha_7+\alpha_8} s_{\alpha_1 +  \alpha_2 + \alpha_3 + 2\alpha_4 + \alpha_5 + \alpha_6 + \alpha_7}) s_3 s_5 s_4 s_2 = \\
   & s_{\alpha_7+\alpha_8} s_{\alpha_1 + 2\alpha_2 + 2\alpha_3 + 3\alpha_4 +  2\alpha_5 + \alpha_6 + \alpha_7}.
 \end{split}
\end{equation*}
And add one more factor:
\begin{equation*}
 \begin{split}
   & s_3 s_4 s_5 s_6(s_{\alpha_7+\alpha_8} s_{\alpha_1 + 2\alpha_2 + 2\alpha_3 + 3\alpha_4 +  2\alpha_5 + \alpha_6 + \alpha_7})s_6 s_5 s_4 s_3 = \\
   & s_{\alpha_3 + \alpha_4 + \alpha_5 + \alpha_6 + \alpha_7+\alpha_8} s_{\alpha_1 + 2\alpha_2 + 3\alpha_3 + 4\alpha_4 +  3\alpha_5 + 2\alpha_6 + \alpha_7}.\\
 \end{split}
\end{equation*}
At last,
\begin{equation*}
 \begin{split}
   w_4 = &
     s_1s_{\alpha_3 + \alpha_4 + \alpha_5 + \alpha_6 + \alpha_7+\alpha_8} s_{\alpha_1 + 2\alpha_2 + 3\alpha_3 + 4\alpha_4 +  3\alpha_5 + 2\alpha_6 + \alpha_7}s_1 = \\
     & s_{\alpha_1 + \alpha_3 + \alpha_4 + \alpha_5 + \alpha_6 + \alpha_7+\alpha_8} s_{2\alpha_1 + 2\alpha_2 + 3\alpha_3 + 4\alpha_4 +  3\alpha_5 + 2\alpha_6 + \alpha_7}. \\
 \end{split}
\end{equation*}
The element $w_4$ is decomposed into $2$ factors. \qed
~\\

Considering all the factors together, the element $w_0$ has $14$ factors.

\subsection{Reducing $w_0$ from $14$ factors to $12$}

 Consider the product of $14$ factors from \eqref{eq_transf_14_factors}:
\begin{equation}
 \label{eq_reduce_1_E8}
 \begin{split}
  w_0 =
      & s_{\alpha_5 + \alpha_6}s_{\alpha_4 + \alpha_5 + \alpha_6 + \alpha_7}
        s_{\alpha_3 + \alpha_4 + \alpha_5 + \alpha_6 + \alpha_7 + \alpha_8} \times \\
      &  s_{\alpha_1 + 3\alpha_2 + 3\alpha_3 + 5\alpha_4 +
                   4\alpha_5 +3\alpha_6+ 2\alpha_7 + \alpha_8}
        s_1 s_{\alpha_2 + \alpha_4} s_{\alpha_3 + \alpha_4 + \alpha_5} \times \\
      &   s_{\alpha_2 + \alpha_3 + 2\alpha_4 +
               \alpha_5 + \alpha_6}
        s_{\alpha_2 + \alpha_3 + \alpha_4 + \alpha_5 +\alpha_6 + \alpha_7}
             s_{\alpha_7 + \alpha_8} s_6 s_{\alpha_2 + \alpha_4 + \alpha_5} \times \\
      & s_{\alpha_1 + \alpha_3 + \alpha_4 + \alpha_5 + \alpha_6 + \alpha_7+\alpha_8}
         s_{2\alpha_1 + 2\alpha_2 + 3\alpha_3 + 4\alpha_4 +  3\alpha_5 + 2\alpha_6 + \alpha_7}.
 \end{split}
\end{equation}
~\\
First, it can be checked that the following roots are orthogonal
\begin{equation*}
 \label{eq_step_1_E8}
 \footnotesize
 \begin{split}
   \alpha_3 + \alpha_4 + \alpha_5 + \alpha_6 + \alpha_7 + \alpha_8 \; \perp \;
      \alpha_1 + 3\alpha_2 + 3\alpha_3 + 5\alpha_4 +
                   4\alpha_5 +3\alpha_6+ 2\alpha_7 + \alpha_8,
  \end{split}
\end{equation*}
i.e., the corresponding reflections in \eqref{eq_reduce_1_E8} can be swapped:
\begin{equation}
 \label{eq_reduce_1_E8_2}
 \begin{split}
   w_0 =
      & s_{\alpha_5 + \alpha_6}s_{\alpha_4 + \alpha_5 + \alpha_6 + \alpha_7}
        s_{\alpha_1 + 3\alpha_2 + 3\alpha_3 + 5\alpha_4 +
                   4\alpha_5 +3\alpha_6+ 2\alpha_7 + \alpha_8} \times \\
      & s_{\alpha_3 + \alpha_4 + \alpha_5 + \alpha_6 + \alpha_7 + \alpha_8}
        s_1 s_{\alpha_2 + \alpha_4} s_{\alpha_3 + \alpha_4 + \alpha_5} \times \\
      &   s_{\alpha_2 + \alpha_3 + 2\alpha_4 +
               \alpha_5 + \alpha_6}
        s_{\alpha_2 + \alpha_3 + \alpha_4 + \alpha_5 +\alpha_6 + \alpha_7}
             s_{\alpha_7 + \alpha_8} s_6 s_{\alpha_2 + \alpha_4 + \alpha_5} \times \\
      & s_{\alpha_1 + \alpha_3 + \alpha_4 + \alpha_5 + \alpha_6 + \alpha_7+\alpha_8}
         s_{2\alpha_1 + 2\alpha_2 + 3\alpha_3 + 4\alpha_4 +  3\alpha_5 + 2\alpha_6 + \alpha_7}.
 \end{split}
\end{equation}
Using conjugation with $s_1$, by Proposition \ref{prop_fact_1} we have:
\begin{equation*}
   s_{\alpha_3 + \alpha_4 + \alpha_5 + \alpha_6 + \alpha_7 + \alpha_8}s_1 =
     s_1 s_{\alpha_1 + \alpha_3 + \alpha_4 + \alpha_5 + \alpha_6 + \alpha_7 + \alpha_8}
\end{equation*}
Then, from \eqref{eq_reduce_1_E8_2} we have
\begin{equation*}
 \begin{split}
  w_0 =
      & s_{\alpha_5 + \alpha_6}s_{\alpha_4 + \alpha_5 + \alpha_6 + \alpha_7}
        s_{\alpha_1 + 3\alpha_2 + 3\alpha_3 + 5\alpha_4 +
                   4\alpha_5 +3\alpha_6+ 2\alpha_7 + \alpha_8} \times \\
      & s_1 s_{\alpha_1 + \alpha_3 + \alpha_4 + \alpha_5 + \alpha_6 + \alpha_7 + \alpha_8}
        s_{\alpha_2 + \alpha_4} s_{\alpha_3 + \alpha_4 + \alpha_5} \times \\
      &   s_{\alpha_2 + \alpha_3 + 2\alpha_4 +
               \alpha_5 + \alpha_6}
        s_{\alpha_2 + \alpha_3 + \alpha_4 + \alpha_5 +\alpha_6 + \alpha_7}
             s_{\alpha_7 + \alpha_8} s_6 s_{\alpha_2 + \alpha_4 + \alpha_5} \times \\
      & s_{\alpha_1 + \alpha_3 + \alpha_4 + \alpha_5 + \alpha_6 + \alpha_7+\alpha_8}
         s_{2\alpha_1 + 2\alpha_2 + 3\alpha_3 + 4\alpha_4 +  3\alpha_5 + 2\alpha_6 + \alpha_7}.
 \end{split}
\end{equation*}
Since $(\alpha_2 + \alpha_4, \alpha_3 + \alpha_4 + \alpha_5) = -1$,
using conjugation with $s_{\alpha_3 + \alpha_4 + \alpha_5}$ and
by Proposition \ref{prop_fact_1}, we obtain
\begin{equation*}
  s_{\alpha_2 + \alpha_4} s_{\alpha_3 + \alpha_4 + \alpha_5} =
     s_{\alpha_3 + \alpha_4 + \alpha_5}s_{\alpha_2 + \alpha_3 + 2\alpha_4 + \alpha_5}
\end{equation*}
and we get the following:
\begin{equation}
 \label{eq_reduce_1_E8_4}
 \begin{split}
  w_0 =
      & s_{\alpha_5 + \alpha_6}s_{\alpha_4 + \alpha_5 + \alpha_6 + \alpha_7}
        s_{\alpha_1 + 3\alpha_2 + 3\alpha_3 + 5\alpha_4 +
                   4\alpha_5 +3\alpha_6+ 2\alpha_7 + \alpha_8} \times \\
      & s_1 s_{\alpha_1 + \alpha_3 + \alpha_4 + \alpha_5 + \alpha_6 + \alpha_7 + \alpha_8}
    s_{\alpha_3 + \alpha_4 + \alpha_5} s_{\alpha_2 + \alpha_3 + 2\alpha_4 + \alpha_5}  \times \\
      &   s_{\alpha_2 + \alpha_3 + 2\alpha_4 +
               \alpha_5 + \alpha_6}
        s_{\alpha_2 + \alpha_3 + \alpha_4 + \alpha_5 +\alpha_6 + \alpha_7}
             s_{\alpha_7 + \alpha_8} s_6 s_{\alpha_2 + \alpha_4 + \alpha_5} \times \\
      & s_{\alpha_1 + \alpha_3 + \alpha_4 + \alpha_5 + \alpha_6 + \alpha_7+\alpha_8}
         s_{2\alpha_1 + 2\alpha_2 + 3\alpha_3 + 4\alpha_4 +  3\alpha_5 + 2\alpha_6 + \alpha_7}.
 \end{split}
\end{equation}
Since  $(\alpha_2 + \alpha_3 + 2\alpha_4 + \alpha_5, \alpha_2 + \alpha_3 + 2\alpha_4 + \alpha_5 + \alpha_6) = 1$, we have
\begin{equation*}
  s_{\alpha_2 + \alpha_3 + 2\alpha_4 + \alpha_5}
     s_{\alpha_2 + \alpha_3 + 2\alpha_4 + \alpha_5 + \alpha_6} =
   s_{\alpha_6}s_{\alpha_2 + \alpha_3 + 2\alpha_4 + \alpha_5}
\end{equation*}
Then, \eqref{eq_reduce_1_E8_4} is transformed as follows:
\begin{equation}
 \label{eq_reduce_1_E8_5}
 \begin{split}
  w_0 =
      & s_{\alpha_5 + \alpha_6}s_{\alpha_4 + \alpha_5 + \alpha_6 + \alpha_7}
        s_{\alpha_1 + 3\alpha_2 + 3\alpha_3 + 5\alpha_4 +
                   4\alpha_5 +3\alpha_6+ 2\alpha_7 + \alpha_8} \times \\
      & s_1 s_{\alpha_1 + \alpha_3 + \alpha_4 + \alpha_5 + \alpha_6 + \alpha_7 + \alpha_8}
 s_{\alpha_3 + \alpha_4 + \alpha_5} s_6 \big ( s_{\alpha_2 + \alpha_3 + 2\alpha_4 + \alpha_5} \times \\
      & s_{\alpha_2 + \alpha_3 + \alpha_4 + \alpha_5 +\alpha_6 + \alpha_7}
             s_{\alpha_7 + \alpha_8} \big ) s_6 s_{\alpha_2 + \alpha_4 + \alpha_5} \times \\
      & s_{\alpha_1 + \alpha_3 + \alpha_4 + \alpha_5 + \alpha_6 + \alpha_7+\alpha_8}
         s_{2\alpha_1 + 2\alpha_2 + 3\alpha_3 + 4\alpha_4 +  3\alpha_5 + 2\alpha_6 + \alpha_7}.
 \end{split}
\end{equation}

\underline{Conjugation with $s_{\alpha_6}$}.
In \eqref{eq_reduce_1_E8_5}, we apply conjugation with $s_{\alpha_6}$  to the
following factors:
\begin{equation*}
    s_{\alpha_2 + \alpha_3 + 2\alpha_4 + \alpha_5}
       s_{\alpha_2 + \alpha_3 + \alpha_4 + \alpha_5 +\alpha_6 + \alpha_7}
             s_{\alpha_7 + \alpha_8}.
\end{equation*}
Since $s_{\alpha_6}$ and $s_{\alpha_2 + \alpha_3 + \alpha_4 + \alpha_5 +\alpha_6 + \alpha_7}$ commute,
after conjugation we get:
\begin{equation}
 \label{eq_reduce_1_E8_6}
 \begin{split}
  w_0 =
      & s_{\alpha_5 + \alpha_6}s_{\alpha_4 + \alpha_5 + \alpha_6 + \alpha_7}
        s_{\alpha_1 + 3\alpha_2 + 3\alpha_3 + 5\alpha_4 +
                   4\alpha_5 +3\alpha_6+ 2\alpha_7 + \alpha_8} \times \\
      & s_1 s_{\alpha_1 + \alpha_3 + \alpha_4 + \alpha_5 + \alpha_6 + \alpha_7 + \alpha_8}
 s_{\alpha_3 + \alpha_4 + \alpha_5} s_{\alpha_2 + \alpha_3 + 2\alpha_4 + \alpha_5 +\alpha_6} \times \\
      & s_{\alpha_2 + \alpha_3 + \alpha_4 + \alpha_5 +\alpha_6 + \alpha_7}
             s_{\alpha_6 + \alpha_7 + \alpha_8}  s_{\alpha_2 + \alpha_4 + \alpha_5} \times \\
      & s_{\alpha_1 + \alpha_3 + \alpha_4 + \alpha_5 + \alpha_6 + \alpha_7+\alpha_8}
         s_{2\alpha_1 + 2\alpha_2 + 3\alpha_3 + 4\alpha_4 +  3\alpha_5 + 2\alpha_6 + \alpha_7}.
 \end{split}
\end{equation}
As a result, there are $12$ factors in \eqref{eq_reduce_1_E8_6}.
~\\

\subsection{Reducing $w_0$ from $12$ factors to $10$}
In  \eqref{eq_reduce_1_E8_6}, roots $\alpha_3 + \alpha_4 + \alpha_5$ and
$\alpha_1 + \alpha_3 + \alpha_4 + \alpha_5 + \alpha_6 + \alpha_7 + \alpha_8$ are orthogonal,
so we can swap corresponding reflections:
\begin{equation}
 \label{eq_reduce_1_E8_7}
 \begin{split}
  w_0 =
      & s_{\alpha_5 + \alpha_6}s_{\alpha_4 + \alpha_5 + \alpha_6 + \alpha_7}
        s_{\alpha_1 + 3\alpha_2 + 3\alpha_3 + 5\alpha_4 +
                   4\alpha_5 +3\alpha_6+ 2\alpha_7 + \alpha_8} \times \\
      & s_1 s_{\alpha_3 + \alpha_4 + \alpha_5} s_{\alpha_1 + \alpha_3 + \alpha_4 + \alpha_5 + \alpha_6 + \alpha_7 + \alpha_8}
 \big ( s_{\alpha_2 + \alpha_3 + 2\alpha_4 + \alpha_5 +\alpha_6} \times \\
      & s_{\alpha_2 + \alpha_3 + \alpha_4 + \alpha_5 +\alpha_6 + \alpha_7}
             s_{\alpha_6 + \alpha_7 + \alpha_8}  s_{\alpha_2 + \alpha_4 + \alpha_5} \big )
             \times \\
      & s_{\alpha_1 + \alpha_3 + \alpha_4 + \alpha_5 + \alpha_6 + \alpha_7+\alpha_8}
         s_{2\alpha_1 + 2\alpha_2 + 3\alpha_3 + 4\alpha_4 +  3\alpha_5 + 2\alpha_6 + \alpha_7}.
 \end{split}
\end{equation}

\underline{Conjugation with $s_{\alpha_1 + \alpha_3 + \alpha_4 + \alpha_5 + \alpha_6 + \alpha_7 + \alpha_8}$}.
In \eqref{eq_reduce_1_E8_7},
we apply conjugation with
 $s_{\alpha_1 + \alpha_3 + \alpha_4 + \alpha_5 + \alpha_6 + \alpha_7 + \alpha_8}$  
 to the following $4$ factors:
\begin{equation*}
    s_{\alpha_2 + \alpha_3 + 2\alpha_4 + \alpha_5 +\alpha_6}
       s_{\alpha_2 + \alpha_3 + \alpha_4 + \alpha_5 +\alpha_6 + \alpha_7}
             s_{\alpha_6 + \alpha_7 + \alpha_8}  s_{\alpha_2 + \alpha_4 + \alpha_5}.
\end{equation*}
Since
\begin{equation*}
 \footnotesize
 \begin{split}
     & (\alpha_2 + \alpha_3 + 2\alpha_4 + \alpha_5 +\alpha_6,  \;\;
        \alpha_1 + \alpha_3 + \alpha_4 + \alpha_5 + \alpha_6 + \alpha_7+\alpha_8) = -1, \\
     & (\alpha_2 + \alpha_3 + \alpha_4 + \alpha_5 +\alpha_6 + \alpha_7, \;\;
        \alpha_1 + \alpha_3 + \alpha_4 + \alpha_5 + \alpha_6 + \alpha_7+\alpha_8) = -1, \\
     & (\alpha_6 + \alpha_7 + \alpha_8, \;\;
       \alpha_1 + \alpha_3 + \alpha_4 + \alpha_5 + \alpha_6 + \alpha_7+\alpha_8) = 1, \\
     & (\alpha_2 + \alpha_4 + \alpha_5, \;\;
        \alpha_1 + \alpha_3 + \alpha_4 + \alpha_5 + \alpha_6 + \alpha_7+\alpha_8) = -1, \\
 \end{split}
\end{equation*}
then
\begin{equation*}
  \begin{split}
   & s_{\alpha_2 + \alpha_3 + 2\alpha_4 + \alpha_5 +\alpha_6}
     s_{\alpha_1 + \alpha_3 + \alpha_4 + \alpha_5 + \alpha_6 + \alpha_7+\alpha_8} = \\
     & \qquad \qquad
     s_{\alpha_1 + \alpha_3 + \alpha_4 + \alpha_5 + \alpha_6 + \alpha_7+\alpha_8}
  s_{\alpha_1 + \alpha_2 +2\alpha_3 + 3\alpha_4 + 2\alpha_5+2\alpha_6+\alpha_7+\alpha_8}, \\
     %% end of 1
     & s_{\alpha_2 + \alpha_3 + \alpha_4 + \alpha_5 +\alpha_6 + \alpha_7}
  s_{\alpha_1 + \alpha_3 + \alpha_4 + \alpha_5 + \alpha_6 + \alpha_7+\alpha_8} = \\
     & \qquad \qquad
    s_{\alpha_1 + \alpha_3 + \alpha_4 + \alpha_5 + \alpha_6 + \alpha_7+\alpha_8}
    s_{\alpha_1 + \alpha_2 + 2\alpha_3 + 2\alpha_4 + 2\alpha_5 + 2\alpha_6 + 2\alpha_7+\alpha_8}, \\
      %% end of 2
     & s_{\alpha_6 + \alpha_7 + \alpha_8} s_{\alpha_1 + \alpha_3 + \alpha_4 + \alpha_5 + \alpha_6 + \alpha_7+\alpha_8} = \\
     & \qquad \qquad
       s_{\alpha_1 + \alpha_3 + \alpha_4 + \alpha_5 + \alpha_6 + \alpha_7+\alpha_8}
       s_{\alpha_1 +\alpha_3 + \alpha_4 + \alpha_5}, \\
      %% end of 3
     & s_{\alpha_2 + \alpha_4 + \alpha_5}
       s_{\alpha_1 + \alpha_3 + \alpha_4 + \alpha_5 + \alpha_6 + \alpha_7+\alpha_8} = \\
     & \qquad \qquad
 s_{\alpha_1 + \alpha_3 + \alpha_4 + \alpha_5 + \alpha_6 + \alpha_7+\alpha_8}
 s_{\alpha_1+\alpha_2+\alpha_3 + 2\alpha_4 + 2\alpha_5 + \alpha_6 + \alpha_7+\alpha_8}.\\
      %%  end of 4
  \end{split}
\end{equation*}
After conjugation with $s_{\alpha_1 + \alpha_3 + \alpha_4 + \alpha_5 + \alpha_6 + \alpha_7+\alpha_8}$, we get:
\begin{equation}
 \label{eq_reduce_1_E8_8}
 \begin{split}
  w_0 =
    & s_{\alpha_5 + \alpha_6}s_{\alpha_4 + \alpha_5 + \alpha_6 + \alpha_7}
        s_{\alpha_1 + 3\alpha_2 + 3\alpha_3 + 5\alpha_4 +
                   4\alpha_5 +3\alpha_6+ 2\alpha_7 + \alpha_8} \times \\
    &  s_1 s_{\alpha_3 + \alpha_4 + \alpha_5}
     s_{\alpha_1 + \alpha_2 +2\alpha_3 + 3\alpha_4 + 2\alpha_5+2\alpha_6+\alpha_7+\alpha_8}
        \times \\
    & s_{\alpha_1 + \alpha_2 + 2\alpha_3 + 2\alpha_4 + 2\alpha_5 + 2\alpha_6 + 2\alpha_7+\alpha_8}  \times \\
   & s_{\alpha_1 +\alpha_3 + \alpha_4 + \alpha_5}
     s_{\alpha_1+\alpha_2+\alpha_3 + 2\alpha_4 + 2\alpha_5 + \alpha_6 + \alpha_7+\alpha_8}
       \times \\
     & s_{2\alpha_1 + 2\alpha_2 + 3\alpha_3 + 4\alpha_4 +  3\alpha_5 + 2\alpha_6 + \alpha_7}.
 \end{split}
\end{equation}
So, we have $10$ factors in \eqref{eq_reduce_1_E8_8}.
~\\

\subsection{Reducing $w_0$ from $10$ factors to $8$} First, we have that 
\begin{equation*}
    s_{\alpha_1} s_{\alpha_3 + \alpha_4 + \alpha_5} = s_{\alpha_3 + \alpha_4 + \alpha_5}
         s_{\alpha_1 + \alpha_3 + \alpha_4 + \alpha_5},
\end{equation*}

\underline{Conjugation with $s_{\alpha_1 + \alpha_3 + \alpha_4 + \alpha_5}$}.
By \eqref{eq_reduce_1_E8_8}, we get the following representation of $w_0$:
\begin{equation}
 \label{eq_reduce_1_E8_9}
 \begin{split}
  w_0 =
    & s_{\alpha_5 + \alpha_6}s_{\alpha_4 + \alpha_5 + \alpha_6 + \alpha_7}
        s_{\alpha_1 + 3\alpha_2 + 3\alpha_3 + 5\alpha_4 +
                   4\alpha_5 +3\alpha_6+ 2\alpha_7 + \alpha_8} \times \\
    &  s_{\alpha_3 + \alpha_4 + \alpha_5} s_{\alpha_1 + \alpha_3 + \alpha_4 + \alpha_5}
     \big (
     s_{\alpha_1 + \alpha_2 +2\alpha_3 + 3\alpha_4 + 2\alpha_5+2\alpha_6+\alpha_7+\alpha_8}
        \times \\
    & s_{\alpha_1 + \alpha_2 + 2\alpha_3 + 2\alpha_4 + 2\alpha_5 + 2\alpha_6 + 2\alpha_7+\alpha_8} \big ) s_{\alpha_1 +\alpha_3 + \alpha_4 + \alpha_5} \times \\
   & s_{\alpha_1+\alpha_2+\alpha_3 + 2\alpha_4 + 2\alpha_5 + \alpha_6 + \alpha_7+\alpha_8}
       %% \times \\
      s_{2\alpha_1 + 2\alpha_2 + 3\alpha_3 + 4\alpha_4 +  3\alpha_5 + 2\alpha_6 + \alpha_7}.
 \end{split}
\end{equation}
~\\
It is easy to check that
\begin{equation*}
 \footnotesize
 \begin{split}
   & \alpha_1 + \alpha_3 + \alpha_4 + \alpha_5 \; \perp \;
      \alpha_1 + \alpha_2 +2\alpha_3 + 3\alpha_4 + 2\alpha_5+2\alpha_6+\alpha_7+\alpha_8, \\
   &  \alpha_1 + \alpha_3 + \alpha_4 + \alpha_5 \; \perp \;
      \alpha_1 + \alpha_2 + 2\alpha_3 + 2\alpha_4 + 2\alpha_5 + 2\alpha_6 + 2\alpha_7+\alpha_8.
 \end{split}
\end{equation*}
Then, in \eqref{eq_reduce_1_E8_9}, the pair of elements
$s_{\alpha_1 + \alpha_3 + \alpha_4 + \alpha_5}$ is reduced:

\begin{equation}
 \label{eq_reduce_1_E8_10}
 \begin{split}
  w_0 =
    & s_{\alpha_5 + \alpha_6}s_{\alpha_4 + \alpha_5 + \alpha_6 + \alpha_7}
        s_{\alpha_1 + 3\alpha_2 + 3\alpha_3 + 5\alpha_4 +
                   4\alpha_5 +3\alpha_6+ 2\alpha_7 + \alpha_8} \times \\
    &  s_{\alpha_3 + \alpha_4 + \alpha_5}
     s_{\alpha_1 + \alpha_2 +2\alpha_3 + 3\alpha_4 + 2\alpha_5+2\alpha_6+\alpha_7+\alpha_8}
        \times \\
    & s_{\alpha_1 + \alpha_2 + 2\alpha_3 + 2\alpha_4 + 2\alpha_5 + 2\alpha_6 + 2\alpha_7+\alpha_8}  \times \\
   & s_{\alpha_1+\alpha_2+\alpha_3 + 2\alpha_4 + 2\alpha_5 + \alpha_6 + \alpha_7+\alpha_8}
      s_{2\alpha_1 + 2\alpha_2 + 3\alpha_3 + 4\alpha_4 +  3\alpha_5 + 2\alpha_6 + \alpha_7}.
 \end{split}
\end{equation}
We get $8$ factors in \eqref{eq_reduce_1_E8_10}.

\subsection{Getting the factor with index $\alpha_{max}$}
We strive to choose the reflection $s_\tau$ for conjugation of factors in \eqref{eq_reduce_1_E8_10}
in such a way that one of the resulting factors is the highest root $\alpha_{max}$.
To this end, let us put
\begin{equation*}
 %% \label{eq_reduce_1_E8_11}
 \tau = \alpha_1 + 2\alpha_2 + 3\alpha_3 + 4\alpha_4 + 3\alpha_5 +
                                      3\alpha_6 + 2\alpha_7 + \alpha_8.
\end{equation*}

Let us find the inner products of $\tau$ with roots of \eqref{eq_reduce_1_E8_10}.
\begin{equation}
 \footnotesize
 \label{eq_reduce_1_E8_12}
 \begin{split}
    & (1) \;\;  \alpha_5 + \alpha_6 \perp \tau, \\
    & (2) \;\; \alpha_4 + \alpha_5 + \alpha_6 + \alpha_7 \perp \tau, \\
    & (3) \;\; (\alpha_1 + 3\alpha_2 + 3\alpha_3 + 5\alpha_4 +
             4\alpha_5 +3\alpha_6+ 2\alpha_7 + \alpha_8, \;\; \tau) =  1, \\
    & (4) \;\; \alpha_3 + \alpha_4 + \alpha_5 \perp \tau,  \\
    & (5) \;\; (\alpha_1 + \alpha_2 +2\alpha_3 + 3\alpha_4 +
               2\alpha_5+2\alpha_6+\alpha_7+\alpha_8, \;\; \tau) = 1, \\
    & (6) \;\; (\alpha_1 + \alpha_2 + 2\alpha_3 + 2\alpha_4 + 2\alpha_5 +
          2\alpha_6 + 2\alpha_7+\alpha_8, \;\; \tau) = 1, \\
    & (7) \;\; (\alpha_1+\alpha_2+\alpha_3 + 2\alpha_4 + 2\alpha_5 +
        \alpha_6 + \alpha_7+\alpha_8, \;\; \tau) = -1, \\
    & (8) \;\; 2\alpha_1 + 2\alpha_2 + 3\alpha_3 + 4\alpha_4 +
                 3\alpha_5 + 2\alpha_6 + \alpha_7 \perp \tau.
   \end{split}
\end{equation}

\underline{Conjugation with $s_\tau$}.
The roots $(1),(2),(4),(8)$ in \eqref{eq_reduce_1_E8_12} do not change,
so the corresponding factors do not change after conjugation with $s_\tau$.
The roots $(3),(5),(6)$  are reduced by $\tau$ after conjugation.
In case $(7)$, after conjugation, the root increases by $\tau$.
\begin{equation}
 \label{eq_reduce_1_E8_13}
 \footnotesize
 \begin{split}
   & (3) \;\; (\alpha_1 + 3\alpha_2 + 3\alpha_3 + 5\alpha_4 +
             4\alpha_5 +3\alpha_6+ 2\alpha_7 + \alpha_8) - \tau =
               \alpha_2 +\alpha_4 +\alpha_5, \\
   & (5) \;\; (\alpha_1 + \alpha_2 +2\alpha_3 + 3\alpha_4 +
               2\alpha_5+2\alpha_6+\alpha_7+\alpha_8) - \tau = \\
       & \qquad \qquad   -(\alpha_2 + \alpha_3 + \alpha_4 + \alpha_5 + \alpha_6 + \alpha_7), \\
   & (6) \;\; (\alpha_1 + \alpha_2 + 2\alpha_3 + 2\alpha_4 + 2\alpha_5 +
          2\alpha_6 + 2\alpha_7+\alpha_8) - \tau = \\
         & \qquad \qquad  -(\alpha_2 + \alpha_3 +2\alpha_4 + \alpha_5 + \alpha_6), \\
   & (7) \;\; (\alpha_1+\alpha_2+\alpha_3 + 2\alpha_4 + 2\alpha_5 +
                \alpha_6 + \alpha_7+\alpha_8) + \tau = \\
   &  \qquad \qquad 2\alpha_1 +3\alpha_2+ 4\alpha_3 + 6\alpha_4 + 5\alpha_5 +
            4\alpha_6 + 3\alpha_7 + 2\alpha_8.
  \end{split}
\end{equation}
 The root (7) in \eqref{eq_reduce_1_E8_13} is the highest root $\alpha_{max}$ in $E_8$,
 The root (8) in \eqref{eq_reduce_1_E8_12} is the highest root $\alpha^{e7}_{max}$ in
 $E_7 \subset E_8$, see \eqref{eq_not_E8_1} and Table \ref{tab_highest_roots_EFG}.
 The decomposition \eqref{eq_reduce_1_E8_10} is transformed as follows:
\begin{equation}
 \label{eq_reduce_1_E8_14}
 \begin{split}
     w_0 =
    & s_{\alpha_5 + \alpha_6}s_{\alpha_4 + \alpha_5 + \alpha_6 + \alpha_7}
        s_{\alpha_2 +\alpha_4 +\alpha_5} s_{\alpha_3 + \alpha_4 + \alpha_5}
        \times \\
    & s_{\alpha_2 + \alpha_3 + \alpha_4 + \alpha_5 + \alpha_6 + \alpha_7}
      s_{\alpha_2 + \alpha_3 + 2\alpha_4 + \alpha_5 + \alpha_6}
     s_{\alpha_{max}}s_{\alpha^{e7}_{max}}. 
     %%  s_{\scriptscriptstyle {e7}, max}.
 \end{split}
\end{equation}

\subsection{Getting factors with indices $\alpha^{d6}_{max}$,
$\alpha^{e7}_{max}$, $\alpha_{max}$}
Note that $\alpha_{max}$ is orthogonal to all simple roots except for $\alpha_8$,
in addition $\alpha^{e7}_{max}$ is orthogonal to all simple roots except for $\alpha_1$ and
$\alpha_8$. Thus, the conjugation associated with any root not containing
$\alpha_1$ and $\alpha_8$ preserves $\alpha_{max}$ and
$\alpha^{e7}_{max}$.  We continue to apply conjugations.
Let us put $\nu = \alpha_4 + \alpha_5 + \alpha_6$.

Since $\nu$ is orthogonal to $\alpha_2 +\alpha_4 +\alpha_5$, ~$\alpha_3 + \alpha_4 + \alpha_5$, 
~$\alpha^{e7}_{max}$ and  $\alpha_{max}$, the conjugation with $s_\nu$ does not change these roots.
Since $(\alpha_5 + \alpha_6, \nu) = 1$, $(\alpha_4 + \alpha_5 + \alpha_6 + \alpha_7, \nu) = 1$
and $(\alpha_2 + \alpha_3 + 2\alpha_4 + \alpha_5 + \alpha_6, \nu) = 1$, then

\begin{equation*}
 %%\label{eq_reduce_1_E8_15}
 \footnotesize
 \begin{array}{lll}
     & \alpha_5 + \alpha_6& \Rightarrow (\alpha_5 + \alpha_6) - \nu = -\alpha_4, \\
     & \alpha_4 + \alpha_5 + \alpha_6 + \alpha_7  & \Rightarrow  (\alpha_4 + \alpha_5 + \alpha_6 + \alpha_7) - \nu = \alpha_7, \\
     & \alpha_2 + \alpha_3 + 2\alpha_4 + \alpha_5 + \alpha_6 & \Rightarrow \\
     &  \qquad    \alpha_2 + \alpha_3 + 2\alpha_4 + \alpha_5 + \alpha_6 - \nu 
            & = \alpha_2 + \alpha_3 + \alpha_4. \\
 \end{array}
\end{equation*}

There remains only one root in \eqref{eq_reduce_1_E8_14} which must be considered.
It is modified by conjugation $s_\nu$ as follows:
\begin{equation*}
 %%\label{eq_obtained_1_E8_D6}
 \footnotesize
 \begin{split}
     & \alpha_2 + \alpha_3 + \alpha_4 + \alpha_5 + \alpha_6 + \alpha_7 \qquad \Rightarrow  \\
     &  \qquad (\alpha_2 + \alpha_3 + \alpha_4 + \alpha_5 + \alpha_6 + \alpha_7) + \nu = 
           \alpha_2 + \alpha_3 + 2\alpha_4 + 2\alpha_5 + 2\alpha_6 + \alpha_7. \\
 \end{split}
\end{equation*}
The obtained root is $s_{\alpha^{d6}_{max}}$.
~\\

\underline{Conjugation with $s_\nu$}. After conjugation with 
$s_{\alpha_4 + \alpha_5 + \alpha_6}$ we get the following decomposition of $w_0$:
\begin{equation}
 \label{eq_reduce_1_E8_16}
 \begin{split}
     w_0 =
    & s_{\alpha_4} s_{\alpha_7}s_{\alpha_2 +\alpha_4 +\alpha_5}
        s_{\alpha_3 + \alpha_4 + \alpha_5} s_{\alpha^{d6}_{max}} s_{\alpha_2 + \alpha_3 + \alpha_4} 
        s_{\alpha_{max}}s_{\alpha^{e7}_{max}}.
 \end{split}
\end{equation}
where $s_{\alpha^{d6}_max}$ (resp. $s_{\alpha^{e7}_{max}}$)
is the reflection corresponding to the highest
root in $D_6 \subset E_7 \subset E_8$ (resp. $E_7 \subset E_8$).
The roots corresponding to factors \eqref{eq_reduce_1_E8_16} are as follows:
\begin{equation}
 \label{eq_final_roots}
 \begin{split}
   & \alpha_4, \;\; \alpha_7, \;\; \alpha_2 + \alpha_4 + \alpha_5,
     \;\; \alpha_2 + \alpha_3 + \alpha_4, \;\; \alpha_3 + \alpha_4 + \alpha_5, \\
   & \alpha^{d6}_{max} =
      \alpha_2 + \alpha_3 + 2\alpha_4 + 2\alpha_5 + 2\alpha_6 + \alpha_7, \\
   & \alpha^{e7}_{max} = 2\alpha_1 + 2\alpha_2 + 3\alpha_3 + 4\alpha_4 +
                 3\alpha_5 + 2\alpha_6 + \alpha_7, \\
   & \alpha_{max} = 2\alpha_1 +3\alpha_2+ 4\alpha_3 + 6\alpha_4 + 5\alpha_5 +
            4\alpha_6 + 3\alpha_7 + 2\alpha_8,
 \end{split}
\end{equation}
see Table \ref{tab_factoriz_all}.

\subsection{Finally: max-orthogonal decomposition in $W(E_8)$}

\begin{proposition}[decomposition in $W(E_8)$]
  \label{prop_factoriz_E8}
 {\rm{(i)}}
 The roots corresponding to factors \eqref{eq_final_roots}
 are mutually orthogonal.

{\rm{(ii)}}
The element $w_0$ is decomposed as the product of $8$ 
reflections associated with the \underline{mutually orthogonal} roots:
\begin{equation}
  \label{eq_final_factoriz_E8}
 \begin{split}
    \alpha_2, \; \alpha_3, \; \alpha_5, \; \alpha_7, \;
    \alpha^{d4}_{max}, \; \alpha^{d6}_{max}, \; \alpha^{e7}_{max}, \; \alpha_{max}.
 \end{split}
\end{equation}
The highest roots in \eqref{eq_final_factoriz_E8} are ordered as follows:
\begin{equation*}
     \alpha^{d4}_{max} < \alpha^{d6}_{max} < \alpha^{e7}_{max} < \alpha_{max}.  \\
\end{equation*}

{\rm{(iii)}}
The element $w_0$ is the longest element in $W(E_8)$.
\end{proposition}

\PerfProof
(i) The root $\alpha_{max}$
orthogonal to the remaining roots in \eqref{eq_final_roots} since they does not contain
the term $\alpha_8$. Exclude $\alpha_{max}$. The root $\alpha^{e7}_{max}$ orthogonal
to the remaining roots since  they does not contain the term $\alpha_1$.
Exclude $\alpha^{e7}_{max}$. The root $\alpha^{d6}_{max}$
orthogonal to the remaining roots since they does not contain the term $\alpha_6$.
It is easy to check that all the remaining roots are mutually orthogonal.

(ii)  Let us conjugate factors \eqref{eq_reduce_1_E8_16} with $s_{\alpha_3 + \alpha_4}$.
This conjugation acts trivially on
$s_7$, $\alpha^{d6}_{max}$, $\alpha^{e7}_{max}$ and
$\alpha_{max}$, since $\alpha_3 + \alpha_4$ is orthogonal to all of them.
Further,
\begin{equation*}
 \begin{split}
    & s_{\alpha_3 + \alpha_4} s_{\alpha_4} s_{\alpha_3 + \alpha_4} = s_{\alpha_3}, \\
    & s_{\alpha_3 + \alpha_4} s_{\alpha_2+\alpha_4 + \alpha_5}  s_{\alpha_3 + \alpha_4} =
                 s_{\alpha_2 + \alpha_3 + 2\alpha_4+\alpha_5} = s_{\alpha^{d4}_{max}}, \\
    & s_{\alpha_3 + \alpha_4} s_{\alpha_2+\alpha_3 + \alpha_4}  s_{\alpha_3 + \alpha_4} =
                       s_{\alpha_2} \\
    & s_{\alpha_3 + \alpha_4} s_{\alpha_3 + \alpha_4 + \alpha_5}  s_{\alpha_3 + \alpha_4} =
                       s_{\alpha_5}.
 \end{split}
\end{equation*}
Thus, we get all factors of \eqref{eq_final_factoriz_E8}.

(iii) Let us construct the determinant by roots of \eqref{eq_final_factoriz_E8},
see (\ref{eq_determ_E8}a). Next, we transform its rows as follows:
\begin{equation*}
 \begin{array}{llll}
   & \alpha_{max} \rightarrow \alpha_{max} - \alpha^{e7}_{max} - \alpha^{d6}_{max}
     - \alpha_7 & = & (0,0,0,0,0,0,0,2), \\
     & \\
   & \alpha^{e7}_{max} \rightarrow \alpha^{e7}_{max} - \alpha^{d6}_{max} - \alpha^{d4}_{max}
      - \alpha_3 & = & (2,0,0,0,0,0,0,0), \\
     & \\
   & \alpha^{d6}_{max} \rightarrow \alpha^{d6}_{max} - \alpha^{d4}_{max}
      - \alpha_7 - \alpha_5 & = & (0,0,0,0,0,2,0,0), \\
     & \\
   & \alpha^{d4}_{max} \rightarrow \alpha^{d4}_{max} - \alpha_2 - \alpha_3 - \alpha_5  & =
      & (0,0,0,2,0,0,0,0). \\
 \end{array}
\end{equation*}
As a result, we obtain a non-degenerate determinant, see (\ref{eq_determ_E8}b).
\begin{equation}
  \label{eq_determ_E8}
 \footnotesize
 \begin{split}
  &\left | \left |
  \begin{array}{llllllll}
     0 & 1 & 0 & 0 & 0 & 0 & 0 & 0 \\
     0 & 0 & 1 & 0 & 0 & 0 & 0 & 0 \\
     0 & 0 & 0 & 0 & 1 & 0 & 0 & 0 \\
     0 & 0 & 0 & 0 & 0 & 0 & 1 & 0 \\
     0 & 1 & 1 & 2 & 1 & 0 & 0 & 0 \\
     0 & 1 & 1 & 2 & 2 & 2 & 1 & 0 \\
     2 & 2 & 3 & 4 & 3 & 2 & 1 & 0 \\
     2 & 3 & 4 & 6 & 5 & 4 & 3 & 2 \\
  \end{array}
  \right | \right |
  \begin{array}{l}
     \alpha_2 \\
     \alpha_3 \\
     \alpha_5 \\
     \alpha_7 \\
     \alpha^{d4}_{max} \\
     \alpha^{d6}_{max} \\
     \alpha^{e7}_{max} \\
     \alpha_{max} \\
  \end{array}
     \Longrightarrow
  \left | \left |
  \begin{array}{llllllll}
     0 & 1 & 0 & 0 & 0 & 0 & 0 & 0 \\
     0 & 0 & 1 & 0 & 0 & 0 & 0 & 0 \\
     0 & 0 & 0 & 0 & 1 & 0 & 0 & 0 \\
     0 & 0 & 0 & 0 & 0 & 0 & 1 & 0 \\
     0 & 0 & 0 & 2 & 0 & 0 & 0 & 0 \\
     0 & 0 & 0 & 0 & 0 & 2 & 0 & 0 \\
     2 & 0 & 0 & 0 & 0 & 0 & 0 & 0 \\
     0 & 0 & 0 & 0 & 0 & 0 & 0 & 2 \\
  \end{array}
  \right | \right | \\
  & \qquad \qquad \qquad \quad (a)
    \qquad \qquad \qquad \qquad \qquad \qquad \qquad \qquad \qquad (b)
  \end{split}
\end{equation}

Thus, the roots of \eqref{eq_final_factoriz_E8} are linearly independent.
Since they are mutually orthogonal, $w_0$ acts as $-1$ on each of these roots,
and hence $w_0$ acts as $-1$ on the whole space. Therefore, $w_0$ acts in the same way
as the longest element. Since the longest element is unique in the Weyl group,
$w_0$ is the longest in $W(E_8)$. \qed
~\\

\section{\bf Decomposition in $W(F_4)$}

\begin{figure}[h]
\centering
   \includegraphics[scale=0.4]{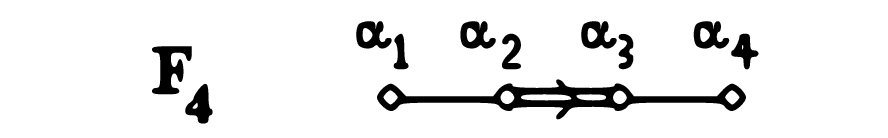}
\caption{\footnotesize{Numbering of simple roots in $F_4$}}
%%%%%% The label must come after caption
\label{fig_F4}
\end{figure}

We introduce the notation for the \underline{highest roots}
in the root subsystems $C_2 \subset C_3 \subset F_4$ as follows:
\begin{equation}
  \label{eq_not_F4_1}
  \begin{split}
     & \alpha^{c2}_{max} := \alpha_2 + 2\alpha_3, \\
     & \alpha^{c3}_{max} := \alpha_2 + 2\alpha_3 + 2\alpha_4,
   \end{split}
\end{equation}
see Table \ref{tab_factoriz_all} and Remark \ref{rem_notations}.
The numbering of vertices of the Dynkin diagram $F_4$ is shown
in Fig. \ref{fig_F4}. Note that indices of simple roots for 
$\alpha^{c3}_{max}$ (resp. $\alpha^{c2}_{max}$)  in \eqref{eq_not_F4_1} differ from the indices
of simple roots that appears as summands in $\alpha^{c3}_{max}$ (resp. $\alpha^{c2}_{max}$)
defined in \eqref{eq_not_Bn_Cn_2}. %%

   As in the case of $E_6$, consider the following $2$ elements in the
Weyl group $W(F_4)$:
\begin{equation}
  \begin{split}
   & w_1 = s_1 s_2 (s_3 s_2 s_4 (s_3 s_2 s_1 s_2 s_3) s_4 s_2 s_3) s_2 s_1 \\
   & w_2 = s_4 s_3 s_2 s_3 s_4 s_3  s_2 s_3 s_2 \\
  \end{split}
\end{equation}

\begin{proposition}[decomposition in $W(F_4)$]
  \label{prop_factoriz_F4}
  {\rm{(i)}}
   The element $w_1$ is the reflection corresponding to the highest root
   $\alpha_{max}$:
\begin{equation}
   \label{eq_how_w1_looks}
   w_1 = s_{\alpha_{max}}.
\end{equation}

  {\rm{(ii)}}
  The element $w_2$ is the product of $3$ reflections corresponding to mutually orthogonal
  positive roots $s_{\alpha_2}$, $s_{\alpha_2 + 2\alpha_3}$ and
  $s_{\alpha_2 + 2\alpha_3 + 2\alpha_4}$:
\begin{equation}
   \label{eq_how_w1_looks_2}
   w_2 = s_{\alpha_2}s_{\alpha_2 + 2\alpha_3}s_{\alpha_2 + 2\alpha_3 + 2\alpha_4}.
\end{equation}

  {\rm{(iii)}}
   The element $w_0 = w_1w_2$ is the product of the following $4$ reflections:
\begin{equation}
  \label{eq_factoris_F4}
   w_0 =  s_{\alpha_2} s_{\alpha^{c2}_{max}} s_{\alpha^{c3}_{max}}s_{\alpha_{max}}.
\end{equation}
The roots corresponding to decomposition \eqref{eq_factoris_F4} are \underline{mutually orthogonal}.
They are ordered as follows:
\begin{equation*}
    \alpha_2 <  \alpha^{c2}_{max} < \alpha^{c3}_{max} < \alpha_{max}.  \\
\end{equation*}

{\rm{(iv)}}
   The element $w_0$ is the longest element in $W(F_4)$.

\end{proposition}

\PerfProof (i) 
For conjugation by reflections corresponding to long roots $\alpha_1$ and $\alpha_2$ 
(resp. short roots $\alpha_3$ and $\alpha_4$), 
we use Proposition \ref{prop_fact_1} (resp. Proposition \ref{prop_fact_2}):
\begin{equation*}
  \begin{split}
   & s_2 s_1 s_2 = s_{\alpha_1 + \alpha_2}, \\
   & s_3 s_2 s_1 s_2 s_3 = s_3 s_{\alpha_1 + \alpha_2} s_3 = s_{\alpha_2 + \alpha_1 + 2\alpha_3},\\
   & s_4s_3 s_2 s_1 s_2 s_3s_4 = s_4 s_{\alpha_2 + \alpha_1 + 2\alpha_3} s_4 =
                 s_{\alpha_1 + \alpha_2 +  2\alpha_3 + 2\alpha_4},\\
   & s_2s_4s_3 s_2 s_1 s_2 s_3s_4s_2 = s_2 s_{\alpha_2 + \alpha_1 + 2\alpha_3 + 2\alpha_4}s_2 =
                 s_{\alpha_1 + 2\alpha_2 +  2\alpha_3 + 2\alpha_4},\\
   & s_3s_2s_4s_3 s_2 s_1 s_2 s_3s_4s_2s_3 =
                 s_{\alpha_1 + 2\alpha_2 +  4\alpha_3 + 2\alpha_4,}\\
   & s_2s_3s_2s_4s_3 s_2 s_1 s_2 s_3s_4s_2s_3s_2 =
                 s_{\alpha_1 + 3\alpha_2 +  4\alpha_3 + 2\alpha_4},\\
   & s_1s_2s_3s_2s_4s_3 s_2 s_1 s_2 s_3s_4s_2s_3s_2s_1 =
                 s_{2\alpha_1 + 3\alpha_2 +  4\alpha_3 + 2\alpha_4}.\\
  \end{split}
\end{equation*}
Since $2\alpha_1 + 3\alpha_2 +  4\alpha_3 + 2\alpha_4$ is the highest root in $F_4$,
we get eq. \eqref{eq_how_w1_looks}.
~\\

(ii) Transform $w_2$ as follows:
\begin{equation*}
  \begin{split}
     & (s_4 s_3 s_2 s_3 s_4)s_3 s_2 s_3 s_2  =
       (s_4 s_{\alpha_2 + 2\alpha_3} s_4) s_3 s_2 s_3 s_2  = \\
     & s_{\alpha_2 + 2 \alpha_3 + 2\alpha_4}(s_3 s_2 s_3) s_2  =
       s_{\alpha_2 + 2 \alpha_3 + 2\alpha_4}s_{\alpha_2 + 2 \alpha_3} s_2.
  \end{split}
\end{equation*}
The orthogonality of the vectors $\alpha_2$, $\alpha_2 + 2 \alpha_3$
and $\alpha_2 + 2 \alpha_3 + 2\alpha_4$ is verified directly.
~\\

(iii) The root $\alpha_{max}$ is orthogonal to all simple root except $\alpha_1$.
So, $\alpha_{max}$  is orthogonal to $\alpha_2$, $\alpha_2 + 2 \alpha_3$
and $\alpha_2 + 2 \alpha_3 + 2\alpha_4$.
~\\

(iv) It suffices to prove that $w_0$ acts on simple roots
in the same way as the longest element. The reflection $s_{\alpha_{max}}$ preserves all
simple roots except $\alpha_1$.
Since $(\alpha_2, \alpha_2 + 2\alpha_3) = 0$ and
$(\alpha_2, \alpha_2 + 2\alpha_3 + 2\alpha_4) = 0$, we get
\begin{equation*}
 \label{F4_w0_on_alph2}
  w_0(\alpha_2) = s_{\alpha_2}(\alpha_2) = -\alpha_2.
\end{equation*}
Since $s_{\alpha_2 + 2\alpha_3}$ maps $\alpha_3$ to $-\alpha_2 -\alpha_3$  and preserves $\alpha_2$,
then
\begin{equation*}
    s_{\alpha_2 + 2\alpha_3}s_{\alpha_2}(\alpha_3) =
      s_{\alpha_2 + 2\alpha_3}(\alpha_2 + \alpha_3) = \alpha_2 -\alpha_2 -\alpha_3 =
        -\alpha_3.
\end{equation*}
Since $(\alpha_2 + 2\alpha_3 + 2\alpha_4, \alpha_3) = -1 + 2 -1 = 0$, we have
\begin{equation}
 \label{F4_w0_on_alph3}
   w_0(\alpha_3) = s_{\alpha_2 + 2\alpha_3 + 2\alpha_4}(-\alpha_3 ) = -\alpha_3.
\end{equation}
Further, $s_{\alpha_2}(\alpha_4) = \alpha_4$, and
$s_{\alpha_2 + 2\alpha_3}(\alpha_4) = \alpha_2 + 2\alpha_3 + \alpha_4$, see Remark \ref{rem_inner_prod}(i).
Then, since $(\alpha_2 + 2\alpha_3 + \alpha_4, \alpha_2 + 2\alpha_3 + 2\alpha_4) = 1$, we have
\begin{equation}
 \label{F4_w0_on_alph4}
 \begin{split}
  w_0(\alpha_4) & =
    s_{\alpha_2 + 2\alpha_3 + 2\alpha_4}(\alpha_2 + 2\alpha_3 + \alpha_4) = \\
  & (\alpha_2 + 2\alpha_3 + \alpha_4) - (\alpha_2 + 2\alpha_3 + 2\alpha_4) =
            -\alpha_4,\\
 \end{split}
\end{equation}
~\\
At last, consider $w_0(\alpha_1)$. First,
$\alpha_{max}(\alpha_1) = \alpha_1 - \alpha_{max}$.
Since $s_2$, $s_{\alpha_2 + 2\alpha_3}$, $s_{\alpha_2 + 2\alpha_3+ 2\alpha_4}$ preserve $\alpha_{max}$ we have
\begin{equation}
 \label{F4_w0_on_alph1}
 \begin{split}
   w_0(\alpha_1) \; = \; &
    s_{\alpha_2 + 2\alpha_3 + 2\alpha_4}s_{\alpha_2 + 2\alpha_3} s_2(\alpha_1) -
       \alpha_{max} = \\
  & s_{\alpha_2 + 2\alpha_3 + 2\alpha_4}s_{\alpha_2 + 2\alpha_3} (\alpha_1+\alpha_2) -
       \alpha_{max} = \\
  & s_{\alpha_2 + 2\alpha_3 + 2\alpha_4}s_{\alpha_2 + 2\alpha_3} (\alpha_1)
     + \alpha_2 - \alpha_{max} = \\
  & s_{\alpha_2 + 2\alpha_3 + 2\alpha_4}(\alpha_1 + \alpha_2 + 2\alpha_3)
     + \alpha_2 - \alpha_{max} = \\
  & s_{\alpha_2 + 2\alpha_3 + 2\alpha_4}(\alpha_1)
     + 2\alpha_2 + 2\alpha_3 - \alpha_{max} = \\
  & \alpha_1 + \alpha_2 + 2\alpha_3 + 2\alpha_4 + 2\alpha_2 + 2\alpha_3
               - \alpha_{max}  = \\
  & \alpha_1 + 3\alpha_2 + 4\alpha_3 + 2\alpha_4 - \alpha_{max}  = -\alpha_1.\\
 \end{split}
\end{equation}
By \eqref{F4_w0_on_alph2}, \eqref{F4_w0_on_alph3}, \eqref{F4_w0_on_alph4}
and \eqref{F4_w0_on_alph1} we get $w_0 = -1$, i.e., $w_0$ is the longest
element in $W(F_4)$, see \cite[Plate VIII]{Bo02}.
\qed
~\\

\section{\bf Decomposition in $W(G_2)$}

There are $6$ positive roots in $G_2$:
\begin{equation*}
  \begin{array}{lll}
      & \text{ short roots: }
              & \alpha_1, \;  \alpha_1 + \alpha_2, \; 2\alpha_1 + \alpha_2, \\
      & \text{ long roots:  }
              & \alpha_2, \;  3\alpha_1 + \alpha_2,\; 3\alpha_1 + 2\alpha_2. \\
   \end{array}
\end{equation*}

\begin{proposition}
  Consider the following element in the Weyl group $W(G_2)$:
\begin{equation}
 \label{eq_longest_G2}
     w_0 = s_2 s_1 s_2 s_1 s_2 s_1.
\end{equation}
  The element $w_0$ is the longest element in  $W(G_2)$.
\end{proposition}

\begin{figure}[h]
\centering
   \includegraphics[scale=0.4]{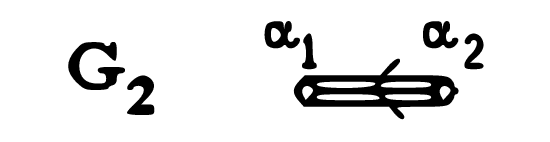}
\caption{\footnotesize{Numbering of simple roots in $G_2$}}
%%%%%% The label must come after caption
\label{fig_G2}
\end{figure}

\PerfProof
By Propositions \ref{prop_fact_1} and \ref{prop_fact_2} we have
\begin{equation}
  \label{eq_first_factoriz_G2}
  \begin{split}
    & s_2 s_1 s_2 = s_{\alpha_1 + \alpha_2}, \quad
      s_1 s_2 s_1 =  s_{3\alpha_1 + \alpha_2}, \\
    & w_0 = s_{\alpha_1 + \alpha_2} s_{3\alpha_1 + \alpha_2}
  \end{split}
\end{equation}
~\\
Consider the action of $w_0$ on $\alpha_1$ and $\alpha_2$.
The following inner products hold:
\begin{equation}
  \footnotesize
  \label{eq_relations_G2}
  \begin{split}
  & (\alpha_1, \alpha_2) = 1\cdot\sqrt{3}\cos \frac{5\pi}{6} =  -\frac{3}{2}, \\
  & (\alpha_1 + \alpha_2, \alpha_2) = -\frac{3}{2} + 3 = \frac{3}{2}
        \; \Longrightarrow \;
        s_{\alpha_1 + \alpha_2}(\alpha_2) = -(3\alpha_1 + 2\alpha_2), \\
  & (\alpha_1 + \alpha_2, \alpha_1) = 1 -\frac{3}{2} = -\frac{1}{2}
        \; \Longrightarrow \; s_{\alpha_1 + \alpha_2}(\alpha_1) = 2\alpha_1 + \alpha_2, \\
  & (3\alpha_1 + \alpha_2, \alpha_1) = 3 - \frac{3}{2} = \frac{3}{2}
    \; \Longrightarrow \; s_{3\alpha_1 + \alpha_2}(\alpha_1) = -(2\alpha_1 + \alpha_2), \\
  & (3\alpha_1 + \alpha_2, \alpha_2) = -\frac{9}{2} + 3 = -\frac{3}{2}
    \; \Longrightarrow \; s_{3\alpha_1 + \alpha_2}(\alpha_2) = 3\alpha_1 + 2\alpha_2.
  \end{split}
\end{equation}
Then,
\begin{equation*}
 \footnotesize
  \begin{split}
  &  w_0(\alpha_1) = -s_{\alpha_1 + \alpha_2}(2\alpha_1 + \alpha_2) =
   -((4\alpha_1 + 2\alpha_2) - (3\alpha_1 + 2\alpha_2)) = -\alpha_1,  \\
  &  w_0(\alpha_2) =  s_{\alpha_1 + \alpha_2}(3\alpha_1 + 2\alpha_2) =
    3(2\alpha_1 + \alpha_2) - 2(3\alpha_1 + 2\alpha_2) = -\alpha_2,
  \end{split}
\end{equation*}
see \cite[Plate IX]{Bo02}. The longest element is unique in $W(G_2)$, so
$w_0$ is the longest. \qed

The root $3\alpha_1 + \alpha_2$ in the decomposition \eqref{eq_first_factoriz_G2} 
is not the highest root in any root subsystem $G_2$.
There is another decomposition of $w_0$ with $2$ factors,  
one of which corresponds to a simple root, and the second to the highest root.

\begin{proposition}[decomposition in $W(G_2)$]
  \label{prop_factoriz_G2}
  The longest element in $W(G_2)$ is decomposed into $2$ factors corresponding to
  the following \underline{orthogonal} roots:
\begin{equation}
  \label{eq_factoriz_G2}
    \alpha_1, \; 3\alpha_1 + 2\alpha_2.
\end{equation}
  Here, $3\alpha_1 + 2\alpha_2$  is the \underline{highest root}
  and $\alpha_1 < \alpha_{max}$.
\end{proposition}

\PerfProof
  By \eqref{eq_longest_G2} and \eqref{eq_first_factoriz_G2} we have
\begin{equation*}
  w_0 = (s_2s_{3\alpha_1 + \alpha_2}s_2)s_1.
\end{equation*}
Further, by \eqref{eq_relations_G2} we have  $(3\alpha_1 + \alpha_2, \alpha_2)  = -\frac{3}{2}$,
then
\begin{equation*}
  s_{\alpha_2}(3\alpha_1 + \alpha_2) =  3\alpha_1 + \alpha_2
      -\frac{2(-\frac{3}{2})}{3}\alpha_2 = 3\alpha_1 + 2\alpha_2.
\end{equation*}
Then, by \eqref{eq_ab_relations} we have
\begin{equation*}
  w_0 = s_{3\alpha_1 + 2\alpha_2}s_{\alpha_1}.
\end{equation*}
At last,  $(3\alpha_1 + 2\alpha_2, \alpha_1) =  3 + 2(-\frac{3}{2}) = 0$,
i.e., roots $3\alpha_1 + 2\alpha_2$ and $\alpha_1$ are orthogonal.
\qed
~\\

\begin{appendix}
\section{\bf The highest roots}

\begin{table}[H]
\centering
\renewcommand{\arraystretch}{1.4}
\begin{tabular}{|c|c|c|}
  \hline
       \footnotesize {Root}          & \footnotesize{The highest} & \footnotesize{The highest root in} \\
           %% \footnotesize{Reflection}   \\
       \footnotesize {system}        & \footnotesize{root} & \footnotesize{ the root subsystems }  \\
  \hline
       $A_n$  & \footnotesize{$\alpha_{max} = \alpha^{an}_{max} = $}
           & \footnotesize$\alpha^{ai}_{max}$ =  \\
              & \footnotesize{$\alpha_1 + \dots + \alpha_n$} &
                \footnotesize{$\alpha_1 + \dots + \alpha_i$} \\
  \hline
       $B_n$  & \footnotesize{$\alpha_{max} = \alpha^{bn}_{max} = $}  & 
          \footnotesize{$\alpha^{bi}_{max}$} = \\ 
              & \footnotesize{$\alpha_1 + 2\alpha_2 + \dots + 2\alpha_n$} &
                \footnotesize{$\alpha_{n-i+1} + 2\alpha_{n-i+2} + \dots + 2\alpha_n$} \\
  \hline
       $C_n$  & \footnotesize{$\alpha_{max} = \alpha^{cn}_{max} = $}  &
         $\footnotesize\alpha^{ci}_{max}$ = \\ 
              & \footnotesize{$2\alpha_1 + \dots + 2\alpha_{n-1} + \alpha_n$} &
                \footnotesize{$2\alpha_{n-i+1} + \dots$  + $2\alpha_{n-1} + \alpha_n$} \\
  \hline
       $D_n$  & \footnotesize{$\alpha_{max} = \alpha^{dn}_{max} = $} 
        & \footnotesize{$\alpha^{di}_{max}$} = \\
              & \footnotesize{$\alpha_1 + 2\alpha_2 + \dots + 2\alpha_{n-2} + $} &
                \footnotesize{$\alpha_{n-i+1} + 2\alpha_{n-i+2} + \dots  + 2\alpha_{n-2}$} \\
              & \qquad \footnotesize{$ \alpha_{n-1} + \alpha_n$} &
                \;\; \footnotesize{$  + \alpha_{n-1} + \alpha_n$}   \\
  \hline
\end{tabular}
\vspace {1mm}
\caption{\footnotesize{The highest roots in the root systems $A, B, C, D$ 
and their root subsystems, see Table \ref{tab_factoriz_all}}.}
\label{tab_highest_roots_ABCD}
\end{table}

\begin{table}[H]
\centering
\renewcommand{\arraystretch}{1.4}
\begin{tabular}{|c|c|c|}
  \hline
  \footnotesize {Root}          & \footnotesize{The highest} & \footnotesize{The highest roots in}  \\
       \footnotesize {system}        & \footnotesize{root} & \footnotesize{ the root susbsystems}    \\
  \hline
       $E_6$  & \footnotesize{$\alpha_{max} = \alpha^{e6}_{max}$} =  &
         \footnotesize{$\alpha^{a3}_{max}$}, \; \footnotesize{${\alpha^{a5}_{max}}$} \\
         & \footnotesize{$\alpha_1 + 2\alpha_2 + 2\alpha_3 + 3\alpha_4 + 2\alpha_5 + \alpha_6$} &   \\
  \hline
       $E_7$  & \footnotesize{$\alpha_{max} = \alpha^{e7}_{max}$} =   &
         \footnotesize{$\alpha^{d4}_{max}$}, \; \footnotesize{$\alpha^{d6}_{max}$}  \\
         & \footnotesize{$2\alpha_1 + 2\alpha_2 + 3\alpha_3 + 4\alpha_4 + 3\alpha_5 + 2\alpha_6 + \alpha_7$}   &  \\
  \hline
       $E_8$  & \footnotesize{$\alpha_{max} = \alpha^{e8}_{max}$}  =  &
         \footnotesize{$\alpha^{d4}_{max}$}, \;  $\alpha^{d6}_{max}$, \; $\alpha^{e7}_{max}$  \\
         & \footnotesize{$2\alpha_1 + 3\alpha_2 + 4\alpha_3 + 6\alpha_4 + 5\alpha_5 + 4\alpha_6 +
         3\alpha_7 + 2\alpha_8$}   &  \\
  \hline
       $F_4$  & \footnotesize{$\alpha_{max} = \alpha^{f4}_{max}$} =
          \footnotesize{$2\alpha_1 + 3\alpha_2 + 4\alpha_3 + 2\alpha_4$} &
         \footnotesize{$\alpha^{c2}_{max}$, \; $\alpha^{c3}_{max}$} \\
         & & \\
  \hline
       $G_2$  & \footnotesize{$\alpha_{max} = \alpha^{g2}_{max}$} =
            \footnotesize{$3\alpha_1 + 2\alpha_2$} &   \\
            & & \\
  \hline
\end{tabular}
\vspace {1mm}
\caption{\footnotesize{The highest roots in the root systems $E, F, G$ and 
 their root subsystems,  see Table \ref{tab_factoriz_all}}.}
\label{tab_highest_roots_EFG}
\end{table}

\section{\bf Two conjugation cases in the Weyl group}
  \label{sec_app_1}
  
  Let $\mathcal{E}$ be the linear space spanned by simple roots of $\Delta$, 
and $(\;,\;)$ be the symmetric bilinear form on $\mathcal{E}$ invariant under Weyl 
group $W(\varPhi)$, \cite[Ch.VI, $\S$1, Prop.7]{Bo02}.  

We distinguish two cases of reflections used in conjugation: those corresponding to long roots
(Proposition \ref{prop_fact_1})
and those corresponding to short roots (Proposition \ref{prop_fact_2}).

\begin{lemma}
For any two roots $\alpha$ and $\beta$ we have
~\\
\begin{equation}
 \label{eq_ab_relations}
 s_{\alpha}s_{\beta}s_{\alpha} = s_{s_{\alpha}(\beta)}.
\end{equation}
\end{lemma}
~\\
\PerfProof
Really,
\begin{equation*}
  \begin{split}
   & s_{\beta}s_{\alpha}(x) = s_{\alpha}(x) - \frac{2(\beta, s_{\alpha}(x))}{(\beta, \beta)}\beta =
      s_{\alpha}(x) - \frac{2(s_{\alpha}(\beta), x)}{(s_{\alpha}(\beta), s_{\alpha}(\beta))}\beta, \text{ and } \\
   & s_{\alpha}s_{\beta}s_{\alpha}(x) =
      x -  \frac{2(s_{\alpha}(\beta), x)}{(s_{\alpha}(\beta), s_{\alpha}(\beta))}s_{\alpha}(\beta) =
            s_{s_{\alpha}(\beta)}(x), \text{ i.e., } \\
   & s_{\alpha}s_{\beta}s_{\alpha} = s_{s_{\alpha}(\beta)}.
  \end{split}
\end{equation*}
\qed

For $\alpha, \beta \in \varPhi$, put
\begin{equation}
     n(\alpha, \beta) := \frac{2(\alpha, \beta)}{(\beta, \beta)}.
\end{equation}

Then, 
\begin{equation}
 \label{eq_action_sa}
   s_{\alpha}(\beta) = \beta - n(\beta, \alpha)\alpha.
\end{equation}

\begin{proposition}[case of a long root]
 \label{prop_fact_1}
  Let $\alpha$ and $\beta$ be 
  two roots in $\varPhi$, such that $\norm{\alpha} \geq \norm{\beta}$.
  If $n(\beta, \alpha) = \pm1$,
  the conjugation of $s_{\beta}$ by $s_{\alpha}$
  is as follows:
  \begin{equation}
    \label{eq1_fact_1}
       s_{\alpha}s_{\beta}s_{\alpha} =
      \begin{cases}
         s_{\alpha+\beta} \text{ if } n(\beta,\alpha) = -1, \\
         s_{\alpha-\beta} \text{ if } n(\beta,\alpha) = 1. \\
      \end{cases}
  \end{equation}
~\\
  The length of the root $\alpha+\beta$ (resp. $\alpha-\beta$) coincide with
  $\norm{\beta}$:
  \begin{equation*}
    \begin{split}
     & \norm{\alpha + \beta} = \norm{\beta} \;  \text{ if } n(\beta,\alpha) = -1,  \\
     & \norm{\alpha - \beta} = \norm{\beta} \; \text{ if } n(\beta,\alpha) = 1.
    \end{split}
  \end{equation*}
 In other words, the length $\norm{\alpha + \beta}$ (resp. $\norm{\alpha - \beta}$)
 either coincides with lengths of $\norm{\alpha}$ and $\norm{\beta}$ (when they are the same length)
 or coincides with length of the shortest of them. 
\end{proposition}  

\PerfProof
For all cases, we substitute $\alpha$ and $\beta$ in \eqref{eq_ab_relations}.
If $n(\beta, \alpha) = \pm{1}$ then 
 \begin{equation*}
 \label{eq_action_sa_proof_2}
 \begin{split}
   s_{\alpha}(\beta) = &
      \begin{cases}
        & \alpha + \beta \; \text{ if } \; n(\beta, \alpha) = -1, \\
        & \alpha - \beta \; \text{ if } \; n(\beta, \alpha) = 1.
      \end{cases} 
     \quad
      2(\beta, \alpha) =
      \begin{cases}
        & \;\;\, \norm{\alpha}^2  \; \text{ if } \; n(\beta, \alpha) = 1, \\
        & -\norm{\alpha}^2  \; \text{ if } \; n(\beta, \alpha) = -1. \\
      \end{cases} \\      
    \\
   & s_{\alpha}s_{\beta}s_{\alpha} = s_{s_{\alpha}(\beta)} = 
      \begin{cases}
        & s_{\alpha + \beta} \; \text{ if } \; n(\beta, \alpha) = -1, \\
        & s_{\alpha - \beta} \; \text{ if } \; n(\beta, \alpha) = 1.
      \end{cases}
  \end{split}
 \end{equation*}
 ~\\
 Then,
 \begin{equation*}
  \begin{array}{llll}
   & \text{ If } n(\alpha, \beta) = -1  & \text{then } &
      \norm{\alpha + \beta}^2 = \norm{\alpha}^2 + \norm{\beta}^2 + 2(\beta, \alpha) = \norm{\beta}^2,  \\
   & \text{ If } n(\alpha, \beta) =  1  & \text{then } &
      \norm{\alpha - \beta}^2 = \norm{\alpha}^2 + \norm{\beta}^2 - 2(\beta, \alpha) = \norm{\beta}^2.  
         \qed  
   \end{array}    
 \end{equation*}

\begin{remark} {\rm
 \label{rem_inner_prod}
  Let $\alpha$ and $\beta$ be roots in $\varPhi$ such that $n(\alpha, \beta) = \pm{1}$.
  All such cases are presented in \cite[Ch.VI,$\S1$,$n^{\rm o}3$]{Bo02}.
  
 {\rm{(i)}}
  If $\norm{\alpha} = \norm{\beta}$ then $(\beta, \alpha) = n(\alpha, \beta) = \pm{1}$.
  This holds for 
   \begin{itemize}
     \item two roots in $A_n$, $D_n$, $E_n$,     
     \item two long roots in $B_n$, $F_n$,    
     \item two short roots in $C_n$, $F_n$.
   \end{itemize} 

 {\rm{(ii)}}
  If $\norm{\alpha} = \sqrt{2}\norm{\beta}$ and $\norm{\beta} = 1$ 
  then $(\beta, \alpha) = n(\beta, \alpha) = \pm{1}$.
  This  holds for two roots of different lengths in $B_n$, $F_n$.
~\\
  
 {\rm{(iii)}} 
 If $\norm{\alpha} = \sqrt{2}\norm{\beta}$ and $\norm{\beta} = \sqrt{2}$
 then $(\beta, \alpha) = 2n(\beta, \alpha) = \pm{2}$.
 This  holds for the roots of different lengths in $C_n$.
~\\
  
 {\rm{(iv)}}
 Let  $\norm{\alpha} = \sqrt{3}\norm{\beta}$ and $\norm{\beta} = 1$.  
 Then $(\beta, \alpha) = \frac{3}{2}n(\beta, \alpha) = \pm\frac{3}{2}$.
 This  holds for the roots of different lengths in $G_2$.
 }
\end{remark}

\begin{proposition}[case of a short root]
 \label{prop_fact_2}
  Let $\alpha$ and $\beta$ be
  two roots in $\varPhi$, such that $\norm{\alpha} > \norm{\beta}$.

 {\rm{(i)}}
   If $\norm{\alpha} = \sqrt{2}\norm{\beta}$, $n(\alpha, \beta) = \pm{2}$,
   the conjugation of $s_{\alpha}$ by $s_{\beta}$
   is as follows:
  \begin{equation}
    \label{eq1_fact_2}
       s_{\beta}s_{\alpha}s_{\beta} =
      \begin{cases}
         s_{\alpha+2\beta} \text{ if } n(\alpha, \beta) = -2, \\
         s_{\alpha-2\beta} \text{ if } n(\alpha, \beta) = 2. \\
      \end{cases}
  \end{equation}

  The root $\alpha+2\beta$ (resp. $\alpha-2\beta$) in the case $n(\alpha, \beta) = -2$
  (resp. $n(\alpha, \beta) = 2$) is long.
  This case occurs for pairs of roots $\norm{\alpha} > \norm{\beta}$ in
  $B_n$, $C_n$ and $F_4$. 
~\\

 {\rm{(ii)}}
   If $\norm{\alpha} = \sqrt{3}\norm{\beta}$,  $n(\alpha, \beta) = \pm{3}$,
   the conjugation of $s_{\alpha}$ by $s_{\beta}$
   is as follows:
  \begin{equation}
    \label{eq1_fact_4}
       s_{\beta}s_{\alpha}s_{\beta} =
      \begin{cases}
         s_{\alpha+3\beta} \text{ if } (\alpha, \beta) = -\frac{3}{2}, \\
         s_{\alpha-3\beta} \text{ if } (\alpha, \beta) = \frac{3}{2}. \\
      \end{cases}
  \end{equation}
\end{proposition}

  This case occurs for pairs of roots $\norm{\alpha} > \norm{\beta}$ in $G_2$.
~\\

\PerfProof
(i)
Eq. \eqref{eq1_fact_2} follows from \eqref{eq_action_sa}:
\begin{equation*}
  s_{\beta}(\alpha) = \alpha - n(\alpha, \beta)\beta =
  \begin{cases}
     \alpha + 2\beta \; \text{ if } n(\alpha, \beta) = -2,\\
     \alpha - 2\beta \; \text{ if } n(\alpha, \beta) = 2.
  \end{cases}
\end{equation*}

 For cases $B_n$, $C_n$, $F_4$,
 since $(\beta, \beta) = 1$, we have $(\alpha, \beta) = \displaystyle\frac{n(\alpha, \beta)}{2} = \pm{1}$.
 By \eqref{eq1_fact_2} $\alpha+2\beta$ (similarly for $\alpha-2\beta$) is long, because
\begin{equation*}
 \footnotesize
\begin{split} 
  & (\alpha+2\beta, \alpha+2\beta) = (\alpha, \alpha) +
      4(\alpha,\beta) + 4(\beta, \beta) =  2 - 4 + 4 = 2 \; \text{ if } (\alpha, \beta) =  -1, \\
  & (\alpha-2\beta, \alpha-2\beta) = (\alpha, \alpha) 
      -4(\alpha,\beta) + 4(\beta, \beta) =  2 - 4 + 4 = 2 \; \text{ if } (\alpha, \beta) =  1.
\end{split}      
\end{equation*}

(ii) Here, $(\beta, \beta) = 1$ and
 $(\beta, \alpha) = \displaystyle\frac{n(\alpha, \beta)}{2} = \pm\frac{3}{2}$, then
\begin{equation*}
  s_{\beta}(\alpha) = \alpha - n(\beta, \alpha)\beta =
  \begin{cases}
     \alpha + 3\beta \; \text{ if } n(\alpha, \beta) = -3,\\
     \alpha - 3\beta \; \text{ if } n(\alpha, \beta) = 3.
  \end{cases}
\end{equation*}
 By \eqref{eq1_fact_4} the root of $\alpha+3\beta$ (similarly for $\alpha-3\beta$) is long, because
\begin{equation*}
\footnotesize
 \begin{split}
  & (\alpha+3\beta, \alpha+3\beta) = (\alpha, \alpha) +
      6(\alpha,\beta) + 9(\beta, \beta) =  3 - 6(\frac{3}{2}) + 9 = 3 \text { if } n(\alpha, \beta) = -3, \\
  & (\alpha-3\beta, \alpha-3\beta) = (\alpha, \alpha) -
      6(\alpha,\beta) + 9(\beta, \beta) =  3 - 6(\frac{3}{2}) + 9 = 3 \text { if } n(\alpha, \beta) = 3.
      \qed
 \end{split}     
\end{equation*}

\end{appendix}

\bibliographystyle{plain}
{\small\bibliography{Decompos_longest_elem}}
\end{document}